\documentclass{amsart}

\usepackage{amsmath}
\usepackage{amsthm}
\usepackage{amsfonts}
\usepackage{amssymb}

\newtheorem{theorem}{Theorem}[section]
\newtheorem{proposition}{Proposition}[section]
\newtheorem{lemma}[theorem]{Lemma}
\newtheorem{corollary}[theorem]{Corollary}
\newtheorem{definition}[theorem]{Definition}

\title[Stationary phase method for oscillatory RHPs]{A nonlinear stationary phase method for oscillatory Riemann-Hilbert problems}

\author
[Yen Do]
{Yen Do}

\address{Yen Do, Department of Mathematics,
UCLA, Los Angeles, CA 90095-1555, USA}
\email{qdo@math.ucla.edu}

\begin{document}

\begin{abstract}
We study the asymptotic behavior of oscillatory Riemann-Hilbert problems arising in the AKNS hierarchy of integrable nonlinear PDE's. Our method is based on the Deift-Zhou nonlinear steepest descent method in which the given Riemann-Hilbert problem localizes to small neighborhoods of stationary phase points. In their original work, Deift and Zhou only considered analytic phase functions. Subsequently Varzugin extended the Deift-Zhou method to a certain restricted class of non-analytic phase functions. In this paper, we extend Varzugin's method to a substantially more general class of non-analytic phase functions. In our work real variable methods play a key role.
\end{abstract}

\maketitle

\tableofcontents

\section{Introduction}

In many studies of asymptotical behaviors of nonlinear systems, the following Riemann-Hilbert problem (RHP) appears: given $J$ an invertible matrix-valued function on $\mathbb R$, we look for a matrix-valued function
$M(\lambda)$ which is analytic on $\mathbb C \setminus \mathbb R$ and satisfies the following jump condition:
\begin{eqnarray*}
    M_+(\lambda) &=& M_-(\lambda)J(\lambda), \; \; \text{for a.e. } \lambda \in \mathbb R
\end{eqnarray*}
$M_\pm(\lambda)$ denote the upper and lower parts of $M$, with appropriate normalization at $\infty$ and non-tangential limits on $\mathbb R$ (denoted by the same symbols). The jump matrix $J$ could be thought of as the multiplicative gain of $M$ when moving from the lower half plane to the upper part.

In this paper, we will be only interested in the $2 \times 2$ case since this appears in most of the applications. In the natural $L^2$ setting, the following normalization conditions are often imposed:
\begin{equation}\label{oscRHP}
M_+(\lambda) - I \in H_2(\mathbb C_+),  \ M_-(\lambda)-I \in H_2(\mathbb C_-),
\end{equation}
where $H_2(\mathbb C_+)$ is the complex Hardy space $H_2$ of the upper half plane and $H_2(\mathbb C_-)$ is the complex Hardy space $H_2$ of the lower half plane, and $I$ is the identity matrix. It is not hard to see that these conditions are equivalent to the existence of $h\in L^2(\mathbb R)$ such that $M_\pm = I + C_\pm(h)$. An overview of $L^p$ Riemann-Hilbert theory is available in \cite{deiftzhouNLSyd}.

For many equations in the AKNS hierarchy of integrable nonlinear PDEs, the jump matrix $J$ has the following oscillatory structure \cite{bealscoifmanyd}:
\begin{eqnarray}
\label{oscjumpmatrix}   J(\lambda,t) &:=& \begin{pmatrix}
                              &1+p(\lambda)q(\lambda) & p(\lambda)e^{-it\theta(\lambda)} \cr
                              & q(\lambda)e^{it\theta(\lambda)} & 1
                       \end{pmatrix}
\end{eqnarray}
Here $p(\lambda)$ and $q(\lambda)$ are the reflection coefficients of the initial data for the corresponding PDE (often $q = \pm\overline{p}$), and $p,q$ will be sufficiently smooth with sufficient decay if the initial data is sufficiently nice (see \cite{zhouL2Sobolevyd} for further details). The phase $\theta$ is real-valued and depends on the algebraic structure of the current PDE. The precise technical assumptions on $p,q,\theta$ used in this paper will be specified later. The oscillatory factors $e^{\pm it\theta}$ originate from the evolution of the reflection coefficients when we let the initial data evolves according to the PDE. Several examples include $\theta(\lambda) = (\lambda-\lambda_0)^2$ for NLS
and $\theta(\lambda) =  4(\lambda^3 - 3\lambda_0^2\lambda)$ for mKdV \cite{bealscoifmanyd}.

In the setting of (\ref{oscRHP}) and (\ref{oscjumpmatrix}), the question of interests is the long-time behavior of the \emph{potentials}
\begin{equation}\label{recoverlimit}
u(t) = \lim_{\lambda \to \infty} \lambda M_{12}(\lambda,t)
\end{equation}
(and similarly $v(t) = \lim_{\lambda \to \infty} \lambda M_{21}(\lambda,t)$).
Here in the limits $\lambda\to\infty$ non-tangentially in $\mathbb C$. In the AKNS setting, the function $u(t)$ recovered from (\ref{recoverlimit}) is a constant multiple of the solution $u(x,t)$ to the respective nonlinear PDE, therefore the above question corresponds to the long-time behavior $u(x,t)$. We note that the spatial parameter $x$ in $u(x,t)$ has been encoded in the phase $\theta$ and hence it is encoded in the stationary points $\{\lambda: \theta'(\lambda)=0\}$, and in this paper the stationary points will be assumed constant (see also the technical assumptions on $\theta$ and $p,q$ below). This is the reason why we will simply write $u(t)$ and $v(t)$ instead of $u(x,t)$ and $v(x,t)$. For simplicity of notation, we'll often suppress $t$ when writing $M$ in this paper.

Oscillatory RHPs also appear in other settings, in which the structure of the oscillatory jump matrix $J$ may be different from (\ref{oscjumpmatrix}) and the Riemann-Hilbert solution may no longer have the $L^2$ normalization (\ref{oscRHP}). For instance, the small-dispersion limit $\epsilon \to 0^+$ of the KdV equation
$$u_t - 6uu_x + \epsilon^2 u_{xxx} = 0$$
corresponds to the large $1/\epsilon$ asymptotics of an oscillatory RHP with a different $L^2$-normalization \cite{GGKMyd, Shabatyd}. Other examples include the (non $L^2$)Riemann-Hilbert formulation of Fokas, Its, Kitaev for orthogonal polynomials with varying weights $e^{-nV(x)}dx$ \cite{FokasItsKitaevyd} whose asymptotics plays an important role in random matrix theory (see \cite{DeiftetalUniformOPyd, McLaughlinMiller08yd} and the references therein, for related results involving random matrix theory and orthogonal polynomials see \cite{Lubinskyyd, LevinLubinskyyd} and the references therein).

The study of long-time behaviors of solutions to nonlinear integrable PDEs goes back at least to the work of Zakharov and Manakov \cite{zamayd}, who were probably the first to write down a correct formula for the leading asymptotics of $u(t)$ (and $v(t)$) in the NLS case. Using monodromy theory, Its \cite{itsyd} was able to reduce the RHP formulation for NLS to a model case, which can then be solved explicitly, giving the desired asymptotics for $u$ (and $v$).

Deift and Zhou \cite{deiftzhouMKdVyd} developed a rigorous nonlinear steepest descent method to study the oscillatory RHP associated with the mKdV equation. Their argument demands analyticity, at least for $\theta(\lambda)$ (the reflection coefficients can be approximated by suitable
analytic functions \cite{deiftzhouNLSyd}). For the mKdV equation there is more than one stationary point and Deift and Zhou were able to separate their contributions using an operator formulation for RHPs that goes back to Beals and Coifman \cite{bealscoifmanyd}.

The first work in the non-analytic setting is due to Varzugin \cite{varzuginyd} in the setting (\ref{oscRHP},\ref{oscjumpmatrix}) for $\theta$ with stationary points of first order. Varzugin's asymptotics was previously ansatzed by Kitaev \cite{kitaevyd} using a method of isomonodromic deformation.

The steepest descent method of Deift and Zhou has been extended in various ways.
In their study of the small-dispersion limit of the KdV equation, Deift, Venakides, and Zhou  \cite{DeiftZhouVenakidesyd} discovered that the contribution to the asymptotics of the solution $u(x,t,\epsilon)$ recovered from the solution of the corresponding RHP (via a limit similar to (\ref{recoverlimit})) comes from a family of intervals instead of the isolated stationary points of the phase in the original RHP. The key idea is to introduce a $g$-function such that after conjugating the given RHP with $e^{ig(\lambda)\sigma_3/\epsilon}$ where $\sigma_3 = \begin{pmatrix} 1 & 0 \cr 0 & -1\end{pmatrix}$ is the third Pauli matrix,
one arrives at an RHP whose jump matrix
$$J_{new} = e^{-i\sigma_3 g_-(\lambda)/\epsilon}J_{old} e^{i\sigma_3 g_+(\lambda)/\epsilon}$$
is \emph{ready} for a steepest descent argument (in Section~\ref{essencesect} we shall discuss this further). Furthermore, a systematic procedure for determining the intervals contributing to the asymptotics of the Riemann-Hilbert problem was also developed. The methods of Deift, Venakides, Zhou have been successfully applied to the orthogonal polynomial setting \cite{DeiftetalUniformOPyd}, which led to deep results in random matrix theory. The argument however still requires (real) analyticity of the phases of the new oscillatory jump matrices. In the dispersion KdV setting, this goes back to analyticity of the initial data and in the orthogonal polynomial setting it goes back to analyticity of the underlying weights.

The second extension of the steepest descent method of Deift and Zhou is the $\overline{\partial}$-steepest descent method of McLaughlin and Miller, which first appeared in the orthogonal polynomial setting \cite{McLaughlinMiller06yd, McLaughlinMiller08yd}. The $\overline{\partial}$ method follows the general scheme of the steepest descent argument of Deift-Venakides-Zhou, however non-analytic data are now continued to the desired contours via the solution of a $\overline{\partial}$ equation. In particular, Green's theorem (which works for smooth functions) is used as a remedy for Cauchy's theorem (which requires analyticity), and the error terms are under control if appropriate constraints are imposed on the $\overline{\partial}$ problem. The $\overline{\partial}$-steepest descent method of McLaughlin and Miller allows for non-analytic phases with two Lipschitz derivatives near stationary points. In \cite{DiengMcLaughlinyd}, this method was adapted to the AKNS setting and was used to obtain long-time asymptotics for solutions to the defocusing NLS, with sharp error bound for initial data in $H^{1,1}(\mathbb R) = \{f\in L^2: f', xf \in L^2\}$.

The main goal of this paper is to further demonstrate that the asymptotics of oscillatory RHPs can be studied using analogue of real-variable tools from the linear theory of oscillatory integrals. We shall focus on the $L^2$-normalized RHPs (\ref{oscRHP}) with oscillatory jump matrix of the form (\ref{oscjumpmatrix}) (which corresponds to the study of long-time behavior of solutions to many nonlinear integrable PDEs), however the author anticipates applications of the argument developed in this paper to other settings (especially those mentioned above where nonanalytic data are frequently encountered). Due to the length and the scope of the current paper, the extent of the applications of the argument is not fully explored, but we shall discuss the underlying ideas that can be used for these potential adaptations in Section~\ref{essencesect}. The analysis used in this paper extends the work of Varzugin \cite{varzuginyd}, and for the RHP (\ref{oscRHP}, \ref{oscjumpmatrix}) our argument allows for non-analytic phases with a finite number of stationary points (of arbitrary orders), under fairly minimal regularity assumptions.

To see a simple connection between the setting of (\ref{oscRHP}, \ref{oscjumpmatrix}) and the linear theory of oscillatory integrals, consider a degenerate situation when $q(\lambda) = 0$. In this case, the oscillatory RHP (\ref{oscRHP}, \ref{oscjumpmatrix}) becomes {\it abelian} and can be solved using the Hilbert transform. Recovering $u(t)$ via the limit (\ref{recoverlimit}), we arrive at the linear problem: study the asymptotics behavior of $\int p(\lambda) e^{-it\theta(\lambda)}d\lambda$, and this could be achieved by the classical stationary phase method \cite{steinyd}, which exploits cancellation resulted from rapid oscillation of $e^{-it\theta(\lambda)}$ away from the stationary points of $\theta$. For this reason, the presentation of the proof of the main result of this paper will follow the spirit of the linear stationary phase method and is summarized in Section~\ref{outlinesect}.

In this paper, the following assumptions will be made on $\theta$:

(A) $\theta$ is real valued and continuously differentiable, and on the complement of a finite set of points $\theta$ has three locally integrable derivatives.

(B) $\theta$ has stationary points $\lambda_1, \dots, \lambda_N$ of orders $k_1, \dots, k_N$ i.e. $\theta$ is $(k_j+1)$-time differentiable near $\lambda_j$ and $\theta'(\lambda_j) =\dots =\theta^{(k_j)}(\lambda_j)=0\neq \theta^{(k_j+1)}(\lambda_j).$

(C) If $\lambda$ is a stationary point of order $k$ of $\theta$ then we require $\theta^{(k+1)}(x)$ to be H\"older continuous at $\lambda$ with some exponent $\beta>0$, i.e. near $\lambda$
$$\theta^{(k+1)}(x) = \theta^{(k+1)}(\lambda) + O(|x-\lambda|^{\beta}).$$
If $k=1$ then we assume $\theta^{(3)}$ is integrable near $\lambda$  in some high $L^r$, more precisely $r> \frac{1}{\beta}$ (by H\"older's inequality $\beta \geq 1/r'$, so this is automatic if $r>2$). If $k=2$ and $p(\lambda)q(\lambda)<0$ then we assume $\theta^{(3)}$ is Lipschitz near $\lambda$.

For convenient, we say that a stationary point $\lambda$ is focusing if $p(\lambda)q(\lambda)\geq 0$ and defocusing if $p(\lambda)q(\lambda)<0$. We note that (C) is weaker than having two Lipschitz derivatives near primary stationary points, which by Rademacher's theorem is equivalent to having a bounded $\theta^{(3)}$ there. The Lipschitz assumption of $\theta^{(3)}$ near defocusing secondary stationary points can also be weakened to a similar high $L^r$ integrability condition of $\theta^{(4)}$, and is chosen here for simplicity of the argument.

For $k,j\in \mathbb{Z}_+$ we recall the Sobolev space $H^{k,j}:= \{f\in L^2: f', \dots, f^{(k)} \in L^2, xf, \dots, x^j f\in L^2\}$. 

The general assumption for $p,q$ will be:

(D) $p,q\in H^{1,0}(\mathbb R)$ such that $0< 1+p(x)q(x) = O(1)$, and $p,q$ have sufficiently decay (which depends on $\theta$).

(E) $p,q$ have two $L^2$ derivatives near every stationary point $\lambda$, and if its order $k \geq 3$ then $p(\lambda)q(\lambda)<1$. If $k\geq 3$ or if $k=2$ and $p(\lambda)q(\lambda)<0$ then we require a third $L^2$ derivative for $p,q$ near $\lambda$, and the assumption in (D) will be changed to $p,q \in H^{2,0}$ overall.

The regularity assumptions on $p,q$ can be improved. For example, if $N=1$ and $\theta^{(k_1+1)}$ is constant near $\lambda_1$ (such as in the case of the NLS) then we only need one $L^2$ derivative for each $p$ and $q$ if $\lambda_1$ is a primary stationary point or a focusing secondary stationary point. For a stationary point $\lambda$ of order $k$, the requirement $p(\lambda)q(\lambda)<1$ may be removed if a certain model RHP associated with stationary points of order $k$ is better understood; the author plans to revisit this issue in a future manuscript.

In this paper we will not try to optimize the decay requirements for $p$ and $q$, this however could be done for explicit $\theta$.

For convenience of notation, denote
\begin{eqnarray*}
\nu_j = -\frac{1}{2\pi}\ln[1+p(\lambda_j)q(\lambda_j)],\;\;\; \epsilon_j =  \begin{cases}
         0, &\text{$k_j$ even;}\\
         \text{sgn}(\theta^{(k_j+1)}(\lambda_j)) , &\text{$k_j$ odd.}
       \end{cases}
\end{eqnarray*}

\begin{theorem}
\label{maintheorem}
$\exists t_0$ such that the RHP (\ref{oscRHP}, \ref{oscjumpmatrix}) has unique solution for $t\geq t_0$. Furthermore, as $t\to\infty$:

(i) If $N=0$ then the recovered potentials $u(t)$ and $v(t)$ satisfies
$$u(t), v(t) = O_\epsilon(t^{-1+\epsilon}), \;\;\;\epsilon >0$$
This can be arbitrarily improved if stronger regularity assumption on $p$, $q$, $\theta$ are given;

(ii) If $N\geq 1$ then there are $d_j>0$ such that
\begin{eqnarray*}
\begin{pmatrix}0 & u(t)\cr v(t) & 0 \end{pmatrix}
= \sum_{j=1}^N \Big[\begin{pmatrix}0 & u_j(t)\cr v_j(t) & 0 \end{pmatrix} + O(t^{-(\frac{1}{k_j+1} + d_j)})\Big]
\end{eqnarray*}
The contributions of $\lambda_j$ are of the form:
\begin{eqnarray}
\label{u_j}   u_j(t) &=& U_j p(\lambda_j)t^{-\frac{1}{k_j+1}}\exp\Big(-i\big[t\theta(\lambda_j)  + \alpha_j\ln t - 2\omega_j\big]\Big)\\
\label{v_j}    v_j(t) &=& V_j q(\lambda_j)t^{-\frac{1}{k_j+1}}\exp\Big(i\big[t\theta(\lambda_j)  + \alpha_j\ln t - 2\omega_j\big]\Big)
\end{eqnarray}
here $U_j,V_j$ depend on $\lambda_j$, $k_j$, $\theta^{(k_j+1)}(\lambda_j)$ and $p(\lambda_j)q(\lambda_j)$, and $U_j,V_j$ are anti-complex conjugates, i.e.
$$U_j = -\overline{V_j}$$
while $\alpha_j,\omega_j$ are real numbers and can be explicitly computed by
\begin{eqnarray*}
  \alpha_j &=& \frac{2\epsilon_j\nu_j}{k_j+1}\\
  \omega_j &=& \frac{1}{2\pi} \int_{D_-} \ln|\lambda_j-y|d\ln[1+pq](y)+\sum_{1\leq k \leq N, \lambda_k \neq \lambda_j} \epsilon_k\nu_k\ln|\lambda_j-\lambda_k|
\end{eqnarray*}
(the right-hand side for $\omega_j$ makes sense if $\ln(1+pq)$ has sufficient decay, say $H^{1,1}$, but $\omega_j$ still exists under weaker assumption, see Section~\ref{factorizesection} and in particular (\ref{omegajdefinition})).
\end{theorem}

\noindent \emph{Notes:} 1. Theorem~\ref{maintheorem} is related to a question in \cite{deiftopenprobyd}. The unique solvability for large $t$ of our oscillatory RHP can be proved for a larger class of $p$, $q$ (see Theorem~\ref{uniquesolvability}). A similar form of the above asymptotics when $\theta$ has one primary and one secondary phase point was previously given as an ansatz by Kitaev \cite{kitaevyd}.

2. Compared to linear theory, there is an extra $\ln t$ term in the exponent of the leading asymptotics of $u(t)$ and $v(t)$. There are also interactions between stationary points, which can be observed in the definition of $\omega_j$.

3. The error estimate $O(t^{-(\frac{1}{k_j+1} + d_j)})$ can be made explicit. For instance, if $r=\infty$ and $\beta_1=\dots = \beta_N =1$ then we can take $d_j = \frac{1}{2(k_j+1)}-\epsilon$, furthermore if $N=1$ then this can be improved to $\frac{1}{k_j+1}-\epsilon$. In the NLS case, these estimates are well-known from the work of Deift and Zhou \cite{deiftzhouNLSyd}.

4. If $k_j=1$ then $U_j$ (hence $V_j$) can be computed explicitly (see for instance \cite{varzuginyd}):
$$U_j = \frac{i\nu_j}{\sqrt{|\nu_j\theta''(\lambda_j)t|}} \frac{1}{\sqrt{|p(\lambda_j)q(\lambda_j)|}}\exp\Big(\frac{i\pi \epsilon_j}{4} + i\text{ arg }\Gamma(i\epsilon_j\nu_j)\Big)$$
(when $p(\lambda_j)=0$ or $q(\lambda_j)=0$ this should be understood in the limiting sense).

\subsection{Notational conventions} In this paper, we make the following conventions:

1. Absolute value of a matrix: For any matrix $M$ define $|M| = (\text{tr} M^*M)^{1/2}$ and for any matrix valued function $A$ on $\mathbb R$ define $\|A\|_p := \big\| |A| \big\|_p$. This makes the respective $L^p$ space a Banach space for $1<p<\infty$. We'll abuse notation and refer to this space also as $L^p$ - it should be clear from the context what is being referred to. Notice that the above absolute value $|.|$ satisfies the triangle inequality, furthermore $|AB| \lesssim |A||B|$ for any $A,B$.

2. Action of operators: The action of any operator T on matrices is done
entry-wise. If $T$ is a bounded operator on $L^p$ then it is also bounded as an operator on the space of $2\times 2$ matrices having $L^p$ entries, with comparable norm.

3. Differentiability: $f$ is said to be $k$-time differentiable if $f$ is $(k-1)$-time continuously differentiable and $f^{(k)}$ is locally integrable.

4. Inequalities up to a constant: For two quantities $A$ and $B$ we say $A \lesssim B$ if there exists an absolute constant $C > 0$ such that $|A| \le CB$. If the constant $C$ depends on the parameters $p_1, \dots, p_n$ we shall say that $A \lesssim_{p_1,\dots, p_n} B$. In some situations, when the dependence of $C$ on certain parameters are not important for the proof or the discussion, we shall abuse notation and suppress those parameters in writing, in particular if the dependence of $C$ on all the parameters are not important we will simply write $A \lesssim B$. The implicit constants used in the manuscript are not necessarily the same and may differ from line to line. The meaning and conventions used for $\gtrsim$ and $\gtrsim_{p_1,\dots, p_n}$ are understood similarly.

5. $\infty-$, $p+$, $p-$: For the sake of brevity of the argument, we shall often write $\infty-$ to denote a finite number that is very large. For any finite number $p$ we shall use $p+$ to denote a finite number $q>p$ such that $q-p$ is sufficiently small. Similarly, $p-$ shall denote a finite number $q<p$ such that $p-q$ is sufficiently small. This notation will be useful when applying H\"older inequality, for instance we can write
$\|fg\|_2 \leq \|f\|_{2+} \|g\|_{\infty-}$.

6. Sobolev spaces: For $k,j\in \mathbb{Z}_+$ we shall denote $H^{k,j}:= \{f\in L^2: f', \dots, f^{(k)} \in L^2, xf, \dots, x^j f\in L^2\}$.

\section{Outline of the proof of Theorem~\ref{maintheorem}}\label{outlinesect}

In this section, we describe the main ideas used to prove Theorem~\ref{maintheorem}.

One fundamental tool that will be used throughout is the operator formulation by Beals and Coifman for solutions to RHPs. Given a pair of $L^\infty$ weights $(w^-, w^+)$, Beals and Coifman introduced the following operator \cite{bealscoifmanyd}
$$C_w f = C_+(fw^-) + C_-(fw^+)$$
here $C_+ \equiv H_+$, $C_- \equiv - H_-$ and $H_+, H_-$ are respectively the projections into the complex Hardy space $H_p(\mathbb C_+), H_p(\mathbb C_-)$ of the upper and lower half planes. This operator can be thought of as a weighted version of the identity operator in the Hardy decomposition of $L^2$
$$f = C_+(f) - C_-(f)$$
Usually, the pair of weights $(w^-, w^+)$ is obtained from a factorization of the jump matrix
$$J = (I-w^-)^{-1}(I+w^+)$$
Except for model cases that will be discussed below, the weights in this paper are in $L^2$ (in addition to being in $L^\infty$). In that case, if $1-C_w$ is invertible on $L^2$ then the respective RHP will be uniquely solvable and its solution can be formulated using $C_w$ \cite{bealscoifmanyd} (see Section~\ref{opsection}). Beals-Coifman's operator formulation is useful in studying perturbation and unique solvability of RHPs.

Our proof of Theorem~\ref{maintheorem} will be a sequence of reductions following the spirit of the classical stationary phase method. In each reduction, we reduce one RHP to another RHP such that, heuristically, if Theorem~\ref{maintheorem} is true for the new RHP then it is true for the current RHP. More precisely, we will show that, for $t$ large, unique solvability of the new RHP implies that of the current RHP, and the respective potentials $u(t),v(t)$ differ by terms of highly decaying order of $t$. Note that for some of the intermediate RHPs that arise during our reductions, the defining limits for $u,v$ may not exist (most likely the limits in the upper and lower half-planes don't agree). Denoting by $\Delta M$ the difference between the solutions of the old and new RHPs, we overcome this issue by showing that the nontangential limits $\limsup_{z\to\infty} \Big|z\Delta M_{12}(z)\Big|$ and $\limsup_{z\to\infty} \Big|z\Delta M_{21}(z)\Big|$ decay strongly as $t\to\infty$. Since $\limsup$ is sub-additive, these estimates are sufficient.

We will start from the given oscillatory RHP whose jump matrix is defined in (\ref{oscjumpmatrix}) and finally end up at $N$ model cases. These model cases are similar and a general study of model cases will be done in Section~\ref{modelsection}. In particular, there we'll show the conclusions of Theorem~\ref{maintheorem} for each of the model cases.

In many reductions, we will use a variant of Varzugin's argument \cite{varzuginyd}. Our variant will be summarized in two perturbation schemes described in Section~\ref{localizationsection}. These schemes are useful when the weights are in $L^2$ (in addition to being in $L^\infty$), and they are based on the Beals-Coifman operator formulation. Several technical estimates will be needed to make these schemes possible; they will be discussed in Section~\ref{hardyspacesection}.

For technical reasons, a certain phase-weight relation in $(w^-,w^+)$ is essential to our reductions, and a conjugation by a scalar RHP can be used to achieve this relation. This type of RHP was introduced by Deift and Zhou \cite{deiftzhouMKdVyd} for a slightly different purpose, and was also used in \cite{varzuginyd}. If $M$ solves our initial RHP, $\sigma_3$ is the third Pauli matrix, and $\delta$ is the solution to the above scalar RHP then the corresponding RHP for $M^\delta:=M \delta^{-\sigma_3}$ will have the conjugated jump matrix $J_{conj}=\delta_-^{\sigma_3}J \delta_+^{-\sigma_3}$. A good choice of $\delta$ will ensure the existence of a factorization of $J_{conj}$ that gives a pair of weights $(w^-,w^+)$ with the desired phase-weight relation.

The first reduction will be {\it localization to small neighborhoods of stationary points}. We say a RHP with a pair of weight $(w^-, w^+)$ is localized to a set if these weights are supported on this set. After this reduction, the weights will be supported on a small neighborhood of $\lambda_1, \dots,\lambda_N$ (if $\theta$ doesn't have any stationary points, we can directly reduce them to $(0,0)$). In this reduction the respective potentials $u(t),v(t)$ differ by terms of high decaying order.

The next reduction will be {\it phase reduction}. The given RHP now has already been localized to a small neighborhood of the stationary points. In each neighborhood, the phase $\theta$ can be approximated by an analytic function (so overall we have a locally analytic function). Near a stationary point $\lambda$ of order $k$, we'll approximate $\theta$ by the following Taylor approximation:
$$\theta(\lambda) + \frac{\theta^{(k+1)}(\lambda)}{(k+1)!}(x-\lambda)^{k+1}$$

The third reduction will be {\it separation of contributions}. To separate the contribution of a stationary point $\lambda_0$, consider the following setting: Let $w_0$ be the part of $w$ that is supported near $\lambda_0$, and $w_1=w-w_0$ the rest. Each pair of weights $w_i$ gives rise to a new normalized $L^2$-RHP. We'll show that unique solvability of the new RHPs imply that of the current RHP, and the current potentials $u(t)$, $v(t)$ differ the sum of the respective potentials $u_0(t) + u_1(t), v_0(t) +v_1(t)$ by terms of high decaying order.

The separation of contributions is nontrivial in general, as matrix multiplication is noncommutative. The main idea is to handle commutators that arise when commuting terms during the calculation. To separate the contribution of $\lambda_0$, we'll use an a priori estimate involving the behavior of the solution to RHPs localized to small neighborhood of $\lambda_0$. The a priori estimate will be proved using the same reductions, i.e. starting from now we will have to show that our reductions will not harm the validity of these estimates.

For each stationary point, we'll generally carry out phase reduction before separating its contribution. The reason is technical: regarding regularity, the reduction of above-mentioned a priori estimate is more expensive than the reduction of $u,v$. Because two different stationary points might require two different regularity assumptions, we will iterate the phase reduction and contribution separation steps through the list of stationary points.

The final reduction will be a {\it reduction to model case}. To reduce a localized RHP with a nice phase $\Theta(x) = a+bx^{k+1}$ to an appropriate model RHP, we will use a steepest descent argument that goes back to Deift and Zhou \cite{deiftzhouMKdVyd}. The model RHP will not have $L^2$ weights (the weights however remain in $L^\infty$). This lack of $L^2$ integrability prevents direct applications of Beals-Coifman operator formulation, and hence our perturbation schemes.

To overcome the above difficulty, the main idea is to reduce the current RHP to an approximation of the model RHP, which we'll refer to as the pre-model RHP, which has $L^2 \cap L^\infty$ analytic weights. The sub-reduction to the pre-model RHP can be done in the same way as before. To reduce the pre-model RHP to the model RHP, we will follow \cite{deiftzhouNLSyd} and exploit analytic continuation to deform the contour $\mathbb R$ of our RHPs to the complex plane, such that on the new contour, the pre-model and model weights have strong decay. This means $L^2$ integrability of the weights will be available and hence we can use our favorite Beals-Coifman operator formulation.

The second difficulty in our reduction to the model RHP is the fact that the model weights are no longer localized to a small neighborhood of the current stationary point. This means the same argument that was used to reduce the a priori estimates needed in our separation argument will not be directly applicable. The main idea is to exploit analyticity of the model weights (and the pre-model weights) to make up for this lack of localization.

Finally, Theorem~\ref{maintheorem} will be verified for the explicit model RHPs that arise as consequences of the above reductions. For each stationary point, we have one model RHP. Beals-Coifman techniques can be used in the classical settings, i.e. when $-1<pq<1$ at the current stationary point. The quadratic non-classical case can be studied using confluent hypergeometric functions as in \cite{varzuginyd}, and for the cubic non-classical case we'll use a result from isomonodromy deformation theory\cite{ItsNovokshenovNotesyd}.

In every reduction, it will be clear that if the new (i.e. reduced) RHP is uniquely solvable, then so is the old RHP. Thus, the unique solvability of the model cases will automatically imply the unique solvability of the given oscillatory RHP. Indeed, we'll prove a norm estimate which implies the unique solvability of RHPs. The proof of this estimate uses the same reductions in the proof of Theorem~\ref{maintheorem}.

\section{Beals-Coifman operator formulation for RHPs} \label{opsection}
Given two weight functions $w^+, w^- \in L^\infty$, consider the following bounded Beals-Coifman operator (acting on $L^p, 1<p<\infty$):
\begin{eqnarray*}
\label{bcoperator} C_w f = C_+(fw^-) + C_-(fw^+)
\end{eqnarray*}
Here $C_+ \equiv H_+$, $C_- \equiv - H_-$ and $H_+, H_-$ are respectively the projections into the complex Hardy spaces $H_p(\mathbb C_+)$, $H_p(\mathbb C_-)$ of the upper and lower half planes. $C_w$ can take $L^\infty$ input if the weights are also in $L^p$.

Usually, we'll be interested in $L^2$ setting where the weights are obtained from a nice factorization of $J$:
$$J(\lambda) = (I-w^-)^{-1}(I+w^+), \  w^\pm \in L^2 \cap L^\infty$$
The $L^2\cap L^\infty$ condition will always be assumed on the weights in this paper, unless otherwise specified.

\begin{theorem} [Beals-Coifman] Suppose that $w^\pm \in L^2 \cap L^\infty$ are obtained from a factorization of $J$. If $\exists \mu \in I + L^2$ satisfying $\mu  = I + C_w \mu$
then the $L^2$-RHP has a solution:
$$M_\pm(\lambda) = I + C_\pm(\mu (w^+ + w^-))(\lambda)$$
\end{theorem}
\proof Notice that $M_+ = \mu(I + w^+)$ and $M_- = \mu(I-w^-)$.
\endproof

This theorem only gives us existence. Notice that the functional equation of $\mu$ can be rewritten as:
\begin{eqnarray}
\label{basiceqn} (I-C_w)(\mu-I) = C_w I
\end{eqnarray}

Thus, if $(1-C_w)$ is invertible on $L^2$ then (\ref{basiceqn}) always has a solution on $L^2$ (it is not hard to see that the invertibility of $C_w$ is independent of the choice of the factorization).

Under mild assumption on the factorizations, the invertibility of $(1-C_w)$ indeed implies the uniqueness of the solution constructed in the theorem of Beals and Coifman. This was observed by Deift and Zhou \cite{deiftzhouNLSyd} under the additional assumption that $(I+w^+)^{\pm 1}, (I-w^-)^{\pm 1} \in L^\infty$.

The following identity (or its variant) will be used frequently in future computation:  If $h,g,-\frac{1}{2}h(Hg)-\frac{1}{2}g(Hh) \in \bigcup_{1<p<\infty} L^p$ then
\begin{eqnarray}
\label{multiplycauchyop}    C(h)C(g) &=& C\big(-\frac{1}{2}h(Hg)-\frac{1}{2}g(Hh)\big)
\end{eqnarray}
($H$ denotes the Hilbert transform). This identity is a consequence of Privalov's theorem (see for instance Garnett \cite{garnettyd}).

An immediate consequence of (\ref{multiplycauchyop}) is that (under the same assumption on $g,h$) $C_\pm(h)C_\pm(g)$ respectively belong to $H^p(\mathbb C_\pm)$, so vanish under $C_\mp$ thanks to the orthogonality of $C_+,C_-$. Note that this orthogonality relies on the fact that $\mathbb R$ with the obvious orientation is a {\it complete} contour \cite{zhouyd}, see Section~\ref{reductmodelsection2} for more details. Essentially, this means this consequence of (\ref{multiplycauchyop}) remains true if $\mathbb R$ is replaced by any complete contour $\Sigma$.

As an example, we'll show the following important lemma in \cite{deiftzhouMKdVyd} (which was attributed to Zhou \cite{zhouyd}). Let $1<p<\infty$ and let $C_\Phi, C_{\widetilde \Phi}$ respectively denote $L^p$ Beals-Coifman operators with $L^\infty$ weights $(-\Phi_-,\Phi_+), (-\widetilde\Phi_-,\widetilde\Phi_+)$ such that
$$(I+\Phi_\pm)(I+\widetilde \Phi_\pm) = I$$
Let $w,w_\Phi$ be two pairs of $L^\infty$ weights such that:
$$(I+w^+_\Phi) = (I+w^+)(I+\Phi_+)$$
$$(I-w^-_\Phi) = (I-w^-)(I+\Phi_-)$$

\begin{lemma}\label{factorize} If for some $1< p_1,p_2<\infty$ we respectively have $\Phi_\pm \in H^{p_1}(\mathbb C_\pm)$ and $\widetilde \Phi_\pm \in H^{p_2}(\mathbb C_\pm)$, then the following identities are true on $L^p$:
\begin{eqnarray}
\label{factorize1}    1-C_{w_\Phi} &=& (1-C_\Phi)\circ(1-C_w)\\
\label{factorize2}           1-C_w &=& (1-C_{\widetilde\Phi})\circ(1-C_{w_\Phi})\\
\label{factorize3} (1-C_\Phi)\circ(1-C_{\widetilde\Phi}) &=& 1 \;\;=\;\; (1-C_{\widetilde\Phi})\circ(1-C_\Phi)
\end{eqnarray}
\end{lemma}
Below, we'll include a proof of Lemma~\ref{factorize} for future reference (which relies on orthogonality of $C_+$, $C_-$). As we'll see, $L^2$ integrability of $w$ (or $w_\Phi$) won't be needed.

\proof The proof of (\ref{factorize1}) relies on the following identities
\begin{eqnarray}
\label{multiplyCauchyPhi1} (C_\pm h)(x)\Phi_\pm(x)  = C_\pm \Big((C_+h) \Phi_+ - (C_-h) \Phi_-\Big)(x)
\end{eqnarray}
for every $h\in L^p(\mathbb R)$, any $1<p<\infty$. Indeed, using boundedness of $\Phi_\pm$ we can always write $\Phi_\pm = C_\pm(g)$ for some $g\in L^{p_1}$ with $g,Hg\in L^\infty$. Now, using (\ref{multiplycauchyop}) it is not hard to see that $(C_+ h) \Phi_+ \in H_p(\mathbb C_+)$, $(C_- h) \Phi_- \in H_p(\mathbb C_-)$, which easily implies (\ref{multiplyCauchyPhi1}).

The rest of the proof of (\ref{factorize1}) is simply algebraic. Notice that (\ref{factorize1}) can be rewritten as $C_\Phi \circ C_w =C_\Phi + C_w - C_{w_\Phi}$. Using the definition of $w_\Phi^\pm$, it is not hard to see that the right-hand side of this equation is a Beals-Coifman operator with the pair of weights
$$(-w^-\Phi_-, -w^+\Phi_+)$$
So essentially we need to understand $C_\Phi \circ C_w$. Fix any $f\in L^p(\mathbb R)$. Observe that any Beals-Coifman operator $C_w$ can be rewritten in two ways
\begin{eqnarray*}
    C_w(f) &=& fw^- + C_-(f(w^-+w^+))\\
                  &=& -fw^+ + C_+(f(w^-+w^+)),
\end{eqnarray*}
note that $f(w^+ + w^-)\in L^p$. Consequently, (\ref{multiplyCauchyPhi1}) and orthogonality of $C_\pm$ imply
\begin{eqnarray*}
    C_+(C_w(f)(-\Phi_-)) &=& -C_+(fw^-\Phi_-)\\
    C_-(C_w(f)\Phi_+) &=& -C_-(fw^+\Phi_+)
\end{eqnarray*}
which completes the proof of (\ref{factorize1}).

The proof of (\ref{factorize2}) is exactly the same with suitable changes of symbols.

Finally, the first equality in (\ref{factorize3}) can be viewed as a special case of (\ref{factorize1}) when the weights $w^\pm$ are the same as the weights $\pm \widetilde\Phi_\pm$ (i.e. $w^\pm_\Phi \equiv 0$). The second equality in (\ref{factorize3}) follows by symmetry.
\endproof
In some cases, $\Phi_\pm$ and $\widetilde \Phi_\pm$ may not belong to the above type of Hardy spaces, but they are still analytic and uniformly bounded on the respective half planes $\mathbb C_\pm$. In that case, the following lemma will be useful:

\begin{lemma} \label{boundedfactorize} The conclusion of Lemma~\ref{factorize} remains true if we assume $\Phi_\pm$ and $\widetilde \Phi_\pm$ are analytic and uniformly bounded on the respective half planes $\mathbb C_\pm$.
\end{lemma}
\proof The main idea is to directly prove (\ref{multiplyCauchyPhi1}) instead of using (\ref{multiplycauchyop}) (which requires the assumptions in Lemma~\ref{factorize} that we want to avoid). The rest of the argument is exactly the same.
To show (\ref{multiplyCauchyPhi1}), essentially we want to show that $(C_\pm h)\Phi_\pm \in H^p(\mathbb C_\pm)$ for $h\in L^p$. Fix any $z\in\mathbb C_+$. We'll show by contour integration that
\begin{eqnarray}
\label{ChPhi} (C_+ h)(z)\Phi_+(z)  = C\Big((C_+h) \Phi_+\Big)(z)
\end{eqnarray}
Choose the contour $\gamma$ to be a semicircle of radius $R$ centered at the origin lying in the upper half plane. Then for large $R$,
$$(C_+ h)(z)\Phi_+(z) = \frac{1}{2\pi i}\int_{-R}^R \frac{(C_+ h)(x)\Phi_+(x)}{x-z}dx + \text{contribution of $\gamma_R$}$$
where $\gamma_R:=\{Re^{i\beta}:0\leq \beta\leq \pi\}$. Since $p<\infty$, to show that the contribution of $\gamma_R$ vanishes as $R\to \infty$ it suffices to show that
$$\|(C_+ h)(z)\Phi_+(z)\|_{L^p(\gamma_R)} \lesssim \|h\|_{L^p(\mathbb R)}$$
This is a consequence of the boundedness of $\Phi_+$ in the upper half plane and Carleson's measure theorem (see for instance \cite{garnettbafyd}), which gives, for every $1<p<\infty$
$$\|(C_+ h)(z)\|_{L^p(\gamma_R)} \lesssim \|h\|_{L^p(\mathbb R)}$$
Consequently, (\ref{ChPhi}) is proved for $z\in \mathbb C_+$. Sending $z\to \mathbb R$ nontangentially, we see that $(C_+ h)\Phi_+ \in H^p(\mathbb C_+)$ (observe that the boundedness of $\Phi_+$ on $\mathbb R$ ensures that $(C_+h) \Phi_+ \in L^p$).

Similarly, $(C_- h)\Phi_- \in H^p(\mathbb C_-)$, and (\ref{multiplyCauchyPhi1}) can now be deduced easily.
\endproof
\noindent {\it Remarks:} Lemma~\ref{factorize} and \ref{boundedfactorize} remain true when $\mathbb R$ is replaced by a {\it complete} contour, and this fact is needed for contour deformation of model RHPs in Section~\ref{reductmodelsection2}. Using (\ref{factorize3}), we can show a partial converse to the above observation of Deift and Zhou. The following proposition is useful when it is possible to construct a nice solution to (\ref{basiceqn}). Indeed, we are able to do that for model RHPs associated with stationary points of primary or secondary order. The proposition then implies the unique solvability of these RHPs.

\begin{proposition}\label{dzconverse} If $(I+w^+)^{\pm 1}$ and $(I-w^-)^{\pm 1}$ are in $L^\infty(\mathbb R)$, and (\ref{basiceqn}) has a solution $\mu \in I + L^2(\mathbb R)$ such that both $\mu$ and its matrix inverse $\mu^{-1}$ are in $L^\infty(\mathbb R)$, then $(1-C_w)$ is invertible on $L^p(\mathbb R)$ for any $1<p<\infty$. Furthermore,
\begin{eqnarray*}
  \|(1-C_w)^{-1}\|_{L^p\to L^p} &\lesssim_p& \|\mu\|_\infty \|\mu^{-1}\|_\infty \Big( \|(I + w^+)^{-1}\|_\infty + \|(I - w^-)^{-1}\|_\infty \Big)
\end{eqnarray*}
\end{proposition}
\proof[Proof] Let $M_\pm = I+C_\pm\big((w^+ + w^-)\mu\big) \equiv \mu(I \pm w^\pm) \in L^\infty$ as in the theorem of Beals and Coifman. For any matrix valued $f\in L^p$ we have:
\begin{eqnarray*}
  (1-C_w)f &=& C_+\Big(f(1-w^-)\Big) + C_-\Big(-f(1+w^+)\Big)\\
           &=& C_+\Big(f\mu^{-1} M_-\Big) - C_-\Big(f\mu^{-1}M_+\Big)\\
           &=& (C_M \circ\mu^{-1})(f)
\end{eqnarray*}
here $\mu^{\pm 1}$ act by right multiplications on $L^p$ (these are bounded since $\mu^{\pm 1} \in L^\infty$); $C_M$ and $C_{M^{-1}}$ denote Beals-Coifman operators with $L^\infty$ weights $(M_-,-M_+)$ and $(M^{-1}_-,- M^{-1}_+)$.

Using Liouville's theorem, it is not hard to see that $\det(M)\equiv 1$. By (\ref{factorize3}),
\begin{eqnarray*}
  C_M C_{M^{-1}} =  C_{M^{-1}} C_M = 1
\end{eqnarray*}
Consequently, $(1-C_w)^{-1} = \mu \circ C_{M^{-1}}$, and so its norm is controlled by
\begin{eqnarray*}
    \|(1-C_w)^{-1}\| &\lesssim& \|\mu\|_\infty \|C_{M^{-1}}\| \\
                     &\lesssim& \|\mu\|_\infty (\|M^{-1}_+\|_\infty + \|M^{-1}_-\|_\infty)\\
                     &\lesssim& \|\mu\|_\infty \|\mu^{-1}\|_\infty \Big( \|(I + w^+)^{-1}\|_\infty + \|(I - w^-)^{-1}\|_\infty \Big)
\end{eqnarray*}
\endproof

\section{Localization schemes and essence of the nonlinear stationary phase method} \label{localizationsection}

In this section, we'll describe the two schemes used to localize our oscillatory RHP (\ref{oscRHP},\ref{oscjumpmatrix}) to small neighborhood of the stationary points. In Section~\ref{essencesect} we shall discuss potential adaptations of these schemes to other oscillatory Riemann-Hilbert settings.

Recall that the oscillatory jump matrix is given by (\ref{oscjumpmatrix}), with the following canonical factorization:
\begin{eqnarray}
\label{canonical} J &=& \begin{pmatrix}1 & pe^{-it\theta}\cr 0 & 1\end{pmatrix}\begin{pmatrix}1 & 0 \cr qe^{it\theta} & 1\end{pmatrix} = (I-w^-_1)^{-1}(I+w^+_1)
\end{eqnarray}
i.e. $w^-_1= \begin{pmatrix}0 & pe^{-it\theta}\cr 0 & 0\end{pmatrix}$ and $w^+_1= \begin{pmatrix}0 & 0 \cr qe^{it\theta} & 0\end{pmatrix}$. This factorization is not always suitable for our reductions due to the lack of certain phase-weight relation (a suitable scalar RHP is needed to get the right factorization; Section~\ref{factorizesection} will be devoted to this discussion), but here we'll \emph{use it to demonstrate the ideas} behind our schemes, avoiding unnecessary technicalities.

Now, a perturbation of $p$ or $q$ leads to a weight perturbation by a term $e^{\pm it\theta}h$. Such oscillatory term will be decomposed into two components, which are handled separately using two schemes. We assume throughout that the weights are in $L^\infty\cap L^2$.

\subsection{The direct scheme} The Beals-Coifman formulation suggests that small changes in the weights won't have any strong effect on the solution of the RHP and hence $u(t)$, $v(t)$. Let $(\Delta w^- ,\Delta w^+)$ be added to $(w^-_1, w^+_1)$ to get a new pair of weights $(w^-_2, w^+_2)$. The effect on $u,v$ can be quantified by the non-tangential limit $\limsup_{\lambda \to \infty} |\lambda(M_1 (\lambda)- M_2(\lambda))|$ and will be estimated below.

For any pair $(w^-,w^+)$, for convenience we denote $w:=w^+ + w^-$ and
$$\|w\|_p := \|w^+\|_p + \|w^-\|_p$$
$$H_p(w) := \|C_-(w^+)\|_p + \|C_+(w^-)\|_p$$

\begin{lemma}
\label{direct_lem}
Suppose that $w_1$ is a pair of strictly-triangular weights and $\|w_1\|_\infty \lesssim 1$ as $t\to\infty$.

(i) Suppose that $\|(1-C_{w_i})^{-1}\|_{L^2\to L^2} \lesssim 1$ as $t\to\infty$ for each $i= 1,2$. Then the effect on $u(t), v(t)$ can be controlled by
\begin{eqnarray}
\nonumber           && \|\Delta w\|_2\Big(H_2(w_2)+H_2(w_1)\Big) \;\;+\;\; \|\Delta w\|_\infty H_2(w_2)H_2(w_1)\\
\label{direct_est}  &+& |\int \Delta w^+ w^+_1 |  + | \int \Delta w^- w^-_1|\\
\nonumber           &+& \limsup_{\lambda \to \infty} \big| \lambda C (\Delta w)(\lambda) \big|
\end{eqnarray}
Moreover, $\|\mu_1-\mu_2\|_2 \; \lesssim \; \|\Delta w\|_2 \; + \;  \|\Delta w\|_\infty  H_2(w_2)$.

(ii) Suppose that $\|(1-C_{w_1})^{-1}\|_{L^2\to L^2} \lesssim 1$ and $\|(1-C_{w_2})^{-1}\|_{L^p\to L^p} \lesssim_p 1$ for $2\leq p <\infty$ as $t\to\infty$. Then the above estimates remain true if we replace $\|\Delta w\|_\infty H_2(w_2)$ by $\|\Delta w\|_{\frac{2p}{p-2}}H_{p}(w_2)$.
\end{lemma}

\noindent \emph{Notes:} The third term in (\ref{direct_est}) is $0$ if the matrix structure of $w_2$ resembles that of $w_1$. In practice, some noise may affect the structure of $w_2$ so it won't necessarily be $0$, but will be very small for large $t$.

\proof (i) Let $\lambda \not\in\mathbb R$. Write $M_2(\lambda)-M_1(\lambda)$ as
\begin{eqnarray*}
C(\mu_2 w_2) - C(\mu_1 w_1) =
 C\Big((\mu_2-\mu_1)w_1\Big) + C\Big((\mu_2-I)\Delta w\Big)  + C(\Delta w)
\end{eqnarray*}
Note that if $f\in L^1\cap L^2$ then $\lim_{\lambda \to \infty} \lambda (C f)(\lambda) = -\frac{1}{2\pi i}\int f(x) dx$ assuming non-tangential limit, and $\|\mu_2-I\|_2 = \|(1-C_{w_2})^{-1}C_{w_2}I\|_2 \lesssim H_2(w_2)$. Thus, by Cauchy-Schwartz, to show (\ref{direct_est}) it remains to estimate $\int (\mu_2-\mu_1) w_1$.

Since $\mu_2 - \mu_1 = C_{\Delta w} \mu_2 + C_{w_1}(\mu_2-\mu_1)$, we get
$$\mu_2 - \mu_1 = (1-C_{w_1})^{-1}C_{\Delta w}\mu_2 = C_{\Delta w}\mu_2 + C_{w_1} \phi$$
where $\phi = (1-C_{w_1})^{-1}C_{\Delta w}\mu_2$. Observe that for any sign combination
$$\|C_\pm (\mu_2 \Delta w^\pm)\|_2 \;\;\lesssim\;\; \|\Delta w\|_2 + \|\Delta w\|_\infty \|\mu_2-I\|_2 \;\;\lesssim\;\; \|\Delta w\|_2 + \|\Delta w\|_\infty H_2(w_2).$$
Thus, $\|\phi\|_2$ can be controlled by the same estimate. Since $\|w_1\|_\infty \lesssim 1$, the desired estimate for $\mu_1-\mu_2$ follows immediately.

Now, using $L^2$ orthogonality of the two Hardy spaces $H_2(\mathbb C_+), H_2(\mathbb C_-)$, we can rewrite $\int  C_{\Delta w}\mu_2 \; w_1$ as
\begin{eqnarray}
\label{twointegrals}
 \int C_+(\mu_2 \Delta w^-) \Big(C_-(w^+_1) + w^-_1\Big) + \int C_-(\mu_2 \Delta w^+) \Big(w^+_1 + C_+(w^-_1)\Big) \;\;
\end{eqnarray}
Using Cauchy-Schwarz, the contributions of $C_-(w^+_1)$ in the first integral and of $C_+(w^-_1)$ in the second integral are controlled by
\begin{eqnarray*}
\Big(\|\Delta w\|_2 +  \|\Delta w\|_\infty H_2(w_2)\Big)H_2(w_1)
\end{eqnarray*}
Below we estimate the contribution of $w^-_1$ in the first integral of (\ref{twointegrals}). The contribution of $w^+_1$ in the second integral can be estimated similarly. Now,
\begin{eqnarray*}
 \int C_+(\mu_2 \Delta w^-) w^-_1
 &=& \int \mu_2 \Delta w^- w^-_1  + \int  C_-(\mu_2 \Delta w^-) w^-_1\\
 &=& \int \Delta w^- w^-_1  + \int (\mu_2-I) \Delta w^- w^-_1  + \int C_-(\mu_2 \Delta w^-) C_+(w^-_1)
\end{eqnarray*}
The last term can be controlled as before. The second term can be controlled using Cauchy-Schwarz and $\|w_1\|_\infty \lesssim 1$. The eventual estimate is
$$|\int \Delta w^- w^-_1| + \|\Delta w\|_2 H_2(w_2) + \Big(\|\Delta w\|_2 +  \|\Delta w\|_\infty H_2(w_2)\Big)H_2(w_1)$$

Similarly, $\int C_{w_1}\phi \;w_1$ can be estimated by
\begin{eqnarray*}
\Big(\|\Delta w\|_2 + \|\Delta w\|_\infty H_2(w_2)\Big)H_2(w_1) + |\int \phi w^-_1 w^-_1| + |\int \phi w_1^+ w^+_1|
\end{eqnarray*}
Using strict-triangularity of $w_1$, we have $w^-_1 w^-_1 = w^+_1 w^+_1 = 0$. This completes the proof of (i).

(ii) The only difference is in the estimate of $\|C_\pm (\mu_2\Delta w^\pm)\|_2$, where we'll use H\"older's inequality to control it by $\|\Delta w\|_2 + \|\Delta w\|_{\frac{2p}{p-2}}\|\mu_2-I\|_p \lesssim \|\Delta w\|_2 + \|\Delta w\|_{\frac{2p}{p-2}} H_p(w_2)$.
\endproof

\noindent \emph{Remarks:} 1. The above argument also proves that $\|\mu_1 - \mu_2\|_2   \lesssim \|C_{\Delta w}I\|_2 + \|\Delta w\|_\infty\|\mu_2 - I\|_2$. This is useful if we want to \emph{directly} exploit oscillation in $\Delta w$.

2. For $H_p(w)$ to decay as $t\to\infty$, we'll see that certain phase-weight relation in the pair $(w^-,w^+)$ is desired.

3. The last term $\limsup_{\lambda \to \infty} \big| \lambda C (\Delta w)(\lambda) \big|$ can be crudely controlled by $\|\Delta w\|_1$, but in our applications $\Delta w$ may consist of Hardy components (which might not even be in $L^1$). In that case, we will exploit cancellation of these terms under the Cauchy transform.

4. If $\|\Delta w\|_\infty$ decays as $t\to\infty$ then it suffices to assume uniform boundedness (as $t\to\infty$) of $\|(1-C_{w_1})^{-1}\|$. For sharp applications we will try to avoid estimates involving $\|\Delta w\|_\infty$ because it generally requires stronger regularity assumptions on $w$ (hence $p,q$). This means if applicable (ii) is better than (i). For this reason, we'll often use weaker assumptions on the resolvent operators. For instance, in part (ii) we might change the assumption on $(1-C_{w_2})^{-1}$ to
\begin{itemize}
\item $(1-C_{w_2})^{-1}$ is bounded from $L^{p+} \cap L^p$ to $L^p$ uniformly as $t\to\infty$,  for $2\leq p<\infty$.
\end{itemize}
The above argument then shows that $$\|\mu_1-\mu_2\|_2 \; \lesssim \;\|\Delta w\|_2 \; + \;  \|\Delta w\|_{\frac{2q}{q-2}}  \Big(H_{q}(w_2) + H_{q+}(w_2)\Big)$$
and the first two terms in (\ref{direct_est}) should be adjusted to
\begin{eqnarray*}
  \|\Delta w\|_2\Big(H_2(w_2) + H_{2+}(w_2)+H_2(w_1)\Big) \;\;+\;\; \|\Delta w\|_{\frac{2q}{q-2}} \Big(H_{q+}(w_2)+ H_q(w_2)\Big)H_2(w_1)
\end{eqnarray*}
and the eventual effect is only an $\epsilon$ in the order of our decay estimates.

\subsection{The indirect scheme}
The goal of the indirect scheme is to perturb the weights by {\it large but very well-structured} terms. This scheme is intuitively a reverse process of the direct scheme: instead of modifying the weights, we'll modify the Riemann-Hilbert factorization $J=M_-^{-1}M_+$ from \emph{outside} so that
\begin{itemize}
\item The leading asymptotics of $u(t),v(t)$ are not destroyed.
\item The effect on $J$ is essentially an addition of these large terms to the respective weights.
\end{itemize}

Below, we demonstrate this scheme by showing how to perturb $w_1^+$ by $\Phi_+$ a large but well-structured term. Typically, $\Phi_+$ is an analytic function on $\mathbb C_+$ which vanishes as $z\to\infty$ non-tangentially, furthermore its boundary value on $\mathbb R$ is in $L^2 \cap L^\infty$ and oscillates as $t\to\infty$. It is important that $\Phi_+$ continutes analytically to the \textit{upper} half-plane if we want to perturb it from the positive weight $w^+_1$. Now, consider the following normalized $L^2$-Riemann-Hilbert factorization
$$\widetilde{J} = M_-^{-1} \widetilde M_+$$
where $\widetilde{M}_+ := M(I+\Phi_+)$ and $\widetilde J = J(I+\Phi_+)$. Typically, the matrix structure of $\Phi_+$ will be similar to that of $w^+$ (which in our applications will be essentially strictly-upper or strictly-lower triangular); after all we are only working with the only non-zero entry of this weight. This implies
\begin{eqnarray*}
\widetilde J &=& (I-w_1^-)^{-1}(I+w_1^+ + \Phi_+ + w_1^+ \Phi_+)\\
             &\approx& (I-w_1^-)^{-1}(I+w_1^+ + \Phi_+)
\end{eqnarray*}
so the essential effect on $w_1^+$ is an addition by $\Phi_+$. In practice (when we are not using the canonical factorization (\ref{canonical})) it is possible that the triangular structure of $w_1,\Phi_+$ are affected by some ``noise'' terms (that decay strongly as $t\to\infty$), but these noise terms are indeed small and can be easily handled by the direct perturbation scheme.

Now, the effect of this modification on $u(t), v(t)$ can be controlled by
$$\limsup_{z\to\infty} |z\Big(M_+(z) - \widetilde{M}_+(z)\Big)| = \limsup_{z\to\infty} |z\Phi_+(z)|$$
Using Cauchy theorem, the last limit can be written as
$$\limsup_{z\to\infty} \Big|z \frac{1}{\pi i}\int_{\mathbb R} \frac{\Phi_+(x)}{x-z}dx\Big| \approx \frac{1}{\pi}\Big|\int_{\mathbb R} \Phi_+(x) dx\Big|$$
which decays as $t\to\infty$ thanks to oscillation of $\Phi_+$ on $\mathbb R$. The last equation is a heuristics (since $\Phi_+$ may not be in $L^1$) but it captures the essential idea.

\subsection{Combination of two schemes}

To perturb the weights by an oscillatory function $h(x)e^{it\Theta(x)}$ (here $\Theta = \pm \theta$), following \cite{varzuginyd} we'll decompose it using the Hardy decomposition:
$$h(x)e^{it\Theta(x)} = C_+(he^{it\Theta})(x) - C_-(he^{it\Theta})(x)$$
Under mild assumptions on $h$'s regularity, it turns out that if $\Theta'$ keeps the same sign on supp$(h)$ then one Hardy component is small for large $t$, and can be perturbed by the first scheme. The remaining term, consequently, inherits the oscillation of $he^{it\Theta}$ but now enjoys analytic continuation to one half-plane. It turns out that to ensure that this remaining term has continuation to the {\it correct} half plane (so that it can be perturbed by the second scheme), we want to have a certain phase-weight relation in the pair of weights $(w^-,w^+)$.

In Section~\ref{hardyspacesection} estimates on Hardy components of oscillatory functions are proved, which then explains the reason why a correct phase-weight relation is desirable.

We make an important observation that in the second scheme, the solution (if uniquely exists) $\mu$ to the basic equation $\mu = I + C_w \mu$ stays unchanged, thanks to the formula $\mu = M_+(I + w^+)^{-1}  = M_-(I - w^-) $. In the first perturbation scheme, we have an estimate controlling $\|\Delta \mu\|_2$. Consequently, during the reductions, $\|\Delta \mu\|_2$ always remains small. We'll see later that the above two schemes are analogous to the two main steps of the steepest descent method of Deift and Zhou: the second scheme is analogous to a contour deformation and the first scheme is analogous to a steepest descent estimate on the deformed contour.

\subsection{Essence of the nonlinear stationary phase method}\label{essencesect}

In this subsection, we shall discuss further the underlying ideas of the nonlinear stationary phase method developed in this paper and how they may be adapted to other oscillatory Riemann-Hilbert settings, in particular those that have been successfully studied by the steepest descent methods of Deift, Venakides, Zhou. Successful adaptations of these ideas will allow us to reduce the analyticity assumptions to fairly weak regularity assumptions on the relevant Riemann-Hilbert data in these applications.

To understand asymptotics of a given oscillatory RHP, the general strategy is to exploit the oscillation of oscillatory terms to show that the given RHP converges to a model RHP as the large parameter goes to $\infty$, and the solution of the limiting model RHP will provide the desired asymptotics of the quantity of interests. For instance, in the AKNS setting, if the phase $\theta$ has only one stationary point then the limiting RHP is a local model RHP associated with that stationary point, and when there are more than one stationary point we have to work a little harder to separate their contributions. In the small-dispersion KdV and orthogonal polynomial settings, the limiting RHPs may have a multi-interval configuration, more precisely the jump matrix will be locally constant \cite{DeiftZhouVenakidesyd, DeiftetalUniformOPyd}. The finiteness of the number of intervals is implied by real analyticity \cite{DeiftZhouVenakidesyd, DeiftetalUniformOPyd} or certain convexity of relevant Riemann-Hilbert data \cite{Benkoyd} and we shall assume this finiteness in subsequent discussions. For these multi-interval model RHPs, often the contribution of different intervals are not separated and the usual strategy is to construct an explicit solution for the multi-interval configuration \cite{DeiftetalUniformOPyd, DeiftZhouVenakidesyd}.

The steepest descent methods of Deift, Venakides, Zhou achieve the above convergence through exploiting the strong decay of the oscillating terms in appropriate deformation contours. As briefly mentioned in the Introduction section, often the jump matrix of the given oscillatory RHP may not be \emph{ready} for the steepest descent argument. The sense in which the jump matrix is \emph{ready} is better explained through an example. Suppose on a subinterval $I$ of the given Riemann-Hilbert contour (which we assume $\mathbb R$ for simplicity) the jump matrix $J$ has a factorization
\begin{equation}\label{3termfact}
J(\lambda,\epsilon) = v^-(\lambda,\epsilon) v(\lambda,\epsilon) v^+(\lambda,\epsilon)\end{equation}
where $1/\epsilon$ is the large parameter and $v^\pm$ consists of oscillatory terms $e^{i\theta/\epsilon}$ that we want to remove. As part of the Deift-Venakides-Zhou steepest descent scheme, the contour of the given RHP will be deformed, moving $v^+$ to the upper half plane and $v^-$ to the lower half plane where the oscillating terms decays strongly as $1/\epsilon \to\infty$. Consequently $v^+$ and $v^-$ converge to constant matrices as the large parameter goes to $\infty$, therefore if $\lim_{1/\epsilon\to\infty} v(\lambda,\epsilon)$ also exists then we can determine exactly the limiting jump matrix on the current interval $I$ (here we ignore all the endpoint issues where strong decay is not available and the use of certain parametrix might be needed).

In the AKNS setting, an function $\delta$ solving a scalar RHP is used to conjugate the given jump matrix
$$J_{new} = e^{\sigma_3 \delta_-} J_{old} e^{-\sigma_3 \delta_+}$$
and the scalar RHP is chosen so that the new jump matrix is ready for the steepest descent argument. In the small-dispersion KdV and orthogonal polynomial settings, this is achieved by more sophisticated versions of $\delta$ often referred to as the $g$-functions (which goes back at least to \cite{DVZShockCollisionyd}). 

The $g$-function for a given oscillatory RHP is often determined by a reverse-engineering process, as follows. Suppose that $e^{i\theta(\lambda)/\epsilon}$ is an oscillatory term in $J_{new}$ that one plans to move from an interval $I \subset \mathbb R$ to a contour in the upper half plane. In other words in the factorization (\ref{3termfact}) this oscillatory term would be part of $v^+(\lambda,\epsilon)$. It follows from the Cauchy-Riemann equation that if $\theta$ is (real) analytic then having $\theta'(\lambda)>0$ on $I$ will ensure strong decay of $e^{i\theta(\lambda)/\epsilon}$ on the new contour as the large parameter $1/\epsilon$ goes to $\infty$ (this observation goes back at least to \cite{DeiftZhouVenakidesyd}, see also the beginning of Section~\ref{reductmodelsection2} of this paper for a philosophical discussion). Similarly, if one plans to move $e^{i\theta/\epsilon}$ to the lower half plane then it is desirable that $\theta'(\lambda)<0$ on $I$. These conditions on the oscillating phases translate back to conditions on $g$ (there are also other conditions on $g$, for instance in the factorization (\ref{3termfact}) of $J_{new}$ we also want $v(\lambda,\epsilon)$ to converge as $1/\epsilon \to \infty$, but this is a different issue).

As discussed in the last two subsections, it is essential for the two perturbation schemes used in this paper that the weights in our Beals-Coifman factorization have correct phase-weight relation. In Section~\ref{hardyspacesection}, it will be shown that the correct phase-weight relation is 
\begin{itemize}
\item[(i)] The oscillating phase in $w^+$ is locally increasing.
\item[(ii)] The oscillating phase in $w^-$ is locally decreasing.
\end{itemize}
(for details see the remarks after the proof of Lemma~\ref{localprpl}). Writing $J=(I-w^-)^{-1}(I+w^+) \equiv v^-v^+$, one sees that the above phase-weight relation is exactly the monotonicity required of the phases to ensure the success of the steepest-descent argument of Deift and Zhou. One can look at this paper from the following angle: in the AKNS setting, given the setup of the steepest descent argument of Deift and Zhou, one can prove the desired convergence to model RHPs using the above two perturbation schemes. In other oscillatory RH settings, given correct phase monotonicity in the correct factors of $J_{new}$ (which in turn is achieved by the use of $g$-functions), one should be able to establish the desired convergence (to the model RHP) by adapting these two schemes to those settings. In other words, the only place where the method developed in this paper differs significantly from the Deift-Venakides-Zhou methodology is in the way one achieves the convergence to the model RHP, where one exploits the cancellation resulted from rapid oscillation of the oscillatory terms away from their stationary points/intervals by the use of real-variable tools.

To apply the two perturbation schemes described in the previous two sections to other settings, certain adaptations might be required. In settings where the quantity of interests is recovered from the Riemann-Hilbert solution by a limit similar to (\ref{recoverlimit}) (e.g. small-dispersion KdV), the second perturbation scheme should be immediately applicable provided that certain triangularity are available (which is often the case). In settings where one is more interested in the asymptotics of certain entries of the Riemann-Hilbert solution (e.g. orthogonal polynomials), the relative position of the oscillatory factors and the entries of interests should be taken into account. For adaptations of the first perturbation scheme, one observes that this scheme is based on the philosophy that Riemann-Hilbert solutions depend ``continuously'' on the corresponding jump matrices, and in the AKNS setting this philosophy is quantified by the use of the Beals-Coifman operator for normalized $L^2$-RHP. For applications to other settings, we may have to adapt the Beals-Coifman formulation to RHPs with different normalizations or different factorization structures. Note that in the implementation of the first scheme in this paper as previously described, triangularity of the weights has been extensively exploited to minimize the regularity assumptions on the Riemann-Hilbert data, but it is not essential.

The next section provides the estimates needed for our perturbation schemes.

\section{Hardy decomposition of oscillatory functions} \label{hardyspacesection}

Let $k\geq 1$. In this section, $\theta(x)$ will denote a real valued function such that:

(i) $\theta$ is $(k-1)$-time continuously differentiable everywhere. On the complement of a finite set of points, $\theta$ has $(k+1)$ locally integrable derivatives.

(ii) $\theta$ has stationary points $\lambda_1,\dots,\lambda_N$ of order $k_1,\dots,k_N$ i.e. $\theta$ is $(k_j+1)$-time differentiable near $\lambda_j$ with $\theta'(\lambda_j) =\dots =\theta^{(k_j)}(\lambda_j)=0 \neq \theta^{(k_j+1)}(\lambda_j)$, and $\theta^{(k_j+1)}$ is continuous at $\lambda_j$.

The value of $k$ will depend on the actual proposition, but it should be clear from the context. In this section, we'll refer to (i,ii) as the two phase conditions. The assumptions (A), (B) for Theorem~\ref{maintheorem} is for $k=2$.

We'll denote by $x_1<\dots<x_S$ the points mentioned in (i).

We'll need several notions of decay relative to $\theta$. For any $k$-time differentiable function $f$, let $f_{\theta,k}$ be the following weighted sum of its derivatives:
\begin{eqnarray*}
     &\sum_{\beta=0}^k |f^{(\beta)}| \cdot w_{k,\beta}(\theta)&\\
\text{where} &w_{k,\beta}(\theta) := \sum_{\alpha_i} |\theta'|^{-(k+ \alpha_1 + \dots + \alpha_k)} |\theta^{(2)}|^{\alpha_1} \dots |\theta^{(k+1)}|^{\alpha_k}&
\end{eqnarray*}
here the sum is taken over all nonnegative integers $\alpha_1, \dots, \alpha_k$ such that
\begin{eqnarray}
\label{alphaconstraint}     \alpha_1 + 2\alpha_2 + \dots + k\alpha_k =  k -\beta
\end{eqnarray}
(to avoid the $0^0$ situation, the convention is: if some $\alpha_j$ is $0$ then the respective power $|\theta^{(j+1)}|^{\alpha_j}$ is treated as $1$). The terms in the above sum appear naturally when we perform integration by parts on certain oscillatory integrals. If asymptotically $\theta$ and its derivatives behave like a polynomial of degree $d$ then $w_{k,\beta}(\theta) \lesssim |x|^{\beta-kd}$ and we can estimate $f_{\theta,k}$ more explicitly. On the other hand, if $f$ is compactly supported and supported away from the stationary points then
\begin{eqnarray}
\label{fthetak_est} f_{\theta,k} \lesssim |f(x)\theta^{(k+1)}(x)| + \sum_{\beta=0}^k |f^{(\beta)}(x)|
\end{eqnarray}

For any $d\in \mathbb R$, let $D(\theta,d,k)$ be the set of $k$-time differentiable functions $f$ satisfying the following two conditions:
\begin{itemize}
  \item $\lim_{x\to \infty}(\theta')^{-1}f_{\theta,j}(x) = 0$ for every $0\leq j\leq k-1$
  \item $|\theta'|^df_{\theta,k}$ is integrable at infinity
\end{itemize}
In particular, every compactly supported function satisfies the decay requirements $D(\theta,d,k)$ for any $\theta,d$.

The next lemma extends Varzugin's localization principle in \cite{varzuginyd}, here some modifications are made to allow for less stringent assumptions on $\theta$ and for a wider class of $f$.

\begin{lemma} \label{localprpl} Let $2\leq p < \infty$. Suppose $f \in L^2$ is $k$-time differentiable and supported in $\{\theta' \geq 0\}$ s.t. $f, f', \dots, f^{(k-1)}$ vanish on the boundary of $\{\theta' > 0\}$. If $f$ has sufficient decay (say, $D(\theta,1-\frac{1}{p},k)$) and vanishes at every endpoint $\lambda_j$'s of $\{\theta'> 0\}$ with high multiplicity, then
\begin{eqnarray*}
    \|C_-(fe^{it\theta})\|_p &\lesssim& t^{-(k-1+1/p)}
\end{eqnarray*}
If $k\geq 2$, this estimate is also true for the endpoint case $p=\infty$.
\end{lemma}

\noindent \emph{Remarks:} If supp$(f) \subseteq \{\theta' \leq 0\}$, the estimate is true for $C_+(fe^{it\theta})$ by symmetry.

\proof We'll largely follow the argument in \cite{varzuginyd}. Indeed, we'll show, $\forall 2\leq p\leq \infty$
\begin{eqnarray}
\label{HardyLp}    \|C_-(fe^{it\theta})\|_p \;\lesssim_{p,\theta,k}\; t^{-(k-1+1/p)}\Big(\big\||\theta'|^{1-1/p}f_{\theta,k}\big\|_1 + \sum_{n=1}^S \sum_{\beta=0}^{k-1} |f^{(\beta)}(x_n)|\Big)
\end{eqnarray}
as long as $kp'>1$. Basically, this means for $p=\infty$ we need $k\geq 2$, while for $p < \infty$ we only require $k\geq 1$, and these are exactly the conditions on $k$ described above.

First, using Hausdorff-Young inequality, for every $2\leq p \leq \infty$ we have:
\begin{eqnarray*}
       \|C_-(fe^{it\theta})\|_{L^p(\mathbb R)}      &\lesssim& \|\widehat {fe^{it\theta}}(-\xi)\|_{L^{p'}_\xi([0,\infty))}
\end{eqnarray*}

Indeed, this is the only place the condition $2\leq p\leq \infty$  (i.e. $1\leq p'\leq 2$) is really needed. The rest of this argument works for $1\leq p \leq \infty$. This remark will be useful for future computation, in particular the proof of Lemma~\ref{lineartheory}.

For $\xi \geq 0$, write $\widehat {fe^{it\theta}}(-\xi) =(2\pi)^{-1} \int_{-\infty}^\infty e^{i(x\xi +t\theta(x))} f(x)dx$, which looks like a standard oscillatory integral, except for the term $x\xi$. This term is however harmless because for $x \in$ supp($f$), $\xi$ and $\theta'(x)$ are of the same sign and so:
\begin{eqnarray}
\label{phasew} |\frac{d}{dx}(x\xi +t\theta(x))| = |\xi +t\theta'(x)| \geq t|\theta'| \geq 0
\end{eqnarray}
Thus, if $f$ behaves nicely on the boundary of $\{\theta' > 0\}$ (which are essentially stationary points $\lambda_j$ of $\theta$), we expect this integral to decay as $t\to\infty$.

Denote $w(x) = \xi +t\theta'(x)$ (we suppress $t,\xi$ in $w$ for simplicity). Let $D_w$ be the following operator that acts on differentiable functions:
\begin{eqnarray*}
    D_w f = \frac{d}{dx}\bigl(\frac{1}{w(x)}f\bigr)
\end{eqnarray*}
Integration by parts $k$ times gives us:
\begin{eqnarray*}
       \widehat {fe^{it\theta}}(-\xi)
       &=& \text{boundary terms} + i^k\int_{-\infty}^\infty (D_w^k f)(x) e^{i\xi x + it\theta(x)}dx
\end{eqnarray*}
The boundary terms are evaluations of $\frac{C(j)}{w}e^{i\xi x + it\theta(x)}D^j_w f$ at various endpoints for some absolute constant $C(j)$. For $0\leq j< k-1$ they are evaluated at $\pm\infty$ and the endpoints of $\{\theta'\geq 0\}$. For $j=k-1$ there are also evaluations at $x_1,\dots, x_S$ as left/right limits, since $\theta^{(k-1)}$ is not differentiable at these points. We'll see that the first two types of evaluations vanish, while the last is $O(t^{-(k-1+1/p)})$ in $L^{p'}_\xi$.

Indeed, the evaluations at endpoints of $\{\theta'\geq 0\}$ vanish since $f$ vanishes at $\lambda_j$ with high multiplicity. To show that the evaluations at $\pm\infty$ vanish, we will control these terms by showing
\begin{eqnarray}
\label{Dwestimate}    |D^j_w f| &\lesssim& t^{-j} f_{\theta,j}.
\end{eqnarray}
This estimate and the given assumption on decay of $f$ will then give the desired claim. To prove (\ref{Dwestimate}), first notice that every $(D_w^j f)(x)$ is a linear combination of
\begin{eqnarray}
\label{partintterm}    f^{(\beta)}(x)\cdot\frac{1}{w(x)^n}  w'(x)^{\alpha_1} \dots w^{(j)}(x)^{\alpha_j}
\end{eqnarray}
with $\alpha_i \geq 0$, $n\geq 0$. Now, our first observation is: $\xi$ no longer appears in the derivatives of $w$, indeed $w^{(i)} = t\theta^{(i+1)}$. The second observation is
\begin{eqnarray*}
     n = j+ \alpha_1 + \dots + \alpha_j
\end{eqnarray*}
This can be easily seen by a scaling symmetry argument: replace $w(x)$ by $c \cdot w(x)$ and notice that $D^n_{cw} = \frac{1}{c^n}D_w$, then do a counting of $c$. In addition, using a dilation symmetry argument (i.e. replace the pair $(w(x),f(x))$ by $(\widetilde w(x) := w(cx), \widetilde f(x) := f(cx))$ and notice that $D_{\widetilde w} \widetilde f (x) = cD_w f(cx)$) we can also show that:
\begin{eqnarray*}
     j = \beta+\alpha_1 + 2\alpha_2 + \dots + j\alpha_j
\end{eqnarray*}
Using (\ref{phasew}), we then can dominate (\ref{partintterm}) by
\begin{eqnarray*}
    &&|f^{(\beta)}|\cdot\frac{1}{|t\theta'|^n}  |t\theta^{(2)}|^{\alpha_1} \dots |t\theta^{(j+1)}|^{\alpha_j}\\
    &=& t^{-j} |f^{(\beta)}| \cdot \frac{1}{|\theta'|^{j+ \alpha_1 + \dots + \alpha_j}} |\theta^{(2)}|^{\alpha_1} \dots |\theta^{(j+1)}|^{\alpha_j}
\end{eqnarray*}
Thus, from the definition of $f_{\theta,j}$ we get the desired estimate (\ref{Dwestimate}). The last task is to estimate the boundary terms coming from $x_n$, which modulo an absolute constant are of the form
$$\lim_{x\to x^+_n} \frac{1}{w}e^{i\xi x + it\theta(x)}D^{k-1}_w f(x), \text{ and } \lim_{x\to x^-_n} \frac{1}{w}e^{i\xi x + it\theta(x)}D^{k-1}_w f(x)$$
These limits are nontrivial only if $\theta'(x_n)>0$, in which case their $L^{p'}_\xi(\mathbb R_+)$ norms can be controlled by the summation over $\{\alpha_1+\dots+\alpha_{k-1}+\beta=k-1\}$ of
\begin{eqnarray*}
&& |f^{(\beta)}(x_n)| \Big\|\frac{|t\theta^{(2)}(x_n)|^{\alpha_1}\dots |t\theta^{(k+1)}(x_n)|^{\alpha_{k-1}}}{(t\theta'(x_n)+\xi)^{k+\alpha_1+\dots +\alpha_{k-1}}}\Big\|_{L^{p'}_\xi([0,\infty))}\\
&\lesssim& \frac{|f^{(\beta)}(x_n)|t^{\alpha_1+\dots +\alpha_{k-1}}}{t^{k+\alpha_1+\dots+\alpha_{k-1}-\frac{1}{p'}}}\\
&=& t^{-(k-1+\frac{1}{p})}|f^{(\beta)}(x_n)|.
\end{eqnarray*}
Note that the condition $kp'>1$ guarantees convergence of the $L^{p'}_\xi$ norm.

By Minkowski's inequality, we can control $\|\widehat {fe^{it\theta}}(-\xi)\|_{L^{p'}_\xi([0,\infty))}$ by
\begin{eqnarray*}
&& \|\text{boundary terms}\|_{L^{p'}_\xi([0,\infty))} + \|\int_{-\infty}^\infty (D_w^k f)(x) e^{i\xi x + it\theta(x)}dx\|_{L^{p'}_\xi([0,\infty))}\\
&\lesssim& t^{-(k-1+1/p)}\sum_{n=1}^S \sum_{\beta=0}^{k-1} |f^{(\beta)}(x_n)| + \int_{-\infty}^\infty \|(D_w^k f)(x)\|_{L^{p'}_\xi([0,\infty))} dx
\end{eqnarray*}
Using the triangle inequality, $\|(D_w^k f)(x)\|_{L^{p'}_\xi([0,\infty))}$ can be controlled by:
\begin{eqnarray*}
  \sum_{\beta=0}^k |f^{(\beta)}(x)|\sum_{\alpha_i} |t\theta^{(2)}(x)|^{\alpha_1} \cdots|t\theta^{(k+1)}(x)|^{\alpha_k} \|\frac{1}{(t\theta'(x)+\xi)^{k+\alpha_1+\dots+\alpha_k}}\|_{L^{p'}_\xi([0,\infty))}
\end{eqnarray*}
Suppose $x\in$ supp($f$) is not a stationary points (so $\theta'(x) > 0$). Since $kp'>1$, the $L^{p'}$ norms in the last sum converge and we can rewrite this sum as
\begin{eqnarray*}
  &&\sum_{\beta=0}^k |f^{(\beta)}(x)|\sum_{\alpha_i} |t\theta^{(2)}(x)|^{\alpha_1} \cdots|t\theta^{(k+1)}(x)|^{\alpha_k}  \frac{\big(p'(k+\alpha_1+\dots+\alpha_k)-1\big)^{1/p'}}{|t\theta'(x)|^{(k+\alpha_1+\dots+\alpha_k)-1/p'}}\\
  &\lesssim&\sum_{\beta=0}^k |f^{(\beta)}(x)|\sum_{\alpha_i} |\theta^{(2)}(x)|^{\alpha_1} \cdots|\theta^{(k+1)}(x)|^{\alpha_k}  \frac{t^{-(k-\frac{1}{p'})}}{|\theta'(x)|^{(k+\alpha_1+\dots+\alpha_k)-1/p'}}\\
  &=& t^{-(k-1+\frac{1}{p})} |\theta'(x)|^{1/p'}f_{\theta,k}(x)
\end{eqnarray*}
Consequently, $\|\widehat {fe^{it\theta}}(-\xi)\|_{L^{p'}_\xi([0,\infty))}$ is controlled by
\begin{eqnarray*}
    &&t^{-(k-1+\frac{1}{p})}\sum_{n=1}^S \sum_{\beta=0}^{k-1} |f^{(\beta)}(x_n)| + \int_{-\infty}^\infty  t^{-(k-1+\frac{1}{p})} |\theta'|^{1/p'}f_{\theta,k}(x) dx\\
    &=& t^{-(k-1+\frac{1}{p})}\Big(\sum_{n=1}^S \sum_{\beta=0}^{k-1} |f^{(\beta)}(x_n)| + \Big\||\theta'|^{1/p'}f_{\theta,k}\Big\|_1\Big)
\end{eqnarray*}
giving the desired estimates (\ref{HardyLp}).

Clearly if $f\in D(\theta, \frac{1}{p'}, k)$ then $|\theta'|^{1/p'}f_{\theta,k}$ is integrable at $\infty$, so the $L^1$ integrability of this term depends on its behavior at the stationary points, where $\theta'=0$. Below, we'll show that if $f$ vanishes at every $\lambda_j$ with high multiplicity, then $|\theta'|^{1/p'}f_{\theta,k}$ is integrable at $\lambda_j$.

By a simple application of L'Hospital's rule, at $\lambda_j$, $\theta^{(n)}$ contributes a zero of order $k_j+1-n$ for every $1\leq n \leq k_j+1$:
\begin{eqnarray*}
  \lim_{x \to \lambda_j} \frac{\theta^{(n)}(x)}{(x - \lambda_j)^{k_j-n+1}} =  \frac{\theta^{(k_j+1)}(\lambda_j)}{(k_j+1-n)!} \neq 0
\end{eqnarray*}
If $k\leq k_j$, the order at $\lambda_j$ of $|\theta'|^{-(k+\alpha_1+\dots+\alpha_k)+1/p'}|\theta^{(2)}|^{\alpha_1} \cdots|\theta^{(k+1)}|^{\alpha_k}$ is:
\begin{eqnarray*}
    - k_j\Big(k + \sum_{i=1}^k \alpha_i - \frac{1}{p'}\Big)+\sum_{i=1}^k (k_j - i)\alpha_i \
    &=&  \beta  -(1 + k_j)k + \frac{k_j}{p'}
\end{eqnarray*}
using (\ref{alphaconstraint}). If $k>k_j$, the order at $\lambda_j$ of this term is at least that much, as
\begin{eqnarray*}
    \sum_{i=1}^{k_j} (k_j - i)\alpha_i
    &\geq& \sum_{i=1}^k (k_j - i)\alpha_i
\end{eqnarray*}
Since $\theta^{(k+1)}$ is locally integrable, $|\theta'|^{1/p'}f_{\theta,k}$ is integrable at $\lambda_j$ when $f$ ``vanishes'' there with high multiplicity.
\endproof

\noindent \emph{Remarks:} The above lemma indicates that: to ensure the success of the perturbation schemes in Section~\ref{localizationsection} we want the following phase-weight relation:
\begin{itemize}
\item The oscillating phase of $w^+$ is (locally) increasing.
\item The oscillating phase of $w^-$ is (locally) decreasing.
\end{itemize}
To see why, consider for example a perturbation of $w^+$ by $fe^{it\Theta}$, where $\Theta$ is the oscillating phase of $w^+$. Decompose according to our schemes
$$fe^{it\Theta} = C_+(\dots) - C_-(\dots)$$
If $\Theta$ is locally increasing then $C_+(\dots)$ is the large and oscillatory term, and it has analytic continuation to the {\it upper} half-plane and can be perturbed from $w^+$ using the indirect scheme. If $\Theta$ was locally decreasing then the large term will be $C_-(\dots)$, and we will not be able to perturb it away from $w^+$ using our second scheme because it has analytic continuation to the \emph{wrong} half-plane.

The canonical factorization (\ref{canonical}) may lack this phase-weight relation, as $\theta'$ is not always positive on $\mathbb R$. In Section~\ref{factorizesection}, we'll use a scalar RHP to get the right factorization. This scalar RHP is the source of the interaction between stationary points and the $\ln t$ term in the asymptotics in Theorem~\ref{maintheorem}. The introduction of this scalar RHP leads to the following type of functions:

\begin{definition}[$A_k+B_k$ decomposition] A measurable function $\Omega(x)$ has an $A_k+B_k$ decomposition (with respect to $\theta$) if it is $k$-time differentiable and has the following properties:
\begin{itemize}
\item[(i)] Boundedness: $\|\Omega\|_\infty \lesssim 1$

\item[(ii)] $\Omega^{(n)}(x) \lesssim 1 +\max_{1\leq j\leq N}\frac{1}{|x-\lambda_j|^n}$ for $1\leq n\leq k-1$ while $\Omega^{(k)}(x)$ can be decomposed as $A_k(x) +B_k(x)$, with
$\|A_k\|_2 \lesssim_k 1$ and $|B_k(x)| \lesssim_k 1+\max_{1\leq j\leq N}\frac{1}{|x-\lambda_j|^k}$
\end{itemize}
(we understand that $\max_{1\leq j\leq N}\frac{1}{|x-\lambda_j|^k}=0$ when $\theta$ has no stationary point).
\end{definition}

A trivial example is $\Omega \equiv 1$. Observe that if $\Omega$ has the above properties then so does $h\Omega$, for any $h\in H^{k,0}$. Furthermore, if $f \in H^{k,0}$ is compactly supported and supported away from the stationary points then using local integrability of $\theta^{(k+1)}$ it is not hard to see that
$$\|(f\Omega)_{\theta,k}\|_1 \lesssim \|f\|_{H^{k,0}}.$$
We'll sometimes say that a function has an $A_k+B_k$ decomposition on a set. This means that (i) and (ii) are only required to be true inside that set.

\begin{definition}[Vanish with multiplicity] A function $f$ vanishes at $\lambda$ with multiplicity $m \geq 0$ up to the $k^{\emph{th}}$ derivative on a neighborhood $P$ of $\lambda$ if:

(i) $f$ is $k$-time differentiable on $P$, and

(ii) $f^{(i)}(x) \lesssim |x-\lambda|^{m-i}$ for any $0 \leq i\leq k-1$, and

(iii) $\|(x-\lambda)^{k-m-\frac{1}{2}}f^{(k)}(x)\|_{L^2(P)} \lesssim 1$.
\end{definition}
We'll say that $f$ vanishes at $\lambda$ with multiplicity $m$ up to the $k^{\emph{th}}$ derivative if such $P$ exists. Note that $m$ is not necessarily an integer and (ii) is trivial if $i>m$. As an example, $f\in H^{1,0}$ has multiplicity $0$ up to the first derivative. If furthermore $f(\lambda)=0$ then the multiplicity can be improved to $\frac{1}{2}$, and if $f'$ is also bounded near $\lambda$ then the multiplicity can be improved to $1-\epsilon$.

\begin{lemma}\label{mainlemma} Let $f$ be $k$-time differentiable and supported near a stationary point $\lambda_0$ of order $k_0$, where it vanishes with multiplicity $m$ up to the $k^{\emph{th}}$ derivative. Let $\Omega$ have an $A_k+B_k$ decomposition on this neighborhood.

\noindent (i) Let $k\leq k_0$. Then for $2\leq p<\infty$,
\begin{eqnarray*}
 &\|C_-(\chi_{\{\theta' \geq 0\}}  f\Omega e^{it\theta})\|_p + \|C_+(\chi_{\{\theta' \leq 0\}}  f\Omega e^{it\theta})\|_p&\\
&\lesssim\;\;\; \max(t^{-(k-1+\frac{1}{p})}, t^{-(m +\frac{1}{p})\frac{1}{k_0+1}}  \ln t)&
\end{eqnarray*}
here $\ln t$ can be removed if $m +  \frac{1}{p} \neq (k_0+1)(k-1+\frac{1}{p})$. The estimate remains true for $p=\infty$ provided that $m>0$ and $(m +\frac{1}{\infty})\frac{1}{k_0+1} \equiv \frac{m}{k_0+1}$ is replaced by $\frac{m}{k_0+1} - \epsilon$.

\noindent (ii) Let $k\geq k_0+1$. If $\theta^{(k+1)}$ is $L^r$ integrable near $\lambda_0$ for some $1\leq r\leq \infty$ then the above estimates remain true after the following adjustments: $k_0+1$ is replaced by $k_0+1+c(k,p)$, and if $p=\infty$ then $m+\frac{1}{p} = m$ is replaced by $m-c(k,p)$. Here
\begin{eqnarray}
\label{ckp} c(k,p) = \begin{cases}\frac{1}{r(k-1+\frac{1}{p})}, & 2\leq p<\infty\\
\frac{1}{rk}, & p=\infty
\end{cases}
\end{eqnarray}
\end{lemma}

\proof By translation invariant, we can assume that our stationary point is $0$. For convenience, denote by $\{|x|\lesssim 1\}$ a small neighborhood of $0$ where $f$ is supported.

(i) It suffices to show the estimate for the first term on the left hand side.

\noindent \underline{\it Non-endpoint case:} Assume that $2\leq p<\infty$.

Let $\varphi$ be a $C^\infty$ cutoff taking values in $[0,1]$ such that $\varphi(x)=0$ for $|x|\geq 2$ and $\varphi(x)=1$ for $|x|\leq 1$. Let $\epsilon>0$ be a small scale.

Let $\varphi_\epsilon(x) := 1-\varphi(\frac{x}{\epsilon})$ and $g := \chi_{\{\theta'\geq 0\}}\varphi_\epsilon f\Omega$. Using the triangle inequality and applying estimate (\ref{HardyLp}) to $g$, we can control $\|C_-(\chi_{\{\theta'\geq 0\}}f\Omega e^{it\theta})\|_p$ by:
\begin{eqnarray*}
  && \|\varphi(\frac{.}{\epsilon})f\Omega\|_p +  \|C_-(g e^{it\theta})\|_p\\
  &\lesssim& \epsilon^{1/p} \sup_{|x|\leq \epsilon}|f(x)| +  t^{-(k-1+\frac{1}{p})} \big\||\theta'|^{1-1/p}g_{\theta,k}\big\|_1
\end{eqnarray*}
We'll optimize the estimate over $\epsilon$ (by carefully keeping track of their orders). It turns out that a good choice of $\epsilon$ decays as $t\to\infty$, so the assumption ``$\epsilon$ is small'' is harmless. Note that below $p$ is allowed to be in $[1,\infty]$.

For $0\leq \beta\leq k$, expand $g^{(\beta)}(x)$ using the product rule and estimate it by
\begin{eqnarray*}
g^{(\beta)}(x) &\lesssim& |(\Omega f)^{(\beta)}\varphi_\epsilon| + \sum_{(n,\alpha,\gamma)\in S_\beta}|f^{(n)} \Omega^{(\alpha)} \varphi^{(\gamma)}_\epsilon|
\end{eqnarray*}
where $S_\beta := \{(n,\alpha,\gamma)\in \mathbb Z_+^3: n+\alpha+\gamma=\beta, \gamma>0\}$.

Note that $\text{supp}(\varphi_\epsilon) \subset \{\epsilon \leq |x| \lesssim 1\}$, and for $\gamma>0$ we have $\varphi^{(\gamma)} \lesssim \epsilon^{-\gamma}1_{\{\epsilon \leq |x| \leq 2\epsilon\}}$. Thus, we can control $g^{(\beta)}$ by
\begin{eqnarray*}
|(\Omega f)^{(\beta)}|1_{\{\epsilon \leq |x| \lesssim 1\}} + \sum_{S_\beta}|f^{(n)}\Omega^{(\alpha)}| \epsilon^{-\gamma}1_{\{|x| \sim \epsilon\}}
\end{eqnarray*}

The last part of the proof of Lemma~\ref{localprpl} shows that near $0$ the $\theta$-dependent weights of $|g^{(\beta)}|$ in the weighted sum $|\theta'|^{1-\frac{1}{p}}g_{\theta,k}$ is controlled by
\begin{eqnarray*}
|x|^{\beta-(1+k_0)k+\frac{k_0}{p'}}
\end{eqnarray*}
notice that the assumption $k\leq k_0$ implies that $\theta^{(k+1)}$ is bounded near $0$.

For $|x| \sim \epsilon$ the weight of $g^{(\beta)}$ is of size $\epsilon^{\beta-(1+k_0)k+\frac{k_0}{p'}}$, so
$|\theta'|^{1-1/p}g_{\theta,k}$ is controlled by the summation over $0\leq \beta\leq k$ of $\text{I}(\beta) + \text{II}(\beta)$, where
\begin{eqnarray*}
\text{I}(\beta) &=& \sum_{n+\alpha=\beta} |f^{(n)} \Omega^{(\alpha)}| |x|^{\beta-(1+k_0)k+\frac{k_0}{p'}}1_{\{\epsilon \leq |x| \lesssim 1\}}\\
\text{II}(\beta) &=& \sum_{S_\beta}|f^{(n)}\Omega^{(\alpha)}| \epsilon^{-\gamma}\epsilon^{\beta-(1+k_0)k+\frac{k_0}{p'}}1_{\{|x| \sim \epsilon\}}
\end{eqnarray*}
By assumption,
\begin{eqnarray}
\label{derivative_est}
&&f^{(n)}(x) \lesssim |x|^{m-n} \;\;\; \text{ for $0\leq n \leq k-1$;}\\
\nonumber \text{and}&&\|1_{\{|x|\lesssim 1\}} x^{k-m-\frac{1}{2}}f^{(k)}(x)\|_{L^2} \lesssim 1
\end{eqnarray}

\noindent {\it Contribution of $I(\beta)$}: If $\alpha<k$, $n<k$ then we have a pointwise bound for $\Omega^{(\alpha)}$ and $f^{(n)}$. Thus, the contribution in $L^1$ of the associated term in $I(\beta)$ is controlled by:
$$\int_{\epsilon \leq |x| \lesssim 1} |x|^{m-n-\alpha+\beta-(1+k_0)k+\frac{k_0}{p'}} dx$$
Since $n + \alpha = \beta$, the above integral is controlled by
$$O(1)+\epsilon^{m+1-(1+k_0)k+\frac{k_0}{p'}}|\ln \epsilon|,$$
where $|\ln \epsilon|$ is not needed if $m+1\neq (1+k_0)k-\frac{k_0}{p'}$, which is equivalent to
\begin{eqnarray}
\label{nologconstraint} m+\frac{1}{p} \neq (k-1+\frac{1}{p})(k_0+1)
\end{eqnarray}
and $O(1)$ is needed only when $m+1 > (k_0+1)k-\frac{k_0}{p'}$, in which case it corresponds to a positive power of the size of our neighborhood of $0$. This means that we can make this $O(1)$ smaller by making our neighborhood smaller. These comments apply to similar estimates in the future.

Using Cauchy-Schwarz, the contribution of $(n,\alpha) = (k,0)$ can be estimated by
$$\|1_{\{|x|\lesssim 1\}}x^{k-m-\frac{1}{2}}f^{(k)}(x)\|_2 \Big(\int_{\{\epsilon\leq |x| \lesssim 1\}}  x^{2(m+\frac{1}{2}-(1+k_0)k+\frac{k_0}{p'})} dx\Big)^{\frac{1}{2}}$$
which is again controlled by $O(1)+\epsilon^{m+1-(1+k_0)k+\frac{k_0}{p'}}|\ln \epsilon|^{\frac{1}{2}}$.

To estimate the contribution of $(n,\alpha)=(0,k)$, we write $\Omega^{(k)} = A_k + B_k$ and notice first that for $B_k$ (which has nice pointwise bound) we can proceed as before. For $A_k$ (with bounded $L^2$ norm), by Cauchy-Schwarz
$$\int_{|x| \lesssim 1} |x|^{k-1}|A_k(x)| dx \lesssim_M 1$$
Using $m+1\leq (1+k_0)k-\frac{k_0}{p'}$, we can similarly control the contribution of $A_k$ by
$$\sup_{\epsilon\leq |x| \lesssim 1} |x|^{m+1 - (1+k_0)k + \frac{k_0}{p'}} \lesssim O(1) + \epsilon_j^{m+1-(1+k_0)k+\frac{k_0}{p'}}$$
Consequently, the contribution of $\sum_{\beta=0}^k I(\beta)$ in $L^1$ can be controlled by:
$$O(1)+ \epsilon^{m+1-(1+k_0)k+\frac{k_0}{p'}}|\ln \epsilon|$$

\noindent {\it Contribution of $II(\beta)$.} For any term in $II(\beta)$ clearly $\alpha<k$ and $n<k$, so using (\ref{derivative_est}) and pointwise bounds on derivatives of $\Omega$, we can estimate this term by:
\begin{eqnarray*}
\epsilon^{-\gamma + \beta-(1+k_0)k+\frac{k_0}{p'}} \int_{|x| \sim \epsilon} |x|^{m-n -\alpha} dx &\approx&  \epsilon^{m+1-(1+k_0)k+\frac{k_0}{p'}}
\end{eqnarray*}
\noindent {\bf Choice of $\epsilon$.} Now, we can estimate $\|C_-(\chi_{\{\theta'\geq 0\}}f\Omega e^{it\theta})\|_p$ by
\begin{eqnarray*}
\epsilon^\frac{1}{p}\epsilon^{m} + t^{-(k-1+\frac{1}{p})} \Big(O(1)+\epsilon^{m+1-(1+k_0)k+\frac{k_0}{p'}} |\ln \epsilon|\Big)
\end{eqnarray*}
The optimal choice $\epsilon =t^{-\frac{1}{k_0+1}}$ gives the desired estimate.

\noindent \underline{\it Endpoint case:} Consider $p=\infty$ under the extra assumption $m >0$ (this forces $f(0)=0$). Note that $\|C_-(ge^{it\theta})\|_\infty$ can be estimated as before, so the main task is to estimate:
$$\|C_-(\chi_{\{\theta'\geq 0\}}\varphi(\frac{.}{\epsilon})f\Omega e^{it\theta})\|_\infty$$
Let $h:= \chi_{\{\theta'\geq 0\}}\varphi(\frac{.}{\epsilon})f\Omega e^{it\theta}$. By standard Sobolev estimates, we have
$$\|C_-(h)\|_\infty \lesssim_q \|h\|_q + \|h'\|_q \;\;\forall 1<q \leq 2$$
Since $f(0) = 0$, it is not hard to see that $h$ is differentiable, with the following derivative:
$$\chi_{\{\theta'\geq 0\}}\Big(\frac{1}{\epsilon}\varphi'(\frac{.}{\epsilon})f\Omega + \varphi(\frac{.}{\epsilon})f'\Omega + \varphi(\frac{.}{\epsilon_j})f\Omega' + \varphi(\frac{.}{\epsilon})f\Omega t\theta'\Big)e^{it\theta}$$
Notice that $\theta'(x)$ is about the size of $|x|^{k_0}$ for $x$ near $0$. Since both $h$ and $h'$ are supported on $\{|x|\lesssim \epsilon\}$, using (\ref{derivative_est}) we can control $h'(x)$ by:
\begin{eqnarray*}
&&\frac{1}{\epsilon}|x|^{m} + |x|^{m-1} + |x|^{m}|x|^{-1}  + |x|^{m}t|x|^{k_0} \\
&&\lesssim |x|^{m-1} + t\epsilon^{m+k_0}
\end{eqnarray*}
Using $m > 0$ and choose $q$ sufficiently close to $1$, after integrating we get $\|h'\|_q \lesssim \epsilon^{m-1+\frac{1}{q}} + t\epsilon^{m+k_0+\frac{1}{q}}$. On the other hand, $\|h\|_q$ can be easily controlled by $\epsilon^{m+\frac{1}{q}}$. Consequently
$$\|h\|_q + \|h'\|_q \lesssim \epsilon^{m-1+\frac{1}{q}} + t\epsilon^{m+k_0+\frac{1}{q}}$$
The choice $\epsilon = t^{-\frac{1}{k_0+1}}$ gives an overall estimate
$$\|C_-(\chi_{\{\theta' \geq 0\}}  f\Omega e^{it\theta})\|_\infty \;\;\lesssim_q\;\; t^{-(m -1+\frac{1}{q})\frac{1}{k_0+1}}\ln t + t^{-(k-1)}$$
now choosing $q$ sufficiently close to $1$ we get the desired estimate.

(ii) The proof is similar. The estimate for the weight of $g\equiv g^{(0)}$ will involve $\theta^{(k+1)}$ (which is not assumed bounded) and this has to be estimated more carefully. Using H\"older's inequality, for that contribution we can get an estimate of
$$1+\epsilon^{m+\frac{1}{r'} - (k_0+1)k+\frac{k_0}{p'}}|\ln \epsilon|$$
instead of $1 + \epsilon^{m+1 - (k_0+1)k+\frac{k_0}{p'}}|\ln \epsilon|$. Eventually we have an overall estimate of
$$\epsilon^{m+\frac{1}{p}} + t^{-(k-1+\frac{1}{p})}\Big(1+\epsilon^{m+\frac{1}{r'} - (k_0+1)k+\frac{k_0}{p'}}|\ln \epsilon|\Big).$$
then optimizing over $\epsilon$ will give us the desired estimate. The fact that we have $\frac{1}{r'}$ instead of $1$ leads to the adjustment of our decaying order, as compared to (i).

The proof for the endpoint case is entirely similar.
\endproof

In the above estimates, notice that (modulo a positive power) the implicit constant for $t^{-(k-1+\frac{1}{p})}$ is proportional to the size of the given neighborhood.

Recall the notation $k_\theta := \max \{0, k_1,\dots, k_N\}$ and $H_p(w) := \|C_+(w^-)\|_p + \|C_-(w^+)\|_p$. The following corollaries are consequences of Lemma~\ref{localprpl} and Lemma~\ref{mainlemma} for $(m,k) = (0,1), (\frac{1}{2},1)$.

\begin{corollary} [Rough estimate] \label{Hardyrough}  Suppose $w^\pm \in H^{1,0}$ and both $\Omega_1,\Omega_2$ have $A_1+B_1$ decompositions, and $\theta_1,\theta_2$ satisfy the two phase conditions for $k=1$. If $w^\pm$ have sufficient decay at $\infty$, then for $2\leq p <\infty$
\begin{eqnarray*}
H_p(w\Omega) &\lesssim_p&   t^{-\frac{1}{p(k_\theta+1)}}
\end{eqnarray*}
for the pair of weights $w\Omega := \Big(\chi_{\{\theta'_2 \leq 0\}} w^- \Omega_2 e^{it\theta_2}, \chi_{\{\theta'_1 \geq 0\}} w^+\Omega_1 e^{it\theta_1}  \Big)$. If $w^\pm$ are compactly supported, the implicit constant is $O_p\Big(\|w^+\|_{H^{1,0}} + \|w^-\|_{H^{1,0}}\Big)$
\end{corollary}

\begin{corollary} [Vanish at stationary points] \label{Hardyvanish} Let $f \in H^{1,0}$ be supported near a stationary point $\lambda_0$ of order $k_0$, where $\Omega$ has an $A_1+B_1$ decomposition. If $f(\lambda_0)=0$ then for $2\leq p < \infty$
\begin{eqnarray*}
\|C_-(\chi_{\{\theta' \geq 0\}}  f\Omega e^{it\theta})\|_p + \|C_+(\chi_{\{\theta' \leq 0\}} f \Omega e^{it\theta})\|_p &\lesssim_p&
    \max(t^{-\frac{1}{p}}, t^{-(\frac{1}{2}+\frac{1}{p})\frac{1}{k_0+1}+\epsilon})
\end{eqnarray*}
\end{corollary}

Taking linear combination of Corollary~\ref{Hardyrough}, we see that generally $H_p(w)$ decays as $t\to\infty$ if there is a correct phase-weight relation in $w=(w^-,w^+)$ (now the phase is allowed to be piece-wisely defined). In that case, if $\|(1-C_w)^{-1}\|_{L^p\to L^p} \lesssim 1$ then the solution $\mu$ to $\mu = I + C_w \mu$ satisfies:
$$\|\mu-I\|_p = \|(1-C_w)^{-1}C_wI\|_p \lesssim H_p(w) \lesssim t^{-\frac{1}{p(k_\theta+1)}}.$$

\begin{corollary} [Almost orthogonality] \label{almostorthogonal}
Assume that $\Omega^\pm_1, \Omega^\pm_2$ have $A_1+B_1$ decompositions and $\theta_1,\theta_2$ satisfy the two phase conditions for $k=1$. For each $j\in \{1,2\}$, consider the restriction to $\{\theta'_j\geq 0\}$ of a pair of compactly supported $H^{1,0}$ weights
$$w_j = \big(\chi_{\{\theta'_j\geq 0\}}w_j^-e^{-it\theta_j}, \chi_{\{\theta'_j\geq 0\}}w_j^+e^{it\theta_j}\big).$$
If \emph{supp}$(w_1)$ and \emph{supp}$(w_2)$ are disjoint then $\forall 2\leq p <\infty$, the $L^p$ operators $C_{w_1\Omega_1}, C_{w_2\Omega_1}$ are almost orthogonal in the sense:
\begin{eqnarray*}
  \|C_{w_1\Omega_1} C_{w_2\Omega_2}\|_{L^p \to L^p}
  &\lesssim_p&     t^{-\frac{1}{p(k_{\theta_1}+1)}}\|w_1\|_{H^{1,0}}\|w_2\|_{H^{1,0}} \\
  \|C_{w_2\Omega_2} C_{w_1\Omega_1}\|_{L^p \to L^p}
  &\lesssim_p&     t^{-\frac{1}{p(k_{\theta_2}+1)}}\|w_1\|_{H^{1,0}}\|w_2\|_{H^{1,0}}
\end{eqnarray*}
for large $t$. The implicit constants depend on $\theta$ and the support of the weights.
\end{corollary}

\proof
By symmetry, it suffices to show the first estimate. Let $f \in L^p$. Since $w_1, w_2$ have disjoint supports, $\text{\emph{distance}}(\text{supp}(w_1), \text{supp}(w_2)) \gtrsim 1$ and $(C_{w_2}f)(\lambda)$ is analytic at every $\lambda \in w_1$. By H\"older's inequality, for every $n\geq 0$ and $\lambda \in \text{supp}(w_1)$:
\begin{eqnarray*}
  \frac{d^n}{d\lambda^n}(C_{w_2\Omega_2}f)(\lambda)    &\lesssim& \Big|\int f(x)\frac{\chi_{\{\theta'_2\geq 0\}}w^+_2(x)\Omega^+_2(x)e^{it\theta_2}}{(x-\lambda)^{n+1}} dx \Big|\\
   &+& \Big|\int f(x)\frac{\chi_{\{\theta'_2\geq 0\}}w^-_2(x)\Omega^-_2(x)e^{-it\theta_2}}{(x-\lambda)^{n+1}} dx\Big|\\
   &\lesssim& \|f\|_p \|w_2\|_\infty \Big(\int_{|x-\lambda|\gtrsim 1} \frac{1}{|x-\lambda|^{(n+1)p'}}dx\Big)^\frac{1}{p'}\\
   &\lesssim_{n,p}&  \|f\|_p \|w_2\|_\infty
\end{eqnarray*}
Now applying Corollary~\ref{Hardyrough}, we have
\begin{eqnarray*}
       \|C_{w_1\Omega_1}(C_{w_2\Omega_2}f)\|_p
       &\lesssim_p& t^{-\frac{1}{p(k_{\theta_1}+1)}} \Big(\big\|w^+_1 C_{w_2\Omega_2 }f\big\|_{H^{1,0}}+\big\| w^-_1C_{w_2\Omega_2}f\big\|_{H^{1,0}}\Big)\\
       &\lesssim_p& t^{-\frac{1}{p(k_{\theta_1}+1)}} \|w_2\|_\infty \|w_1\|_{H^{1,0}} \|f\|_p
\end{eqnarray*}
\endproof

Again, by taking linear combination, Corollary~\ref{almostorthogonal} remains true if there is a correct phase-weight relation in the pairs of weights $w_1$ and $w_2$.

The next lemma summarizes standard results in linear theory (see for instance Stein \cite{steinyd}) and is included here for convenience. In this lemma, $f$ have enough decay to ensure that relevant tail $L^1$ norms are finite.

\begin{lemma}\label{lineartheory} Let $f$ be $k$-time differentiable and supported near a stationary point $\lambda_0$ where it vanishes with multiplicity $m$ up to the $k^{\emph{th}}$ derivative. Let $\Omega$ have an $A_k+B_k$ decomposition on this neighborhood. Then
\begin{eqnarray}
\label{lineartheoryest}    \int f\Omega e^{it\theta} dx &\lesssim& \begin{cases} \max(t^{-k}, t^{-\frac{m+1}{k_0+1}}\ln t), & k\leq k_0\\
\max(t^{-k}, t^{-\frac{m+1}{k_0+1+\frac{1}{rk}}}\ln t), & k\geq k_0+1
\end{cases}
\end{eqnarray}
\end{lemma}

\begin{corollary} \label{farsupp}
With the same assumptions as in Lemma~\ref{lineartheory}, for any $\lambda_0 \in \mathbb C$ such that $\text{\emph{distance}}(\lambda_0,\text{supp}(f)) \gtrsim 1$ we can control $C(f\Omega e^{it\theta})(\lambda_0)$ by the same estimate as in (\ref{lineartheoryest}).
\end{corollary}
\noindent (write the left-hand side in integral form and apply the result of the lemma.)

\proof[Proof of Lemma~\ref{lineartheory}] Take a smooth cutoff $\varphi$ as before, we'll proceed as in the proof of Lemma~\ref{mainlemma}. By translation invariant we can assume that the given stationary point is $0$. Let $g(x) = \big(1-\varphi(\frac{x}{\epsilon})\big)f(x)\Omega(x)$. Using integration by parts, we have:
\begin{eqnarray*}
  \int_{\mathbb R} f\Omega e^{it\theta} dx
  &\lesssim& \epsilon\sup_{|x|\lesssim \epsilon}|f(x)\Omega(x)| + \int_{\mathbb R} |D^k_{t\theta'}g| dx\\
  &\lesssim&  \epsilon \sup_{|x|\lesssim \epsilon}|f(x)| + t^{-k}\|g_{\theta,k}\|_1
\end{eqnarray*}
Then rest of the argument is exactly the same as before, applied to $p=1$ (recall that the condition $p\geq 2$ was used in that proof only for Hausdorff-Young inequality). We then have the desired bound.
\endproof

Recall that in Lemma~\ref{direct_lem} we have two error terms of the form $\int \Delta w \; w_1$, which vanish if there is agreement in the triangularity structures of $w_1$ and $w_2$. The introduction of the scalar RHP in the next section might create some noise in their triangularity structure, and the next two lemmas will be used to handle these noise. These lemmas are complex variants of Lemma~\ref{localprpl} and Lemma~\ref{mainlemma}. Below, $\Gamma$ is a ray originating from one stationary point that forms a nontrivial angle with the real line, and $\Omega$ has an $A_k+B_k$ decomposition.

\begin{lemma} \label{variantLP} Let $1<p\leq \infty$. If $f$ has $k$ locally integrable derivatives and is supported away from the stationary points of $\theta$, and $f$ has sufficient decay then
$$\|C(fe^{it\theta})\|_{L^p(\Gamma)} \lesssim_p t^{-k}$$
\end{lemma}
\proof Without loss of generality, suppose that $\Gamma$ is originated from $0$. Let $g(x,\lambda) = \frac{f(x)}{x-\lambda}$. Using integration by parts we have
$$\int g e^{it\theta} \lesssim |\text{boundary terms}| + t^{-k} \int g_{\theta,k}$$
It is not hard to see that $g_{\theta,k}$ is controlled by a weighted sum of the derivatives of $f$, where for each $0\leq \beta \leq k$ $f^{(\beta)}$ has the following weight
$$\sum_{n+\beta = j \leq k} |x-\lambda|^{-(n+1)} w_{k,j}(\theta).$$
Consequently, by Minkowski's inequality we have
$$\Big\|\int g_{\theta,k}\Big\|_{L_\lambda^p(\Gamma)} \lesssim \sum_{\beta=0}^k \int  |f^{(\beta)}| \sum_{n+\beta = j \leq k}  |x|^{\frac{1}{p} - n - 1} w_{k,j}(\theta)$$
which is finite if $f$ has sufficient decay. Note that the condition $p>1$ ensures that the $L_\lambda^p(\Gamma)$ norms of $|x-\lambda|^{-(n+1)}$ are finite for $x\in \text{supp}(f)$. Also, similar to the proof of Lemma~\ref{localprpl}, the boundary terms of the integration by parts are either $0$ or can be controlled by the evaluations at $x_1,\dots, x_S$ of
$$\frac{1}{t|\theta'|} |D^{k-1}_{t\theta'}g| \lesssim t^{-k}\frac{1}{|\theta'|} g_{\theta,k-1}$$
Since $x_1,\dots, x_S$ are not stationary points, after taking $L^p(\Gamma)$ in $\lambda$ these boundary terms contribute $O(t^{-k})$.
\endproof

\begin{lemma} \label{variantMainLemma} Let $f$ has $k$ locally integrable derivatives and is supported near a stationary point $\lambda_0$ of order $k_0$. Let $f$ vanish at $\lambda_0$ with multiplicity $m$ up to the $k^{th}$ derivative.

(i) If $k\leq k_0$ then for $1<p< \infty$,
$$\|C(f\Omega e^{it\theta})\|_{L^p(\Gamma)} \;\;\lesssim_p\;\; \max(t^{-k},t^{-(m+\frac{1}{p})\frac{1}{k_0+1}}\ln t)$$

(ii) If $k\geq k_0+1$ and $\theta^{(k+1)}$ is assumed $L^r$ integrable near $\lambda_0$ then the estimate remains true if $k_0+1$ is replaced by $k_0+1+\frac{1}{rk}$.
\end{lemma}

\proof By symmetry can assume that $\lambda_0 = 0$. Take a normalized cutoff function $\varphi$ which is supported on $\{|x|\leq 2\}$ and equals to $1$ on $\{|x|\leq 1\}$ as usual. Then decompose $f = \varphi(\frac{.}{\epsilon}) f + (1-\varphi(\frac{.}{\epsilon})) f$ and estimate
$$\|C(\varphi(\frac{.}{\epsilon}) f\Omega e^{it\theta})\|_{L^p(\Gamma)} \lesssim \|\varphi(\frac{.}{\epsilon}) f\Omega e^{it\theta}\|_{L^p(\mathbb R)} \lesssim \epsilon^{\frac{1}{p}+m}$$
For $g:=(1-\varphi(\frac{.}{\epsilon})) f$ we proceed as in Lemma~\ref{variantLP} and get an estimate of
$$t^{-k} \sum_{n+\beta = j \leq k} \int |(\Omega f)^{(\beta)}| |x|^{\frac{1}{p} - n - 1} w_{k,j}(\theta)$$
As before, if $k\leq k_0$ then for every $\beta$, the weight of $(\Omega f)^{(\beta)}$ can be estimated by
$$|x|^{\frac{1}{p} - n - 1} |x|^{j-(k_0+1)k} = |x|^{\frac{1}{p}+\beta-1-(k_0+1)k}$$
(note that the estimate is independent of the choice of $(n,j)$ as long as they respect $n+\beta = j\leq k$). Thus, the same argument as before gives the following estimate for the contribution in $L^1$ of those terms
$$t^{-k}(1+ \epsilon^{\frac{1}{p}+m-(k_0+1)k})|\ln \epsilon|$$
and consequently we will have an overall estimate of
$$\epsilon^{\frac{1}{p} + m} + t^{-k}(1+ \epsilon^{\frac{1}{p}+m-(k_0+1)k}|\ln \epsilon|)$$
and the choice $\epsilon = t^{-\frac{1}{k_0+1}}$ gives the desired estimate.

When $k>k_0$, the estimate involving $\theta^{(k+1)}$ has to be done more carefully. This appears when $\beta = n = 0$ and the contribution of the respective term in the overall estimate can be estimated by
$$t^{-k}(1+\epsilon^{\frac{1}{p} + m -\frac{1}{r} - (k_0+1)k}|\ln \epsilon|) $$
so again we have an overall estimate of
$$\epsilon^{\frac{1}{p} + m} + t^{-k}(1+\epsilon^{\frac{1}{p} + m -\frac{1}{r} - (k_0+1)k}|\ln \epsilon|) $$
and the optimal choice is now $\epsilon = t^{-\frac{1}{k_0+1+\frac{1}{rk}}}$.
\endproof
Note: The endpoint estimate when $p=\infty$ in Lemma~\ref{variantMainLemma} can be formulated and proved similarly, and is left as an exercise. Lemma~\ref{variantLP} and \ref{variantMainLemma} remain true if $\Gamma \subset \mathbb R$ provided that there is a nontrivial distance from $supp(f)$ to $\Gamma$.

\section{Factorizations of the oscillatory jump matrix $J(\lambda,t)$} \label{factorizesection}

The operator formulation of Beals and Coifman indicates that a good factorization of $J$ is crucial to the study of our oscillatory RHP. Intuitively, a good factorization should separate the two terms $e^{\pm it\theta}$, since they decay in opposite regions. From Lemma~\ref{localprpl}, we want in our factorization the following behavior:
\begin{itemize}
\item The oscillating phase in $w^-$ is locally decreasing;
\item The oscillating phase in $w^+$ is locally increasing.
\end{itemize}
In our applications, the phase is either $\theta$ or $-\theta$ and could be piecewise-defined.

The jump matrix (\ref{oscjumpmatrix}) has the following canonical factorization:
\begin{eqnarray*}
J &=& \begin{pmatrix}1 & pe^{-it\theta}\cr 0 & 1\end{pmatrix}\begin{pmatrix}1 & 0 \cr qe^{it\theta} & 1\end{pmatrix} \equiv (I-w^-)^{-1}(I+w^+)
\end{eqnarray*}

This factorization has the correct weight-phase behavior on $\{\theta' \geq 0\}$ only, thus extra work is needed when both $\{\theta'>0\}$ and $\{\theta'<0\}$ are nonempty. Below, we will discuss an important scalar RHP which will be used to fix this problem.

Let $D_- = \{\theta' <0\}, D_+ = \{\theta' > 0\}$ and consider the scalar $L^2$-normalized RHP with the following jump matrix:
\begin{eqnarray}
\label{scalarRHP}    (1+pq)1_{D_-} + 1_{D_+}
\end{eqnarray}
The existence, boundedness and several important properties of the solution $\delta$ of this RHP will be proved in this section. In particular, its behavior near the stationary points will be studied. This $\delta$ was introduced by Deift and Zhou \cite{deiftzhouMKdVyd} in their studies of mKdV and was also used in \cite{varzuginyd}.

We describe how $\delta$ can be used to \emph{conjugate} the jump matrix $J(x,t)$. Indeed, if $M$ solves (\ref{oscRHP}) then $M^\delta:=M\delta^{-\sigma_3} \equiv M\begin{pmatrix}1/\delta & 0\cr 0 & \delta\end{pmatrix}$ solves the normalized $L^2$ RHP with the following jump matrix (it will be shown that $\delta$ is bounded away from $0$).
\begin{eqnarray}
\nonumber
        J_{conj}(x,t)
                 &=& \begin{pmatrix} \delta_-\delta_+^{-1}(1+pq) & \delta_-\delta_+ pe^{-it\theta} \cr \delta_-^{-1}\delta_+^{-1}qe^{it\theta} & \delta_-^{-1}\delta_+ \end{pmatrix}\\
\label{conjjumpmatrix}
                 &=&     \begin{cases}
                        \begin{pmatrix} 1 & \delta_-\delta_+ pe^{-it\theta} \cr 0 & 1 \end{pmatrix} \begin{pmatrix} 1 & 0 \cr \delta_-^{-1}\delta_+^{-1}qe^{it\theta} & 1   \end{pmatrix} & \text{if $x \in D_+$;}\\
                        \begin{pmatrix} 1 &  0\cr \delta_-^{-1}\delta_+^{-1}qe^{it\theta} & 1 \end{pmatrix} \begin{pmatrix} 1 & \delta_-\delta_+ pe^{-it\theta} \cr 0 & 1   \end{pmatrix} & \text{if $x \in D_-$}
                    \end{cases}
\end{eqnarray}
The factorization (\ref{conjjumpmatrix}) has the desired phase-weight relation. It is not hard to see that the conjugation doesn't change the unique solvability or the recovered potentials $u(t),v(t)$. Furthermore, the resolvent operators $(1-C_w)^{-1}$ associated with these RHPs have comparable norms on $L^p$. The last claim is a consequence of Lemma~\ref{boundedfactorize} and the boundedness of $\delta_\pm$. As we'll see, the above scalar RHP is the source of the $\ln t$ term and the interactions of stationary points in the leading asymptotics of $u(t)$, $v(t)$ stated in Theorem~\ref{maintheorem}.

\subsection{Existence and boundedness}
\begin{proposition} \label{scalarexistence} If $\ln(1+pq) \in L^2(\mathbb R)\cap L^\infty(\mathbb R)$ is real-valued then the $L^2$-normalized scalar RHP associated with the jump matrix (\ref{scalarRHP}) has unique solution, given by:
\begin{eqnarray*}
\label{delta}    \delta_\pm(\lambda) = \exp{C_\pm (1_{D_-}\ln(1+pq))(\lambda)}
\end{eqnarray*}
This solution is bounded, indeed for every $\lambda \in \mathbb C$
\begin{eqnarray*}
  \ln |\delta(\lambda)| \lesssim  \|\ln(1+pq)\|_\infty
\end{eqnarray*}
furthermore
\begin{eqnarray*}
|\delta_\pm(x)| = \begin{cases}(1+p(x)q(x))^{\pm\frac{1}{2}}, & x \in D_-\\
                                    1, & x \in D_+\end{cases}
\end{eqnarray*}
\end{proposition}

\proof Let $\delta$ be as in the conclusion of the proposition. The jump relation is automatic. We only need to show that $\delta_\pm -1 $ belongs to the respective Hardy spaces $H_2$ of the upper and lower half planes (by symmetry only need to show this for $\delta_+$). This will then easily imply the uniqueness of the solution.

For convenience, denote $g=\ln(1+pq)1_{D_-}$.

Let $\lambda$ be in the upper half plane, write $\lambda = x + iy$ with $y >0$. Since $g$ is real-valued,
\begin{eqnarray*}
\text{Re}(C_+g)(\lambda) = \frac{1}{2\pi}\int_{\mathbb R} \frac{g(u)y}{(u-x)^2 + y^2} du \lesssim  \|g\|_\infty
\end{eqnarray*}
In particular, this implies that $\delta_+$ is bounded:
$$\ln|\delta_+(z)| = \text{Re}(C_+g)(\lambda) \lesssim \|g\|_\infty$$

Now, using the elementary inequality $|e^{z} - 1| \lesssim |z| e^{|\text{Re}(z)|}$ and the above estimate for $\text{Re}(C_+g)(\lambda)$ , we easily have:
\begin{eqnarray*}
\|\delta_+(x+iy) - 1\|_2  &\lesssim& \  \|C_+g(.+iy)\|_2 \ e^{c \|g\|_\infty}\\
                             &\lesssim& \|g\|_2 \ e^{c \|g\|_\infty}
\end{eqnarray*}
for some absolute constant $c>0$. Thus, $\delta_\pm(\lambda)-1$ are bounded and belong to the respective Hardy spaces. Since $g$ is real-valued, $\delta(z)\overline{\delta(\overline{z})}=1$ for any $z \in \mathbb {C}\backslash D_-$. Thus, $|\delta_+(x)\delta_-(x)| = 1$ on $\mathbb R$. This plus the jump relation show the desired equalities for $|\delta_\pm(x)|$.

Below we show uniqueness. Suppose $\delta_1(z)$ is another normalized solution to the scalar RHP (\ref{scalarRHP}). Notice that if we replace $g$ by $-g$ in the above calculation, we would get $\delta^{-1}$ instead of $\delta$. So $\delta_\pm^{-1}-1$ are bounded and belong to the respective Hardy spaces. Now the function $m(\lambda) :=\delta_1(\lambda) \delta^{-1}(\lambda)$ satisfies the jump relation:
$$m_+(x) = m_-(x), \;\; x \in \mathbb R$$

Using boundedness of $\delta^{-1}$ and $\delta_1$,  it is not hard to check that $m_\pm(\lambda)$ belongs to the respective Hardy spaces. From there it is clear that $m(\lambda) = 1$.
\endproof

\subsection{Bounds on derivatives of $\delta$} \label{derivativesbound}

We will prove some results about the behavior of $\delta_\pm$ on $\mathbb R$. Part (b) of the next lemma generalizes a result of Deift and Zhou \cite{deiftzhouNLSyd}, where it was proved for $\theta(x)=x^2$.

From now on we will define $\Omega_0(x) := \delta_+(x) \delta_-(x)$.

\begin{lemma}[A-B Decomposition Lemma]\label{ABdecomp} Assume that $\ln(1+pq) \in H^{k,0}$. Then,

(a) $\Omega_0^{\pm 1}$ are $k$-time differentiable inside $\mathbb R\setminus \partial \overline{D_-}$, and for $1\leq n\leq k-1$:
\begin{eqnarray}
\label{deltaest}    \frac{d^n}{dx^n}(\Omega_0^{\pm 1})(x) &\lesssim& 1 +\frac{1}{\text{\emph{distance}}(x, \partial \overline{D_-})^n}
\end{eqnarray}

(b) $\frac{d^k}{dx^k}(\Omega_0^{\pm 1})(x)$ can be decomposed into $A_k(x) +B_k(x)$, with
\begin{eqnarray*}
    \|A_k\|_2 \lesssim_k 1, &&  |B_k(x)| \lesssim_k 1+\frac{1}{\text{\emph{distance}}(x, \partial \overline{D_-})^k}
\end{eqnarray*}
\end{lemma}

\noindent \emph{Remarks:} 1. If $\theta$ doesn't have any stationary points then $D_- = \emptyset$ or $\mathbb R$. In these cases, our convention in the lemma is $\partial \overline{D_-}=\emptyset$ and $\frac{1}{\text{\emph{distance}}(x, \partial \overline{D_-})}=0$.

2. If $\ln(1+pq)$ is only $H^{k,0}$ on a neighborhood of $x_0$ (still $L^2$ elsewhere), the above properties remain valid locally. More precisely, we can use a smooth cutoff centered at $x_0$ to decompose $\ln(1+pq)$ into two components, one with high regularity and supported near $x_0$, and another with low regularity but supported away from $x_0$. Then notice that the contribution of the low regularity part is bounded if $x \in \mathbb R\setminus \partial \overline{D_-}$ is near $x_0$.

\proof

(a) Let $g=\ln(1+pq)1_{D_-}$. Since $\Omega_0^{\pm 1}(x) = e^{\mp (Hg)(x)}$ are bounded, to show (\ref{deltaest}) it suffices to show that $(Hg)(x)$ is $k$-time differentiable on $\mathbb R\setminus \partial \overline{D_-}$, with
\begin{eqnarray}
\label{HTransformEst}   \frac{d^n}{dx^n}(Hg)(x) &\lesssim& 1 +\frac{1}{\text{\emph{distance}}(x, \partial \overline{D_-})^n}, \; \forall 1\leq n \leq k-1
\end{eqnarray}
If $\theta$ doesn't have any stationary points then $g\in H^{k,0}$ (and hence $Hg \in H^{k,0}$), and (\ref{HTransformEst}) follows from a standard Sobolev estimate:
$$\|f^{(k)}\|_\infty \lesssim \|f\|_{H^{k+1,0}}$$
In the rest of the argument, we'll assume that $\theta$ has $N\geq 1$ stationary points.

Below, we'll assume that $1\leq n\leq k$ unless otherwise specified.

Let $\varphi_0(y)$ be a nice $C^\infty$ function supported in $[-\frac{1}{2},\frac{1}{2}]$ and equal to $1$ in $[-\frac{1}{4},\frac{1}{4}]$. Take $M>0$ be a small number and consider $x\not\in \bigcup_{\alpha \in \partial \overline{D_-}} [\alpha-M,\alpha+M]$.

Define $\varphi_1(y), \; \varphi_2(y)$ by:
$$\varphi_1(y) = \sum_{\alpha \in \partial \overline{D_-}} \varphi_0(\frac{y-\alpha}{M}),\;\varphi_2 = (1-\varphi_1)1_{D_-}$$
Notice that $\varphi_2$ is zero near the endpoints of $D_-$, thus it is also $C^\infty$. Decompose $\frac{d^n}{dx^n}(H g)(x)$ into:
\begin{eqnarray}
\label{Hilbertdecomp}\text{I}(x) + \text{II}(x) \equiv \frac{d^n}{dx^n} H(\ln(1+pq)\varphi_1 1_{D_-})(x) +\frac{d^n}{dx^n}H\Big(\ln(1+pq)\varphi_2\Big)(x)
\end{eqnarray}
{\bf Estimation of $\text{II}$:} Since $\varphi_2 \ln(1+pq) \in H^{k,0}$, for all $n\leq k$ we have
\begin{eqnarray*}
\text{II}(x) &=& H\Big([\ln(1+pq)\varphi_2]^{(n)}\Big)(x)\\
             &=& H\Big(\ln(1+pq)^{(n)}\varphi_2\Big)(x) + \sum_{j=1}^n \binom{n}{j} H\Big(\ln(1+pq)^{(n-j)}\varphi_2^{(j)}\Big)(x)\\
             &\equiv& \text{II}_1(x) + \text{II}_2(x)
\end{eqnarray*}
\underline{For $\text{II}_2(x)$:}  Notice that for any $j\geq 1$, we have $\varphi^{(j)}_2 \lesssim \frac{1}{M^j}$ and $\text{supp}(\varphi^{(j)}_2) \subset \{y: \text{\emph{distance}}(y,\partial \overline{D_-}) \leq M/2\}$. This means $x$ is outside the support of $\varphi^{(j)}_2$, furthermore
$$\text{\emph{distance}}(x,\text{supp}(\varphi^{(j)}_2)) \gtrsim M.$$
Thus, we can write every $H\Big(\ln(1+pq)^{(n-j)}\varphi_2^{(j)}\Big)(x)$ in the integral form. Since $\|\ln(1+pq)^{(n-j)}\|_\infty \lesssim \|\ln(1+pq)\|_{H^{n,0}}$, $\text{II}_2(x)$ can be controlled by
\begin{eqnarray*}
\text{II}_2(x)
&\lesssim&\|\ln(1+pq)\|_{H^{n,0}}\sum_{j=1}^n \frac{1}{M}\int |\varphi_2^{(j)}(y)|dy\\
&\lesssim&   \|\ln(1+pq)\|_{H^{n,0}}\Big(\frac{1}{M}+\dots + \frac{1}{M^n}\Big)
\end{eqnarray*}
\underline{For $\text{II}_1(x)$:} Notice that $\|\text{II}_1(x)\|_2 \lesssim \|\ln(1+pq)\|_{H^{n,0}}$. Thus, even when $n=k$ this term exists a.e. When $n<k$ we can estimate it pointwise as follows:
\begin{eqnarray*}
   \text{II}_1(x)  &\lesssim& \|\ln(1+pq)^{(n)}\varphi_2\|_{H^{1,0}}\\
             &\lesssim& \|\ln(1+pq)\|_{H^{n+1,0}} \Big(1+ \frac{1}{M}\Big)
\end{eqnarray*}
{\bf Estimation of $\text{I}$:} Notice that $x\not\in \text{supp}(\varphi_1)$, thus
\begin{eqnarray*}
  \text{I}(x)
  &=& \frac{1}{\pi i}\int_{D_-} \frac{(-1)^n n! \ln(1+pq)(y)\varphi_1(y)}{(y-x)^{n+1}} dy\\
  &=& \text{boundary terms} + \frac{1}{\pi i}\int_{D_-} \frac{\big[\ln(1+pq)\varphi_1\big]^{(n)}(y)}{y-x} dy
\end{eqnarray*}
using repeated integration by parts. Here, the boundary terms are of the forms
\begin{eqnarray*}
C_{i,n} \frac{\big[\ln(1+pq)\varphi_1\big]^{(i)}(\alpha)}{(x-\alpha)^{n-i}},  \;0\leq i\leq n-1,\; \alpha \in \partial\overline{D_-}
\end{eqnarray*}
Since $\varphi_1(\alpha)=1$ and $\varphi^{(j)}_1(\alpha)=0$ for any $j \geq 1, \;\alpha \in \partial\overline{D_-}$, these boundary terms can be controlled by $$\|\ln(1+pq)\|_{H^{n,0}}(\frac{1}{M}+\dots+\frac{1}{M^n}).$$
On the other hand, the remaining term $\int_{D_-} \big[\ln(1+pq)\varphi_1\big]^{(n)}(y)/(y-x) dy$ is easily controlled in $L^2$ by $\|\ln(1+pq)\|_{H^{n,0}}$, even for $n=k$. When $n<k$, it can be estimated a. e. by:
\begin{eqnarray*}
       &&  \sum_{j=0}^n \int_{D_-} \frac{|\ln(1+pq)^{(j)}(y)| |\varphi_1^{(n-j)}(y)|}{|y-x|} dy\\
       &\lesssim& \sum_{j=0}^n\|\ln(1+pq)\|_{H^{j+1,0}}  \frac{1}{M}\int_{D_-} |\varphi_1^{(n-j)}(y)| dy\\
       &\lesssim& \|\ln(1+pq)\|_{H^{n+1,0}} (1 + \frac{1}{M} + \dots + \frac{1}{M^n})
\end{eqnarray*}
Consequently, $\frac{d^n}{dx^n}(H g)(x)$ exists for every $n\leq k$, and if $1\leq n\leq k-1$ then:
$$\frac{d^n}{dx^n}(H g)(x) \lesssim \|\ln(1+pq)\|_{H^{k,0}} (1+\frac{1}{M^n})$$
for $x\not\in \bigcup_{\alpha \in \partial \overline{D_-}} [\alpha-M,\alpha+M]$. Taking $M=\text{\emph{distance}}(x,\partial \overline{D_-})/2$ we get (\ref{HTransformEst}).

(b) We can use the chain rule to compute $\frac{d^k}{dx^k}(\Omega_0^{\pm 1})(x)$ as
\begin{eqnarray*}
\frac{d^k}{dx^k}(\Omega_0^{\pm 1})(x)
&=& \mp\Omega_0(x) \Big(\frac{d^k}{dx^k}(H g)(x) + \text{remaining terms}\Big)\\
&\equiv& A_k + B_k
\end{eqnarray*}
By a symmetry trick $x \mapsto cx$, it is not hard to see that the remaining terms on the right-hand side are linear combinations of those $(Hg)^{(i_1)}(x) \dots (Hg)^{(i_j)}(x)$ with $i_1 + \dots + i_j = k$ and $1\leq i_1, \dots, i_j \leq k-1$. Therefore using part (a) we see that these terms are bounded by
$$1 +\frac{1}{\text{\emph{distance}}(x, \partial \overline{D_-})^k}$$
This gives the desired bounds for $A_k,B_k$.
\endproof

\subsection{Approximation of $\delta$} \label{behaviorsect} Compared to the last section, the added assumption in this section is $p,q\in H^{1,0}$, which makes $\ln(1+pq)\in H^{1,0}(\mathbb R)$. The goal of this section is to find an approximation for $\delta$ at a stationary point $\lambda_j$ of $\theta$. Our approximation will satisfy the following {\it model} scalar RHP:
$$\delta_{j+}(x) = \delta_{j-}(x) \Big(1_{D_{j+}}(x) + 1_{D_{j-}}(x)[1+p(\lambda_j)q(\lambda_j)]\Big),\; x \in \mathbb R,$$
$$\text{ where } D_{j\pm} = \{x\in \mathbb R: \pm \theta'_j(x)> 0\}$$
Here $\theta_j(x) := \theta(\lambda_j) + \frac{\theta^{(k_j+1)}(\lambda_j)}{(k_j+1)!}(x-\lambda_j)^{k_j+1}$ is the Taylor approximation of $\theta(x)$ at $\lambda_j$.

Notice that we don't impose any normalization condition for $\delta_j$, only analyticity on $\mathbb C\setminus \mathbb R$ is required. Intuitively, the above jump matrix can be seen as a limiting approximation of the jump matrix of $\delta$, which explains why such $\delta_j$ plays an important role in modeling RHPs localized to very small neighborhood of  $\lambda_j$ (where $p(x),q(x)$ are essentially $p(\lambda_j),q(\lambda_j)$ and $\theta$ is essentially $\theta_j$). Note that $D_{j\pm}$ can be seen as the local approximation near $\lambda_j$ of $D_\pm = \{\pm   \theta'>0\}$.

A consequence of the result in this section is the existence of such $\delta_j$ so that along any ray $\gamma$ originated from $\lambda_j$ that forms non trivial angle with $\mathbb R$, we'll have
$$|\delta(z) - \delta_j(z)| \lesssim |z-\lambda_j|^\frac{1}{2}, \;\; z\in\gamma$$
The upper bound can be improved  to $|z-\lambda_j|^{1-\epsilon}$ if $\ln(1+pq)\in H^{2,0}$, and the argument is similar with minor adaptations.

The approximation of $\delta$ is based on the approximation of its logarithm, $C(1_{D-}\ln(1+pq))$. Essentially, if $h$ approximates this logarithm and $Re(h)$ is bounded then by the mean value theorem
$$\delta - \exp(h) \lesssim |C(1_{D-}\ln(1+pq)) - h|.$$

For simplicity of notation, we'll denote $p(\lambda_j)$ and $q(\lambda_j)$ by $p_j$ and $q_j$. Consider the following triangular function (supported on $[-1,1]$):
\begin{eqnarray*}
T(x) = \begin{cases}
             0, & \text{ if $|x|\geq 1$;}\\
             x+1, & \text{ if $-1\leq x \leq 0$;}\\
             -x+1, & \text{ if $0\leq x\leq 1$.}
         \end{cases}
\end{eqnarray*}
(note that $T \in H^{1,0}$, this is also the reason for choosing it triangular). Let $c\gtrsim 1$ be such that $c<\min_{i\neq j} |\lambda_i -\lambda_j|$ and let $T_c(x) := T((x-\lambda_j)/c)$.

To get a sense of the desired approximation, we'll compute $C(1_{D_-}\ln(1+p_jq_j)T_c)(z)$ for $z\in \mathbb C\setminus \mathbb R$. Thanks to the constraint on $c$, we can write it as
\begin{eqnarray*}
  \frac{\ln(1+p_jq_j)}{2\pi i}\int_{D_{j-}} \frac{T_c(x)}{x-z} dx
\end{eqnarray*}
To compute this integral explicitly, we have to examine $D_{j-}$, or equivalently the relative position of $\lambda_j$ within $\overline{D_-}$. Let $\epsilon$ be a function on $\mathbb R$ that assigns $1$ to every right endpoint, $-1$ to every left endpoint, and $0$ to every interior/exterior point of $\overline{D_-}$ (note that endpoints of $\overline{D_-}$ are stationary points). Note that for every $j$, $\epsilon(\lambda_j)$ depends only on the parity of $k_j$ and the sign of $\theta^{(k_j+1)}(\lambda_j)$:
\begin{eqnarray*}
  \epsilon(\lambda_j)
  &=&  \begin{cases}
         0, &\text{ if $k_j$ is even;}\\
         \text{sgn}(\theta^{(k_j+1)}(\lambda_j)), &\text{ if $k_j$ is odd.}
       \end{cases}\\
  &\equiv& \epsilon_j
\end{eqnarray*}
Recall that $\nu_j := -\frac{1}{2\pi}\ln[1+p_jq_j]$. Consider the first case when $\epsilon(\lambda_j) = 1$. Then the above integral can be rewritten as
\begin{eqnarray*}
  &&i\nu_j\int^{\lambda_j}_{\lambda_j-c} \frac{x-\lambda_j+c}{c(x-z)} dx \\
  &=& i\nu_j \;\frac{c+ (z-\lambda_j+c)[\ln(z-\lambda_j) -  \ln(z-\lambda_j+c)]}{c}\\
  &=& i\nu_j \;\ln(z-\lambda_j) + i\nu_j \; g_c(z)
\end{eqnarray*}
where $g_c(z) := 1 + \Big((z-\lambda_j)\ln(z-\lambda_j) -  (z-\lambda_j+c)\ln(z-\lambda_j+c)\Big)/c$. It is not hard to see that
$$\|g'_c(z)\|_{L^2(\gamma)} \lesssim 1$$
where the implicit constant depends on $c$ and the angle between $\gamma$ and $\mathbb R$ (indeed near $\lambda_j$, $g'_c$ is $\lesssim |\ln(z-\lambda_j)|$ which is $L^p$-integrable there, while for large $|z-\lambda_j|$ we can estimate $g'_c$ by $\frac{1}{|z-\lambda_j|}$). Consequently by Cauchy-Schwarz's inequality,
$$|g(z_1)-g(z_2)| \lesssim |z_1-z_2|^{1/2}, \; \forall z_1,z_2\in\gamma$$
so $\exists \lim_{z\to\lambda_j, z\in \gamma}g(z)$, which indeed exists nontangentially (i.e. independent of $\gamma$):
$$\lim_{z\to\lambda_j}g(z) = 1 -\ln c$$
On the other hand, it is not hard to see that $\exp(i\nu_j \;\ln(z-\lambda_j))$ satisfies the desired model jump relation. Thus, the desired approximation for $\exp\Big(C(1_{D_-}\ln(1+p_jq_j)T_c)\Big)$ in this case is:
$$\exp\Big(i\nu_j \;\ln(z-\lambda_j) + i\alpha\Big)$$
for some real number $\alpha$ that can be computed by the following nontangential limit:
$$\alpha = \frac{1}{i}\lim_{z\to\lambda_j} \Big(C[1_{D_-}\ln(1+p_jq_j)T_c](z) - i\nu_j \ln(z-\lambda_j)\Big)$$

When $\epsilon(\lambda_j) = -1$, the above argument essentially works, except for one place: in the computation, we need to choose $\ln(x-z)$ as the anti-derivative of $\frac{1}{x-z}$, instead of $\ln(z-x)$. This will ensure that the respective $g_c(z)$ (which is now $1 + \Big((z-\lambda_j)\ln(z-\lambda_j) -  (z-\lambda_j+c)\ln(z-\lambda_j+c)\Big)/c$) has nontangential limit as $z\to\lambda_j$ independent of the approaching direction (more precisely, of $\text{Im}(z-\lambda_j)$). Consequently we will have a slightly different approximation
$$\exp\Big(-i\nu_j \;\ln(\lambda_j-z) + i\alpha\Big), \;\;\text{ $\alpha$ is defined by a similar limit.}$$

When $\epsilon(\lambda_j)=0$ (i.e. $k_j$ is even), there are two possible scenarios: $D_{j-} =\emptyset$ (which corresponds to $\theta^{(k_j+1)}(\lambda_j)>0$ or equivalently $\lambda_j$ is an exterior point of $\overline{D_-}$), or $D_{j-}=\mathbb R$ (the opposite case). In the former case, clearly $C(1_{D_-}\ln(1+p_jq_j)T_c)(z) \equiv 0$. In the latter case, we can ``sum" the above approximations to get the desired approximation:
$$\exp\Big(i\nu_j \;\ln(z-\lambda_j)-i\nu_j \;\ln(\lambda_j-z) + i\alpha\Big), \;\;\text{ $\alpha$ is defined by a similar limit.}$$
This can be further simplified using $\ln(z-\lambda_j)- \ln(\lambda_j-z)= \pi i\text{sgn}[\text{Im}(z)]$. Consequently, the desired approximation of $\exp\Big(C[1_{D_-}\ln(1+p_jq_j)T_c]\Big)$ is of the form $\exp(i\alpha + \beta_j(z))$, where
\begin{eqnarray*}
\beta_j(\lambda) &=& \begin{cases}
0, &\text{if $\lambda_j$ is an exterior point of $\overline{D_-}$;}\\
i\epsilon(\lambda_j)\nu_j\ln\big[\epsilon(\lambda_j)(\lambda-\lambda_j)\big], &\text{if $\lambda_j$ is an endpoint of $\overline{D_-}$;}\\
-\pi\nu_j\text{sgn}[\text{Im}(\lambda)]], &\text{if $\lambda_j$ is an interior point of $\overline{D_-}$.}
\end{cases}\\
\alpha &=& \frac{1}{i}\lim_{z\to\lambda_j} \Big(C[1_{D_-}\ln(1+p_jq_j)T_c](z) - \beta_j(z)\Big)
\end{eqnarray*}

In the next proposition, we'll show that the desired approximation of $C(1_{D_-}\ln(1+pq))(z)$ is of similar form, $$\delta_j(z) := \exp(i\omega_j + \beta_j(z)),$$
where $\omega_j$ is defined by the following nontangential limit:
\begin{eqnarray}
\label{omegajdefinition} \omega_j &:=& \frac{1}{i}\lim_{z\to\lambda_j} \Big(C[1_{D_-}\ln(1+pq)](z) - \beta_j(z)\Big)
\end{eqnarray}
\noindent \underline{\it Remarks:} Clearly, if the limit exists then $\omega_j\in\mathbb R$. From there we can easily see that $\delta_j$ is bounded on $\mathbb C\setminus \mathbb R$. Indeed, the only nontrivial case is when
$\lambda_j$ is an endpoint of $\overline{D_-}$. In that case, $\forall \lambda \in \mathbb C$,
$$|\delta_j(\lambda)| = \exp\Big(-\epsilon(\lambda_j)\nu_j\arg[\epsilon(\lambda_j)(\lambda-\lambda_j)]\Big) \lesssim \exp(\|\ln(1+pq)\|_\infty)$$

\begin{proposition} \label{behaviorprop} Assume $p,q,\ln(1+pq)\in H^{1,0}$ and $\ln(1+pq)$ is real valued. Then:

(i) The limit defining $\omega_j$ exists nontangentially.

(ii) $\forall$ ray $\gamma$ originated at $\lambda_j$ that forms an angle of measure $\gtrsim 1$ with $\mathbb R$:
$$|\delta(\lambda) - \delta_j(\lambda)|   \lesssim |\lambda-\lambda_j|^\frac{1}{2}$$
\end{proposition}
\proof Let $h=\ln(1+pq) - T_c$. It is good enough to show that $\lim_{z\to\lambda_j} C(1_{D_-}h)(z)$ exists nontangentially, and on any such $\gamma$, $C(1_{D_-}h)(z)$ is H\"older continuous of exponent $1/2$. Indeed, we'll show that these claims are true for every $h \in H^{1,0}(\mathbb R)$ such that $h(\lambda_j)=0$.

Decompose $1_{D_-}h$ as $h_1 + h_2 \equiv 1_{D_{j-}}h + 1_{D_-\setminus D_{j-}}h$ and notice that $h_1$ is in $H^{1,0}$ and vanishes at $\lambda_j$, while $h_2$ is in $L^2 \cap L^\infty$ and supported away from $\lambda_j$. We'll show the claims for $h_1$ and $h_2$ using these properties.

First, since $h_1 \in H^{1,0}$ (so is continuous), the nontangential limits $\lim_{z\to\lambda_j}(C_\pm h_1)(z)$ exist and equal to $(C_\pm h_1)(\lambda_j)$. Since $(C_+h_1)(\lambda_j) - (C_-h_1)(\lambda_j) = h_1(\lambda_j)=0$, these limits are the same, so $\lim_{z\to\lambda_j}(C h_1)(z)$ exists as desired. Furthermore, by a simple application of Carleson's measure theorem, we have
$$\|(C h_1)'\|_{L^2(\gamma)}  =  \|(C h'_1)\|_{L^2(\gamma)} \lesssim \|h'_1\|_{L^2(\mathbb R)}$$
so by Cauchy-Schwarz $(C h_1)$ is H\"older continuous of exponent $1/2$ on $\gamma$.

For $h_2$, the fact that $h_2$ is supported away from $\lambda_j$ easily implies the existence of the limit, indeed $(C h_2)$ is now analytic at $\lambda_j$. Again, for the estimate it is good enough to show that $(Ch_2)' \in L^2$. Using the support condition it is not hard to see that for $z\in \gamma$,
\begin{eqnarray*}
(Ch_2)'(z) &=& \frac{1}{2\pi i} \int_{\mathbb R} \frac{h_2(x)}{(x-z)^2} dx\\
           &\lesssim& \|h_2\|_\infty \frac{1}{1 + |z-\lambda_j|}
\end{eqnarray*}
notice the implicit constant depends on the angle between $\gamma$ and $\mathbb R$. Consequently,
\begin{eqnarray*}
\|(Ch_2)'(z)\|_{L^2(\gamma)}  &\lesssim& \|h_2\|_\infty
\end{eqnarray*}
as desired.
\endproof
\noindent \underline{\it Remarks:} When $\ln(1+pq)\in H^{1,1}$, we can compute $\omega_j$ more explicitly by
\begin{eqnarray}
\label{omega_j} \frac{1}{2\pi} \int_{D_-} \ln|\lambda_j-y|d\ln[1+pq](y)+\sum_{1\leq k \leq N, \lambda_k \neq \lambda_j} \epsilon(\lambda_k)\nu(\lambda_k)\ln|\lambda_j-\lambda_k|
\end{eqnarray}
using suitable integration by parts. The condition $\ln(1+pq)\in H^{1,1}$ in particular ensures that the infinity boundary terms of the partial integration vanish. To see this it might be convenient to distinguish between values of $\epsilon(\lambda_j)$.

\section{Unique solvability of oscillatory RHPs}

\subsection{Unique solvability of the original RHP} In this section, we sketch the main ideas to show that for $t$ large enough, the RHP (\ref{oscRHP}, \ref{oscjumpmatrix}) is uniquely solvable. In light of the discussion in Section~\ref{opsection}, it suffices to show that for large $t$ there is a factorization of $J(\lambda,t)= (I-w^-(\lambda,t))^{-1}(I+w^+(\lambda,t))$ such that $I-C_w$ is invertible.

\begin{theorem} \label{uniquesolvability} Let $\theta $ satisfies (A) and (B) and $p$,$q$ are continuous and vanish at $\infty$, with $0<1+pq \lesssim 1$. Assume that $pq<1$ at stationary points of order $\geq 3$. Then there exists a factorization of $J(t)$ such that $\|(I-C_{w(t)})^{-1}\|_{L^r \to L^r} \lesssim 1$ uniformly as $t\to\infty$ for $r\geq 2$ sufficiently close to $2$. If $k_\theta \leq 2$ and $pq\geq 0$ at every secondary stationary point then we allow $2\leq r<\infty$.
\end{theorem}

Along the reductions used in the proof of Theorem~\ref{maintheorem}, Theorem~\ref{uniquesolvability} can be proved for a nicer class of $p,q$, say $p, q \in C^3_o$; furthermore for large $t$ the resolvent norms $\|(1-C_w)^{-1}\|$ are bounded by constants depending essentially on the values of $|p|$ and $|q|$ at stationary points.  To extend this to more general $p,q$, just approximate them by $p_1, q_1 \in C^3_o$ such that
\begin{itemize}
\item $0<1+p_1q_1 \lesssim 1$;
\item $(p,q)$ agrees with $(p_1,q_1)$ at every stationary point;
\item $p-p_1$ and $q-q_1$ are uniformly small on $\mathbb R$;
\end{itemize}
and invoke an application of Newman series.

Let $p,q\in C^3_o$ now. Conjugate $J$ using the solution $\delta$ to our scalar RHP as in the last section to get $J_{conj}$; this does not affect the conclusion of Theorem~\ref{uniquesolvability} thanks to Lemma~\ref{boundedfactorize}. After factorizing the jump matrix $J_{conj}(\lambda,t)$ nicely as in Section~\ref{factorizesection}, we use Lemma~\ref{localprpl} described in Section~\ref{localizationsection} to reduce the corresponding weights to a sufficiently small neighborhood of stationary points. In this reduction, we will show that for large $t$ and for any $2\leq p<\infty$:
$$\|(1-C_{w_{old}})^{-1}\|_{L^p \to L^p} \lesssim_p 1+\|(1- C_{w_{new}})^{-1}\|_{L^p \to L^p}$$

Near each stationary point $\lambda_0$, $\theta$ will be reduced to an analytic phase, after that the contribution of $\lambda_0$ is separated. Eventually we arrive at a number of simpler RHPs, one localized to a small neighborhood of each stationary point with a nice analytic phase. Each such RHP will be further reduced to a model RHP. The validity of Theorem~\ref{uniquesolvability} for the model RHPs will be shown separately in Section~\ref{modelsection}.

For the separation of contribution, let $w = w_0 + w_1$ where $w_0 = (w^-_0, w^+_0)$ is supported near $\lambda_0$ and $w_1=(w^-_1,w^+_1)$ is supported near the other stationary points. The main idea is to use the following parametrix, introduced by Varzugin \cite{varzuginyd}:
\begin{eqnarray*}
(I-C_w)(I + \sum_{j=0}^1 C_{w_j}(I-C_{w_j})^{-1}) &=& I - \sum_{0\leq
j\neq n \leq 1}
C_{w_j}C_{w_n} (I-C_{w_n})^{-1}\\
(I + \sum_{j=0}^1 C_{w_j}(I-C_{w_j})^{-1})(I-C_w) &=& I - \sum_{0\leq
j\neq n \leq 1}  (I-C_{w_n})^{-1}C_{w_j}C_{w_n}
\end{eqnarray*}
Thanks to the scalar RHP~\ref{scalarRHP} which places the oscillating terms $e^{\pm it\theta}$ in the right places, we can apply Corollary~\ref{almostorthogonal} to see that $\|C_{w_n}C_{w_j}\|_{L^p \to L^p}$'s decay as $t\to\infty$ (whenever $n\neq j$). This enables the success of the above parametrix.

\subsection{Unique solvability of intermediate RHPs}

In our reductions, we want to ensure that:

\begin{corollary}\label{coruniqsolv}
If $p,q \in H^{1,0}$ then for large $t$ the intermediate RHPs that appear as consequences of the reductions are uniquely solvable, with uniformly bounded resolvent norms (as $t\to\infty$) in relevant $L^p$.
\end{corollary}

\proof Thanks to the above parametrix, it suffices to show the resolvent bound for an RHP associated with a pair of weights $(w^-,w^+)$ localized to a neighborhood $P$ of one stationary point. Let $\phi$ be a normalized smooth cutoff function supported on $P$. Then we can assume
\begin{eqnarray*}
    (w^-,w^+) =     \begin{cases}
                        \Bigl(\begin{pmatrix} 0 & \phi \delta_-\delta_+ p e^{-it\theta} \cr 0 & 0 \end{pmatrix}, \begin{pmatrix} 0 & 0 \cr \phi \delta^{-1}_-\delta^{-1}_+ q e^{it\theta} & 0   \end{pmatrix}\Bigr), & \text{if $\lambda \in D_+$;}\\
                        \Bigl(\begin{pmatrix} 0 &  0\cr \phi \delta^{-1}_-\delta^{-1}_+ q e^{it\theta} & 0 \end{pmatrix}, \begin{pmatrix} 0 & \phi \delta_-\delta_+ p e^{-it\theta} \cr 0 & 0   \end{pmatrix}\Bigr), & \text{if $\lambda \in D_-$}
                    \end{cases}
\end{eqnarray*}

Now, approximate $p,q$ by $p_1,q_1 \in C^3_o$ such that $\|p-p_1\|_{H^{1,0}}$ and $\|q-q_1\|_{H^{1,0}}$ are small, furthermore $p$ and $p_1$ agree at the current stationary point, and so do $q$ and $q_1$. For $\delta = e^{C(1_{D_-}\ln(1+pq))}$, we'll approximate $\ln(1+pq)$ by $h \in C^3_o$ real valued, such that they agree at \emph{every} stationary point and $\|\ln(1+pq) - h\|_{H^{1,0}}$ is small. It is then not hard to see that $\delta_1 := e^{C(1_{D_-}h)}$ is close to $\delta$ in $L^\infty$. Consequently, the approximated weights $(W^-,W^+)$ (where $p,q,\delta$ are replaced by $p_1,q_1,\delta_1$) is close to $(w^-,w^+)$ in $L^\infty$.

Now, note that $\Omega_1 = \delta_{1-} \delta_{1+}$ and its inverse have $A_3+B_3$ decompositions on $P$. The same reductions as in the proof of Theorem~\ref{uniquesolvability} can be used to show that the RHP associated with $W$ satisfies the desired resolvent bound, and for large $t$ the bound depends essentially on the value of $|p_1|,|q_1|$ at the current stationary point. Consequently if the above approximation errors are small enough then we get a comparable resolvent bound for $(w^-,w^+)$.
\endproof

\section{Reduction to model cases (I): Localization and phase reduction} \label{reductmodelsection}
In this section and the following section, we will use the localization
schemes described in Section~\ref{localizationsection} to reduce the RHP (\ref{oscRHP},\ref{oscjumpmatrix}) to a number of model RHPs, at the expense of modifying the solution and the potentials $u(t), v(t)$ by terms having sufficient decay as $t\to\infty$.

As discussed in Section~\ref{factorizesection}, we first conjugate our RHP (\ref{oscRHP},\ref{oscjumpmatrix}) by the scalar RHP (\ref{scalarRHP}) and end up with the RHP $(M^\delta,J_{conj})$. For simplicity of notation, let $P=\Omega_0 p \equiv \delta_- \delta_+ p$, $Q=\Omega_0^{-1} q \equiv \delta_-^{-1} \delta_+^{-1} q$. From (\ref{conjjumpmatrix}), the weights for $J_{conj}$ are:
\begin{eqnarray*}
    (w^-,w^+) =     \begin{cases}
                        \Bigl(\begin{pmatrix} 0 & Pe^{-it\theta} \cr 0 & 0 \end{pmatrix}, \begin{pmatrix} 0 & 0 \cr Qe^{it\theta} & 0   \end{pmatrix}\Bigr), & \text{if $x \in D_+$;}\\
                        \Bigl(\begin{pmatrix} 0 &  0\cr Qe^{it\theta} & 0 \end{pmatrix}, \begin{pmatrix} 0 & Pe^{-it\theta} \cr 0 & 0   \end{pmatrix}\Bigr), & \text{if $x \in D_-$}
                    \end{cases}
\end{eqnarray*}
It is clear that $P,Q,\ln(1+PQ)$ are still in $L^2\cap L^\infty$. Indeed, $\ln(1+PQ) \equiv \ln(1+pq)$. The reduction order is:
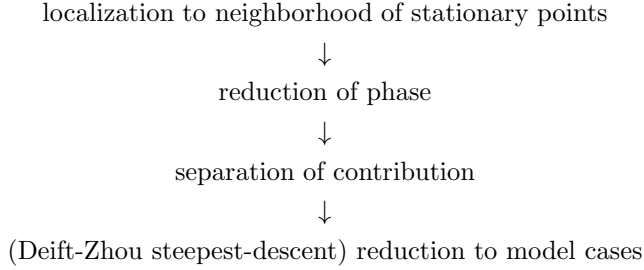
\begin{eqnarray*}
&\text{localization to neighborhood of stationary points}&\\
&\downarrow&\\
&\text{reduction of phase}&\\
&\downarrow&\\
&\text{separation of contribution}&\\
&\downarrow&\\
&\text{(Deift-Zhou steepest-descent) reduction to model cases}&
\end{eqnarray*}
The middle two reductions are iterated through the list of the stationary points. In general, the order of these two reductions can be interchanged; we choose the above order to avoid undesired regularity assumption on $\theta$. Part I (this section) describes of the first three reductions.

For simplicity, we'll assume that $p,q \in H^{k,0}$ with $k \geq 2$, but overall it will be clear from the argument that $k\geq 2$ is not always required. In particular, $k=1$ is always allowed in the separation of contributions.

We'll always assume that $p$ and $q$ and their relevant derivatives have sufficient decay at $\infty$.

\subsection{Localization to neighborhood of stationary points} \label{nbhreductsect} The argument used here will be used as a model for future reductions, and we will often refer to it whenever the needed proofs are essentially similar.

Let $1=\phi_+ + \phi_- + \phi_0$ be a $C^\infty$ partition of unity, where
\begin{itemize}
\item[(a)] $0\leq \phi_+, \phi_-,\phi_0 \leq 1$,
\item[(b)] $\phi_+, \phi_-$ are respectively supported in $D_+, D_-$, and
\item[(c)] supp($\phi_0$) consists of small disjoint neighborhoods of $\lambda_1, \dots, \lambda_N$, and $\phi_0(x)=1$ for $x$ sufficiently close to any stationary point.
\end{itemize}

Our goal in this subsection is to reduce $P,Q$ to $P\phi_0, Q\phi_0$. The new weights are denoted by $w_L^\pm$ and the new solution is denoted by $M_L$.

In the schemes described in Section~\ref{localizationsection}, we decompose $\Delta w:=w - w_L$ using the Hardy decomposition. Recall the notation $H_p(w) = \|C_+(w^-)\|_p + \|C_-(w^+)\|_p$ for a pair $w=(w^-,w^+)$. Thanks to correct phase-weight relation, the following estimates are true for $2\leq p< \infty$:
\begin{eqnarray}
\label{1stschemecomp}  H_p(1_{D_-}\Delta w), H_p(1_{D_+}\Delta w) &\lesssim& t^{-(k-1+\frac{1}{p})}
\end{eqnarray}
When $k\geq 2$, $p=\infty$ is allowed in (\ref{1stschemecomp}); these $L^\infty$ decays will be useful for the reduction of unique solvability and resolvent bound (i.e. for Theorem~\ref{uniquesolvability}). On the other hand, if $k_\theta \leq 2$ and there are no defocusing stationary point then the reduction of $u(t)$, $v(t)$ in this section can be done without requiring these $L^\infty$ estimates, consequently $k=1$ is allowed.

For simplicity of notation, let weights $w^\pm_K$ be defined by
$$w^+_K- w^+_L = C_-(\Delta w^+), \;\;\; w^-_K- w^-_L = C_+(\Delta w^-)$$
(\ref{1stschemecomp}) implies that for $2\leq p<\infty$
$$\|w_K - w_L\|_p \;\lesssim\; t^{-(k-1+\frac{1}{p})}$$
so in principle we can go from $w_K$ to $w_L$ using the first scheme. To go from $w$ to $w_K$, in the spirit of the second scheme we introduce $M^\Phi = M^\delta(I+\Phi)$, where $\Phi_\pm$ are respectively defined on the upper half and lower half planes of $\mathbb {C\setminus R}$. Heuristically,
\begin{eqnarray*}
(\Phi)_+      &\approx& -C_+(\Delta w^+) = - C_+(1_{D_+}\Delta w^+) - C_+(1_{D_-}\Delta w^+)  \equiv (\Phi_1)_+ + (\Phi_2)_+\\
(\Phi)_-      &\approx& -C_-(\Delta w^-) = - C_-(1_{D_+}\Delta w^-) - C_-(1_{D_-}\Delta w^-)  \equiv (\Phi_1)_- + (\Phi_2)_-
\end{eqnarray*}
The triangularity of $((\Phi_1)_-,(\Phi_1)_+)$ and $(\Delta w^-,\Delta w^+)$ agree on $D_+$ and disagree on $D_-$ (opposite story for $\Phi_2$). Thus, at a time only one term on the right-hand side has correct triangularity needed for the second scheme. It turns out that, as $t\to\infty$, the other term is always small: $(\Phi_1)_\pm$ are small on $D_-$ and $(\Phi_2)_\pm$ are small on $D_+$. This observation indicates that, modulo small ``noise",  the matrix structure of $\Phi_1$ and $\Phi_2$ are good enough for our applications.

The simple choice $\Phi = \Phi_1 + \Phi_2$ turns out to be inconvenient since it is not guaranteed that $I+\Phi$ is invertible on $\mathbb {C\setminus R}$. Although when $k\geq 2$ this would not be a problem since it can be shown that $\Phi$ is asymptotically strictly-triangular, it is better to avoid this restriction. A better ``superposition" of $\Phi_1, \Phi_2$ is:
\begin{eqnarray}
\label{Phieqn} \Phi &=& \Phi_1 + \Phi_2 + \Phi_2\Phi_1
\end{eqnarray}
(this has the algebraic advantage that $\det(I+\Phi)=\det(I+\Phi_2)\det(I+\phi_1)=1$ since $\Phi_i$ are strictly triangular.) This choice of $\Phi$ goes back to \cite{varzuginyd} but can be naturally interpreted as successive applications of the second scheme by $\Phi_2$ and $\Phi_1$. The jump matrix for $M^\Phi$ is:
\begin{eqnarray*}
  J^\Phi &=& (I+\Phi_-)^{-1}J(I+\Phi_+) = (I-w_\Phi^-)^{-1} (I  + w_\Phi^+)
\end{eqnarray*}
here the weights $w_\Phi$ are defined via
$$(I+w^+_\Phi) = (I+w^+)(I+\Phi_+)$$
$$(I-w^-_\Phi) = (I-w^-)(I+\Phi_-)$$
We can compute them explicitly below:
$$w_\Phi^- = I-(I-w^-)(I+\Phi_-) = w^- - \Phi_- + w^- \Phi_-$$
$$w_\Phi^+ = (I+w^+)(I+\Phi_+) - I = w^+ + \Phi_+ + w^+ \Phi_+$$
Because of the complication arose from lack of uniform triangularity, $w_\Phi$ differs $w_K$ by some noise, and the sequence of subreductions will be $w \to w_\Phi \to w_K \to w_L$. For convenience we group the last two into Proposition~\ref{1stschemeprop}. In the following two propositions, part (i) should always be viewed as the reduction of Theorem~\ref{uniquesolvability} for nice $p(x),q(x)$, while the other parts are for the reduction of $u(t),v(t)$ assuming the validity of Theorem~\ref{uniquesolvability} and Corollary~\ref{coruniqsolv}.

\begin{proposition}[$w \to w_\Phi$] \label{2ndschemeprop} (i) If $k\geq 2$ then for any $1< p<\infty$ the following estimates are equivalent:
$$\|(1-C_w)^{-1}\|_{L^p \to L^p}= O_p(1)  \text{ as $t\to\infty$}$$
$$\|(1-C_{w_\Phi})^{-1}\|_{L^p \to L^p} =O_p(1)  \text{ as $t\to\infty$}$$
(ii) Assume that $k \geq 1$ and both resolvent operators exist on $L^2$ (although their norms may not be uniformly bounded as $t\to\infty$). Then
$$\limsup_{\lambda \to \infty} |\lambda (M^\delta(\lambda) - M^\Phi(\lambda))| \lesssim t^{-k}$$
\end{proposition}

\begin{proposition}[$w_\Phi \to w_L$] \label{1stschemeprop} (i) If $k\geq 2$ then for any $1< p<\infty$ the following estimates are equivalent:
$$\|(1-C_{w_L})^{-1}\|_{L^p \to L^p}= O_p(1) \text{ as $t\to\infty$}$$
$$\|(1-C_{w_\Phi})^{-1}\|_{L^p \to L^p} =O_p(1) \text{ as $t\to\infty$}$$

(ii) Suppose that as $t\to\infty$ both $(1-C_{w_L})^{-1}$ and $(1-C_w)^{-1}$ are uniformly bounded on $L^2$. Then for $k\geq 2$ we have
$$\limsup_{\lambda \to \infty} |\lambda (M^\Phi(\lambda) - M_L(\lambda))| \lesssim_\epsilon t^{-(k-1+\frac{1}{k_\theta+1})+\epsilon}$$

(iii) Suppose that as $t\to\infty$ both $(1-C_{w_L})^{-1}$ and $(1-C_w)^{-1}$ are uniformly bounded on $L^p$ for every $2 \leq p <\infty$. Then for $k \geq 1$ we have
$$\limsup_{\lambda \to \infty} |\lambda (M^\Phi(\lambda) - M_L(\lambda))| \lesssim_\epsilon t^{-(k-\frac{1}{2}+\frac{1}{2(k_\theta+1)})+\epsilon}$$
\end{proposition}

Recall that by the A-B decomposition lemma, $\Omega_0^{\pm 1}$ have $A_k+B_k$ decompositions.

\proof[Proof of Proposition~\ref{2ndschemeprop}]

(i) Let $\widetilde{\Phi} = -\Phi_1 - \Phi_2 + \Phi_1 \Phi_2$, then
$$(I+ \widetilde{\Phi})(I+\Phi) = (I+\Phi)(I+ \widetilde{\Phi}) = I$$
The desired claim is now a consequence of (\ref{factorize1}),(\ref{factorize2}), and (\ref{factorize3}). Here, note that the  $L^\infty(\mathbb R)$ norms of $\Phi_\pm$ and $\widetilde{\Phi}_\pm$ can be controlled by $\|(\Phi_1)_\pm\|_\infty + \|(\Phi_2)_\pm\|_\infty + \|(\Phi_1)_\pm\|_\infty \|(\Phi_2)_\pm\|_\infty$ which are bounded as $t\to\infty$ (for $k\geq 2$).

(ii) Let $k\geq 1$. Oscillation of $\Delta w$ ensures that
\begin{eqnarray*}
\limsup_{\lambda\to\infty} |\lambda C(1_{D_\pm}\Delta w^\pm)(\lambda)| = \frac{1}{2\pi} |\int 1_{D_\pm}\Delta w^\pm(x) dx| \lesssim t^{-k}
\end{eqnarray*}
The above nontangential limits converge thanks to integrability at $\infty$ of $p,q$. This is the only place where this requirement is used and future implicit constants will not depend on these two $L^1$ norms. The above estimates imply
$$\lim_{\lambda \to \infty} |\lambda \Phi_j(\lambda)|     \lesssim \ t^{-k}$$
Consequently,
$$\limsup_{\lambda \to \infty} |\lambda (M^\delta(\lambda) - M^\Phi(\lambda))| \;\lesssim\; \limsup_{\lambda \to \infty} |\lambda \Phi(\lambda)| \lesssim t^{-k}$$
here recall that $\lim_{\lambda \to \infty} M^\delta(\lambda) = I$.
\endproof

\proof[Proof of Proposition~\ref{1stschemeprop}] (i) We'll show that, for $k\geq 1$ and $2\leq q < \infty$,
\begin{eqnarray}
\label{wPhiwL} \|w_\Phi - w_L\|_q \lesssim t^{-(k-1+\frac{1}{q})+\epsilon},
\end{eqnarray}
and if $k\geq 2$ then $q=\infty$ is allowed. This will automatically imply (i). We'll show the estimate for $\|1_{D_+}(w^+_\Phi - w^+_L)\|_p$ and the argument for other sign combinations are entirely similar. Write
\begin{eqnarray}
\nonumber 1_{D_+}(w^+_\Phi - w^+_L)
                    &=& 1_{D_+}(w^+_\Phi - w^+) + 1_{D_+}(w^+ - w^+_L)\\
\label{1DplusPhiL}  &=& 1_{D_+}(\Phi_+ + w^+\Phi_+) + 1_{D_+}\Delta w^+
\end{eqnarray}
Notice that as $t\to\infty$, $w^+$ is uniformly bounded, $(\Phi_1)_+$ is uniformly bounded in $L^{\infty-}$, and $1_{D_+}(\Phi_2)_+ \equiv -1_{D_+}C_-(1_{D_-}\Delta w^+)$ decays like $O(t^{-(k-1+\frac{1}{p})})$ in $L^p$. Consequently by H\"older's inequality, $(\Phi_1)_+(\Phi_2)_+$ decays like $O(t^{-(k-1+\frac{1}{p})+\epsilon})$ in $L^p$. Substitute $\Phi_+$ using (\ref{Phieqn}) we see that the significant terms in $1_{D_+}(\Phi_+ + w^+\Phi_+)$ are:
$$1_{D_+}(\Phi_1)_+ , 1_{D_+}w^+(\Phi_1)_+$$

Now, since both $w^+, (\Phi_1)_+$ has the same strictly-upper-triangular matrix structure on $D_+$ (this is the way the second perturbation scheme was designed), the term $ 1_{D_+}w^+(\Phi_1)_+$ is effectively $0$. Consequently, the main contribution in (\ref{1DplusPhiL}) comes from
\begin{eqnarray*}
1_{D_+}(\Phi_1)_+ + 1_{D_+}\Delta w^+
&=& 1_{D_+}\Big((\Phi_1)_+ + 1_{D_+}\Delta w^+\Big)\\
&=& 1_{D_+}\Big(-C_+(1_{D_+}\Delta w^+) + C_+(1_{D_+}\Delta w^+) - C_-(1_{D_+}\Delta w^+)\Big)\\
&=& -1_{D_+}C_-(1_{D_+}\Delta w^+)
\end{eqnarray*}
which is small $L^p$. This gives us the desired estimate. When $k\geq 2$ the argument can be repeated and the endpoint case $q=\infty$ is allowed.

(ii) In this part we assume $k\geq 2$. By part (i), we also have uniform boundedness of $(1-C_{w_\Phi})^{-1}$ in $L^2$ as $t\to\infty$. Using Lemma~\ref{direct_lem} for $w_1 = w_L$ and $w_2 =w_\Phi$, we can control $\limsup_{\lambda \to \infty}|\lambda(M^\Phi (\lambda)- M_L(\lambda))|$ by:
\begin{eqnarray}
\nonumber       && \|w_\Phi -w_L\|_2\Big(H_2(w_\Phi)+H_2(w_L)\Big) \;\;+\;\; \|w_\Phi -w_L\|_\infty H_2(w_\Phi)H_2(w_L)\\
\label{MPhiML}  &+& |\int (w^+_\Phi -w^+_L) w^+_L |  + | \int (w^-_\Phi -w^-_L) w^-_L|\\
\nonumber       &+& \limsup_{\lambda \to \infty} \big| \lambda C (w_\Phi -w_L)(\lambda) \big|
\end{eqnarray}

\noindent \underline{First two terms}: By Corollary~\ref{Hardyrough}
$$H_2(w_L) \lesssim t^{-\frac{1}{2(k_\theta+1)}}$$
thus using (\ref{wPhiwL}) this estimate is also true for $H_2(w_\Phi)$. Consequently, using Lemma~\ref{localprpl}, the first two terms in (\ref{MPhiML}) are controlled by
$$t^{-(k-\frac{1}{2})}t^{-\frac{1}{2(k_\theta+1)}} + t^{-(k-1)}t^{-\frac{1}{2(k_\theta+1)}}t^{-\frac{1}{2(k_\theta+1)}} \lesssim t^{-(k-1+\frac{1}{k_\theta+1})}$$
recall that $k_\theta:=\max\{0,k_1,\dots, k_N\}$).

\noindent \underline{Two middle terms}: We'll estimate $\int 1_{D_+} (w^+_\Phi - w^+_L)w^+_L$, the other cases can be done similarly. By triangularity,
$$1_{D_+}(w^+_\Phi - w^+_L)w^+_L = 1_{D_+}\Big((\Phi_2)_+ w^+_L \;\;+\;\; (\Phi_1)_+(\Phi_2)_+ w^+_L \;\;+\;\; w^+ (\Phi_2)_+ w^+_L\Big)$$
The contribution of $1_{D_+}(\Phi_2)_+ w^+_L$ can be decomposed into
\begin{eqnarray*}
\int 1_{D_+}(\Phi_2)_+ C_+(w^+_L) - \int 1_{D_+} (\Phi_2)_+ C_-(w^+_L)
\end{eqnarray*}
Since $(\Phi_2)_+$ is small on $D_+$, using Cauchy-Schwarz the second integral is controlled by $t^{-(k-\frac{1}{2})} t^{-\frac{1}{2(k_\theta+1)}}$. To estimate the first integral, we'll use analytic continuation to $\mathbb{C}_+$ and Cauchy theorem. Eventually it suffices to show that along any vertical ray $\Gamma$ in $\mathbb{C}_+$ originated from arbitrary stationary point, $\int_{\Gamma}(\Phi_2)_+ C_+(w^+_L)$ is small. Using the complex variants of Lemma~\ref{localprpl} and Lemma~\ref{mainlemma} (i.e. Lemma~\ref{variantLP} and Lemma~\ref{variantMainLemma}), this follows easily: in $L^2(\Gamma)$, $(\Phi_2)_+$ have strong decay and $C_+(w_L^+)$ has some nontrivial decay as $t\to\infty$. Thus, by Cauchy-Schwarz we get an estimate of size
$$t^{-(k+\frac{1}{2(k_\theta+1)})}$$

The contribution of other terms can be estimated similarly. It might be convenient to notice that $1_{D_+}w^+ (\Phi_2)_+ w^+_L$ has at most one nonzero entry and the above argument should be applied to that entry.

\noindent \underline{Last term}: We'll show that
$$\limsup_{\lambda \to \infty} |\lambda C(w^+_\Phi - w^+_L)(\lambda)| \lesssim t^{-(k-\frac{1}{2}+\frac{1}{2(k_\theta+1)})}$$
the estimate for $w^-_\Phi - w^-_L$ is similar. Now,
\begin{eqnarray*}
  w^+_\Phi - w^+_L
  &=& (\Phi_1)_+ + (\Phi_2)_+ + \text{other terms}\\
  &=& C_+(\Delta w^+) + \text{other terms}
\end{eqnarray*}
where the other terms are in $L^1$. Notice that $C_+C_- = 0$ and $C_-C_-=C_-$, so
\begin{eqnarray*}
\limsup_{\lambda \to \infty} \big| \lambda C(C_+(\Delta w^+))(\lambda) \big|
   &=&  \limsup_{\lambda \to \infty} \big|\lambda C_+(\Delta w^+)(\lambda)\big|\\
   &=& O(t^{-k}), \;\;\text{by linear theory }
\end{eqnarray*}
The other terms are $(\Phi_1)_+(\Phi_2)_+$, $w^+ (\Phi_1)_+(\Phi_2)_+$, $w^+ (\Phi_1)_+$, $w^+ (\Phi_2)_+$, which are in $L^1$. Thus,
$$\limsup_{\lambda\to\infty} |\lambda C(\text{other terms})(\lambda) |= \Big| \frac{1}{2\pi}\int_{\mathbb R} \text{other terms}\Big|$$

By orthogonality, the contribution of $(\Phi_1)_+(\Phi_2)_+$ in $\int_{\mathbb R}$ is $0$, and the contribution of $w^+ (\Phi_1)_+(\Phi_2)_+$ becomes
$$\int_{\mathbb R} C_-(w^+) (\Phi_1)_+(\Phi_2)_+ \;\;\lesssim\;\; H_2(w) \|(\Phi_1)_+(\Phi_2)_+\|_2$$
using Cauchy-Schwarz. Since $(\Phi_1)_+$ is small on $D_-$ and  $(\Phi_2)_+$ is small on $D_+$, the above inequality gives an estimate of size $t^{-\frac{1}{2(k_\theta+1)}}t^{-(k-\frac{1}{2})}$.

Notice that by triangularity, $\int 1_{D_+}w^+ (\Phi_1)_+ = \int 1_{D_-}w^+ (\Phi_2)_+ = 0$. Now, the argument used to estimate the two middle terms of (\ref{MPhiML}) can be applied to estimate $\int 1_{D_-}w^+(\Phi_1)_+$ and $\int 1_{D_+} w^+(\Phi_2)_+$. Consequently, we have an estimate of size $t^{-\frac{1}{2(k_\theta+1)}}t^{-(k-\frac{1}{2})}$ for the contribution of $w^+ (\Phi_1)_+$ and $w^+ (\Phi_2)_+$.

(iii) Suppose that $k\geq 1$ and $\|(1-C_{w_L})^{-1}\|_{L^p \to L^p}, \|(1-C_{w})^{-1}\|_{L^p \to L^p} \lesssim_p 1$ (as $t\to\infty$) for any $2\leq p <\infty$. For $p,q \in H^{1,0}$ with sufficient decay, we'll have $\Phi_1,\Phi_2\in L^\infty$. This allows us to use (\ref{factorize1}),(\ref{factorize2}) and (\ref{factorize3}), but these decay requirements of $p$ and $q$ will not contribute to the finiteness of the implicit constants.

Under the assumption $k\geq 1$, $\|\Phi_1\|_\infty$,$\|\Phi_2\|_\infty$ are however not known to be uniformly bounded as $t\to\infty$. Consequently, we won't have uniform boundedness of $(1-C_{w_\Phi})^{-1}$ in $L^p$.  On the other hand, we know that $\Phi_1,\Phi_2$ are uniformly bounded in $L^{\infty-}$ as $t\to\infty$. Using (\ref{factorize1}),(\ref{factorize2}) and (\ref{factorize3}), it is then not hard to show that $(1-C_{w_\Phi})^{-1}$ is uniformly bounded from $L^{p+} \cap L^p$ to $L^p$ for any $2\leq p<\infty$. This allows us to use a variant of Lemma~\ref{direct_lem} (see the remarks after this lemma for details).

Now, the same argument as in (ii) can be repeated, where to improve the estimate of (ii) to
$$O_\epsilon(t^{-(k-\frac{1}{2} + \frac{1}{2(k_\theta+1)})+\epsilon})$$
it suffices to show that for $2 \leq q < \infty$
\begin{eqnarray*}
 \|w_\Phi -w_L\|_{\frac{2q}{q-2}} \Big(H_q(w_L) + H_{q+}(w_L)\Big) \lesssim t^{-(k-\frac{1}{2}-\frac{1}{q})}
\end{eqnarray*}
(after that choosing $q$ sufficiently large will give us the desired estimate). We note that in other places an $\epsilon$ in the decay order might also be lost, one reason is the lack of an uniform bound as $t\to\infty$ of $\|\Phi_1\|_\infty$ and $\|\Phi_2\|_\infty$ (we need to use $L^{\infty-}$ instead of $L^\infty$). Now, the above inequality is a consequence of (\ref{wPhiwL}).
\endproof

\noindent \emph{Remarks:} Let $A$, $B$, $C$ respectively be the decay estimates for $H_2(w_L)$, $\|1_{D_\pm}(w_\phi - w_L)\|_2$ and $\|w_\phi - w_L\|_\infty$. Then typically (modulo an $\epsilon$) our argument gives an overall estimate of $AB$ in part (iii), and an overall estimate of $\max(AB, A^2 C)$ in part (ii). We'll always have $A = t^{-\frac{1}{2(k_\theta+1)}}$, thanks to Corollary~\ref{Hardyrough} and the correct phase-weight relation in our weights.

\subsection{Reduction of phase} \label{reductphasesect} In this section, we study a RHP localized to a small neighborhood of the stationary points of $\theta$. Let $\lambda_0$ be a stationary point of $\theta$. Below, we'll show the phase reduction of $\theta$ near $\lambda_0$.

By translation symmetry, assume that $\lambda_0 = 0$ and is of order $k_0$. Let the weights of our RHP be denoted by $w^\pm$. On the neighborhood of $0$ where the weights are supported, approximate $\theta$ by
\begin{eqnarray*}
    \theta(0) +  \frac{\theta^{(k_0+1)}(0)}{(k_0+1)!}x^{k_0+1}
\end{eqnarray*}
and keep the value of $\theta$ near other stationary points. Let $\theta_0$ denote this approximation. Without loss of generality we can assume that the portion of $supp(w)$ near $0$ is small enough so on a neighborhood of this part $\theta$ is $(k_0+1)$-time  differentiable, with
$$\theta'(0) = \dots = \theta^{(k_0)}(0) = 0,\; \theta^{(k_0+1)}(0) \neq 0$$
$$\theta^{(k_0+1)}(x) = \theta^{(k_0+1)}(0) + O(|x|^{\beta}), \;\beta >0$$

Below we show how to reduce $\theta \to \theta_0$. The following assumption is required only if $\theta\not\equiv \theta_0$ in the above neighborhood (since otherwise we don't need this reduction).

(*) $p,q$ have two $L^2$ derivatives on a neighborhood of $0$ if $k_0=1$. Three $L^2$ derivatives are required if $k_0\geq 3$, or $k_0=2$ and $p(0)q(0)<0$.

The regularity assumption on $p,q$ for $k_0\geq 3$ may be improved given better understanding of model RHPs associated with a stationary point of such orders. More precisely, we don't know if the corresponding resolvent operator $(1-C_w)$ associated with such model RHP is invertible on $L^p$ for large $p$, therefore in our direct scheme we need some nontrivial decay for $\|\Delta w\|_\infty$, which is the source of this extra requirement. For this reason, if $|pq|$ are very small at a stationary point of such order then this extra requirement can be removed.

Notice that $0$ is also a stationary point of order $k_0$ for $\theta_0$: $\theta_0^{(k_0+1)}(0) = \theta^{(k_0+1)}(0) \neq 0$. Furthermore, by L'Hospital's rule and continuity of $\theta^{(k_0+1)}$ at $0$ we have:
\begin{eqnarray*}
    \lim_{x\to0} \frac{\theta'_0(x)}{\theta'(x)} = \dots = \lim_{x\to 0} \frac{\theta^{(k_0+1)}_0(x)}{\theta^{(k_0+1)}(x)} = 1
\end{eqnarray*}

Thus, our neighborhood of $0$ can be chosen such that $\frac{\theta'_0(x)}{\theta'(x)} \gtrsim 1$ for every $x \neq 0$ in this neighborhood. In particular, this means $\theta(x,s):=s\theta_0(x) + (1-s)\theta(x)$ has a stationary phase point at $0$ of order $k_0$ for every $s\in [0,1]$; furthermore, for $x$ near $0$, $\frac{d}{dx}\theta(x,s)$ has the same signs as $\theta'(x)$.

Let $w^\pm_0$ be the weight obtained from $w^\pm$ after replacing $\theta$ by $\theta_0$. For simplicity of notation, let $\Delta w = w_0  - w$, which is supported near $0$. To reduce from $w$ to $w_0$, we will proceed exactly as in the reduction from $w^\Phi$ to $w_L$. For convenient, in the following propositions we'll denote by $k\geq 2$ the number of $L^2$ derivatives that $p$ and $q$ are required to have in the assumption (*). The argument of Section~\ref{nbhreductsect} works once we show the following estimates and its complex variant:

\begin{proposition} \label{phasereduct}  For any sign combinations,
\begin{eqnarray}
\label{phasereductpw} \limsup_{\lambda \to \infty} |\lambda C \Big(1_{D_\pm}\Delta w^\pm \Big)| &\lesssim& t^{-(\beta+1)\frac{1}{k_0+1}+\epsilon}
\end{eqnarray}
Furthermore, if $k_0>1$ then for $2\leq p\leq \infty$
\begin{eqnarray}
\label{phasereductLp} H_p(1_{D_+}\Delta w), \; H_p(1_{D_-}\Delta w) &\lesssim& \max \Big(t^{-(k-2+\frac{1}{p})}, t^{-(\beta+\frac{1}{p})\frac{1}{k_0+1}+\epsilon}\Big)
\end{eqnarray}
If $k_0=1$ and $\theta^{(3)}$ is assumed $L^r$ integrable near $0$ then the above estimates remains true after the following adjustments: $k_0+1$ is replaced by $k_0+1+c(2,p)$, and $\beta$ is replaced by $\beta-c(2,p)$ for $2\leq p<\infty$, or $\beta-2c(2,p)$ for $p=\infty$.
\end{proposition}
\noindent Note: $c(k,p)$ is defined in (\ref{ckp}). In particular, $c(2,2) = \frac{2}{3r}$ and $c(2,\infty) = \frac{1}{2r}$. Recall the assumption (B) on $\theta$, which in particular says that
$$\beta>\begin{cases}\frac{1}{r}, & \text{if $k_0=1$}\\ 0,& \text{if $k_0>1$.}\end{cases}$$
The reason for this assumption will be explained shortly.

\begin{proposition} \label{phasereductvariant} Let $1<p<\infty$. If $\Gamma$ is a ray originating from a stationary point that forms a nontrivial angle with $\mathbb R$ and $k_0>1$, then
$$\|C(1_{D_\pm} \Delta w)\|_{L^p(\Gamma)} \lesssim \max(t^{-(k-1)}, t^{-(\beta+\frac{1}{p})\frac{1}{k_0+1}+\epsilon})$$
If $k_0=1$ then this estimate has to be adjusted as in Proposition~\ref{phasereduct}.
\end{proposition}
Before proving these propositions, we make some comments. Recall that the phase reduction and separation of contribution will be iterated through the list of stationary points. Using the above estimates, the effect on $u(t), v(t)$ can always be controlled by $AB + A^2 C$, where $A=t^{-\frac{1}{2(k_\theta+1)}}$, $B$ is the estimate for $H_2(\Delta w)$, and $C$ is the estimate for $H_\infty(\Delta w)$, and all we need is $L^2$ boundedness of $(1-C_w)^{-1}$, which is always the case (posteriori). However, if $(1-C_w)^{-1}$ is also bounded in $L^p$ for large $p$, we can improve this estimate to $AB$ (modulo an $\epsilon$ in the decay exponent). For this reason, if there are more than one stationary point we'll always start out with stationary points of order $\geq 3$ and defocusing secondary stationary points \textit{before} the rest.
The main reason is whenever these bad stationary points are involved, we will not have $L^p$ boundedness, but after separating them out $L^p$ boundedness suddenly becomes available, which enable us to avoid strong assumptions on regularity of $\theta,p,q$ near the good stationary points.

Below we show that our estimate of the effect on $u(t),v(t)$ can be controlled by something decay better than $t^{-\frac{1}{k_\theta+1}}$. If we only use $L^2$ boundedness of $(1-C_w)^{-1}$ then we'll have a decay of $t^{-(\frac{1}{k_\theta+1}+d_0)+\epsilon}$.
$$d_0=\frac{\beta}{k_0+1}>0$$
assuming three $L^2$ derivatives of $p,q$ near $0$ and boundedness of $\theta^{(4)}$. These conditions are satisfied when $k_0\geq 3$ or $k_0=2$ with $p(0)q(0)<0$ by our assumptions; we note that the boundedness of $\theta^{(4)}$ can be changed to some high $L^r$ integrability condition similar to the assumption on primary stationary points; but we avoid doing that here for the sake of simplicity.

On the other hand, if we assume high $L^p$ boundedness of $(1-C_w)^{-1}$ then the effect on $u(t),v(t)$ can be controlled by $t^{-(\frac{1}{2(k_\theta+1)}+d_0)+\epsilon}$, where $d_0>\frac{1}{2(k_\theta+1)}$. More specifically, assuming only two $L^2$ derivatives of $p,q$ near $0$ we can get
$$d_0 = \min(\frac{1}{2}, (\beta+\frac{1}{2})\frac{1}{k_0+1}) > \frac{1}{2(k_0+1)}$$
if $k_0>1$. When $k_0=1$, if $\theta^{(3)}$ is assumed to be $L^r$ integrable near $0$ then
$$d_0 = \min\Big(\frac{1}{2}, (\beta+\frac{1}{2}-\frac{2}{3r})\frac{1}{2+\frac{2}{3r}}\Big)$$
In this case, using the assumption $\beta>\frac{1}{r}$ we have
$$(\beta+\frac{1}{2}-\frac{2}{3r})\frac{1}{2+\frac{2}{3r}})>(\frac{1}{2}+\frac{1}{3r})\frac{1}{2+\frac{2}{3r}}) \geq \frac{1}{4}$$
Consequently, only two $L^2$ derivatives is needed for $p,q$ near $0$ if $0$ is a good stationary point (i.e. primary or focusing secondary). On the other hand, in that case the reduction regarding unique solvability and resolvent bound needs to be looked at more carefully since we won't have decay for $H_\infty(\Delta w)$ as $t\to\infty$. In this case, the above assumption on $\beta$ saves us. Indeed, by (\ref{phasereductLp}) we have
$$H_\infty(1_{D_\pm} \Delta w) = O(1) + O(t^{-c})$$
as $t\to\infty$ for some $c>0$. More specifically,
$$c=\begin{cases}\frac{\beta}{k_0+1}-\epsilon, & \text{if $k_0>1$;}\\
(\beta - \frac{1}{r})\frac{1}{2+\frac{1}{2r}}-\epsilon, & \text{if $k_0=1$.}
\end{cases}$$
so $c>0$. Now, the implicit constant of $O(1)$ is (modulo a positive power) proportional to the size of the current neighborhood. Thus, it can be made arbitrarily small. For large $t$, \emph{a posteriori} the resolvent norm $\|(1-C_{w_0})^{-1}\|$ on $L^p$ can be controlled by an absolute $p-$constant depending essentially on the values of $|p(x)|$ and $|q(x)|$ at the stationary points. For $p$ in any given compact subset of $[2,\infty)$, if $t$ is large then unique solvability follows by choosing our neighborhood small enough and using Newman series. Furthermore, on those $L^p$ we can control $\|(1-C_w)^{-1}\|$ by a comparable $p$-constant, which is good enough for our applications.

Below we'll prove Proposition~\ref{phasereduct}, the proof for Proposition~\ref{phasereductvariant} is similar. We'll use a lemma whose proof is a simple application of Fubini's theorem:
\begin{lemma} \label{phasereductlemma} Let $f(x,s)$ be in $L^1_s H^{1,1}_x((0,1) \times \mathbb R)$, i.e. $\int_0^1 \|f(x,s)\|_{H^{1,1}_x} ds < \infty$. Then $C_\pm$ and $\int_0^1 ds$ commute:
\begin{eqnarray}
\label{lemma1}   \Big(C_\pm \int_0^1 f(.,s) ds\Big)(y) =  \int_0^1 \Big(C_\pm f(.,s)\Big)(y) ds
\end{eqnarray}
\end{lemma}

\proof[Proof of lemma]
Since $H^{1,1}_x \subset L^1_x$, we know $f\in L^1_s L^1_x$. Fubini's theorem gives us:
\begin{eqnarray*}
   \widehat{LHS}(\xi) =  1_{\{\pm\xi\geq 0\}} \;\frac{1}{2\pi}\int_{\mathbb R}\int_0^1 f(x,s) e^{-ix\xi}\; ds  dx = \int_0^1 1_{\{\pm\xi\geq 0\}} \widehat{f}(\xi,s) ds
\end{eqnarray*}

Similarly (use the above argument for $g(\xi,s) := \widehat{f}(\xi,s)1_{\{\pm\xi\geq 0\}} \in L^1_s L^1_\xi$ instead of $f(x,s)$ and the inverse Fourier transform instead of the Fourier transform),
\begin{eqnarray*}
   F^{-1} \Big(\int_0^1 g(.,s) ds\Big)(y) = \int_0^1 F^{-1}(g(.,s))(y) ds = \int_0^1 \big(C_\pm f(.,s)\big)(y) ds
\end{eqnarray*}
\endproof
\proof[Proof of Proposition~\ref{phasereduct}] To show (\ref{phasereductLp}) and (\ref{phasereductpw}), we want to estimate the respective Hardy projection of $\Delta w$. Below we'll estimate $$|C_+(1_{D_+}\Delta w^-)| = |C_+(1_{D_+}\phi_0 P \big(e^{-it\theta} - e^{-it\theta_0}\big))|,$$
the other cases are similar. For simplicity we'll assume that $k_0>1$, the case $k_0=1$ is entirely similar and the adjustments in this case followed from the corresponding adjustments in Lemma~\ref{mainlemma}. Below we'll show (\ref{phasereductLp}), the argument is entirely similar for (\ref{phasereductpw}).

Write $1_{D_+}\phi_0 P \big(e^{-it\theta} - e^{-it\theta_0}\big)$ as:
$$it\int_0^1 f(x)e^{-it\theta(x,s)} 1_{D_+}(x) \Omega_0(x) ds$$
where recall that $\theta(x,s) = s \theta_0(x) + (1-s)\theta(x)$ and
$$f(x) :=   \phi_0(x) (\theta_0(x) - \theta(x)) p(x) $$
Note that $f(x)$ is compactly supported, continuously differentiable, furthermore $f$ vanishes at the endpoint of $D_+$. Thus it is not hard to see that $f1_{D_+}\Omega_0 \in L^1_s H^{1,1}_x$. Apply the above lemma, we have
$$C_+\big(\phi_0 P (e^{-it\theta} - e^{-it\theta_0})1_{D_+}\big)(\lambda) = it\int_0^1 C_+ \big(1_{D_+}f(.,s)\Omega_0e^{-it\theta(.,s)}\big)(\lambda) ds$$
Under assumption (*) and by the A-B decomposition lemma, $\Omega_0$ has an $A_k+B_k$ decomposition in a neighborhood of the stationary point $0$. Now, $f$ has multiplicity $k_0+1+\beta$ at $0$ up to the $k^{th}$ derivative. The last claim is a consequence of the H\"older condition of $\theta^{(k_0+1)}$ at $0$. By Lemma~\ref{mainlemma} we have
\begin{eqnarray*}
\|C_+ \big(1_{D_+}f(.,s)\Omega_0\big)(y)\|_p
&\lesssim_s& \max (t^{-(k-1+\frac{1}{p})}, t^{-(k_0+1+\beta+\frac{1}{p})\frac{1}{k_0+1}}) \\
&=& t^{-1} \max (t^{-(k-2+\frac{1}{p})}, t^{-(\beta+\frac{1}{p})\frac{1}{k_0+1}})
\end{eqnarray*}
Below we'll show that the dependence on $s$ in the above inequalities are mild, so that after integrating over $s\in [0,1]$, the implicit constants remain finite.

Since $f$ is compactly supported, the $s$-dependence of the implicit constants is originated from the phase dependence weights in the weighted sum $f_{\theta,k}$:
$$|\frac{d^2}{dx^2}\theta(x,s)|^{\alpha_1} \cdots|\frac{d^{k+1}}{dx^2}\theta(x,s)|^{\alpha_k}  \frac{1}{|\frac{d}{dx}\theta(x,s)|^{(k+\alpha_1+\dots+\alpha_k)-1/p}}$$
Since $s\in [0,1]$, we are not concerned about those $s$'s appearing in the numerators.

For those $s$'s that appear in the denominator, we note that: our choice of $\phi_0$ ensures that on its support, $|\frac{d}{dx}\theta(x,s)| \sim |\theta'(x)|$ uniformly over $s\in [0,1]$. Therefore we can remove the $s$-dependence in our estimates and complete the proof of this proposition. \endproof

\subsection{Separation of contributions} \label{separationsubsect} In this section, let $\lambda_0$ is a stationary point of order $k$ of $\theta$ and $w$ is a pair of weights supported near the stationary points. Let $w_0$ denote the part of $w$ supported near $\lambda_0$ and $w_1 = w - w_0$ be the remaining. We'll assume that both $w_0,w_1$ have correct-phase weight relation and the RHPs associated with $w_0$ and with $w_1$ are uniquely solvable. Furthermore, we assume quantitative bounds $\|(I-C_{w_j})^{-1}\|_{L^p \to L^p} \lesssim 1$ as $t\to\infty$ for some $2\leq p < \infty$.

Regarding unique solvability, we'll prove that

\begin{proposition} The normalized RHP with weight $w = w_0 + w_1$ is uniquely solvable for large $t$, indeed
\begin{eqnarray}
\label{contribnormest} \|(1-C_w)^{-1}\|_{L^p \to L^p} \lesssim 1+ \sum_{j=0}^1 \|(1-C_{w_j})^{-1}\|_{L^p \to L^p}
\end{eqnarray}
\end{proposition}

\proof To show that $I-C_w$ is invertible for large $t$, we'll use the parametrix $I + \sum_{j=0}^1 C_{w_j}(I-C_{w_j})^{-1} \equiv I + \sum_{j=0}^1 (I-C_{w_j})^{-1}C_{w_j}$, introduced in \cite{varzuginyd}:
\begin{eqnarray*}
(I-C_w)\Big(I + \sum_{j=0}^1 C_{w_j}(I-C_{w_j})^{-1}\Big) &=& I - \sum_{0\leq
i\neq j \leq 1}
C_{w_i}C_{w_j} (I-C_{w_j})^{-1}\\
\Big(I + \sum_{j=0}^1 (I-C_{w_j})^{-1}C_{w_j}\Big)(I-C_w) &=& I - \sum_{0\leq
i\neq j \leq 1}  (I-C_{w_j})^{-1}C_{w_j}C_{w_i}
\end{eqnarray*}

By Corollary~\ref{almostorthogonal}, for $i\neq j$, $\|C_{w_i}C_{w_j}\|_{L^p \to L^p}$ decays as $t\to\infty$:
$$\|C_{w_i}C_{w_j}\|_{L^p \to L^p}\lesssim t^{-\frac{1}{p(k_\theta+1)}}$$

Thus, for large $t$ we know that $\sum_{0\leq i\neq j \leq 1} C_{w_i}C_{w_j} (I-C_{w_j})^{-1}$ and $\sum_{0\leq
j\neq i \leq 1}  (I-C_{w_j})^{-1}C_{w_j}C_{w_i}$ have small norm. Consequently, $I-C_w$ is invertible, and
\begin{eqnarray*}
  \|(I-C_w)^{-1}\|_{L^p \to L^p} &\lesssim& \|I + \sum_{j=0}^1 C_{w_j}(I-C_{w_j})^{-1}\|_{L^p \to L^p} \\
                   &\lesssim& 1 + \sum_{j=0}^1 \|(I-C_{w_j})^{-1}\|_{L^p \to L^p}
\end{eqnarray*}
\endproof
For $0\leq j \leq 1$, denote by $M_j,\mu_j$ the respective ingredients of the normalized $L^2$ RHP with weights $(w^-_j, w^+_j)$. To separate the contributions of different stationary points, we'll need the following a priori estimates:

\begin{lemma}
\label{apriori}
Assuming that near $\lambda_0$ the phase $\theta$ is of the form $a+b(x-\lambda_0)^{k+1}$. Assume that locally $p,q$ have two $L^2$ derivatives. Then if $P \subset \mathbb R$ so that $\emph{distance}(P,\text{supp}(w_0)) \gtrsim 1$ and $\emph{distance}(P,\lambda_0)\gtrsim 1$ then for a.e. $\lambda\in P$ we have:
\begin{eqnarray}
\label{pwmujest} \mu_0(\lambda,t) = I + \begin{pmatrix} 0 & \frac{u_0(t)}{\lambda-\lambda_0}\cr \frac{v_j(t)}{\lambda-\lambda_0} & 0\end{pmatrix} + O_\epsilon(t^{-\frac{3}{2(k+1)}+\epsilon})
\end{eqnarray}
Here $u_0$ and $v_0$ are the contribution of $\lambda_0$ to the leading asymptotics of $u(t),v(t)$, defined as in Theorem~\ref{maintheorem}.
\end{lemma}
We note that the above assumption on $\theta$ is automatic because at this moment $\theta$ has already been reduced to a nice analytic phase near $\lambda_0$. Furthermore, if $k=1$ then the Lemma remains true for $p,q$ with one $L^2$ derivative. On the other hand, for purely harmonic analysis interests, the above asymptotics remains true for more general class of $\theta$, for instance those having three Lipschitz derivatives near $\lambda_0$ and satisfying
$$|\theta^{(k+1)}(x)-\theta^{(k+1)}(\lambda_0)| \lesssim |x-\lambda_0|^\beta, \;\;\beta > \frac{1}{2}$$
(the error term has weaker decay if $\beta<1$), however the main issue here is that the Lemma requires higher regularity assumptions on $p,q$ to allow for this generality. Indeed, the phase reduction of this Lemma is fairly expensive, it may require $p,q$ to have three $L^2$ derivatives near $\lambda_0$. This is the reason why for every stationary point we always carry out phase reduction before separating its contribution.

It would be interesting to see if the decay estimate $O_\epsilon(t^{-\frac{3}{2(k+1)}+\epsilon})$ can be improved to $O_\epsilon(t^{-\frac{2}{k+1}+\epsilon})$ assuming enough regularity of $p,q$. That improvement would allow us to obtain a stronger decay estimate for the error term in the leading asymptotics of $u(t)$ and $v(t)$, when there are more than one stationary point.

The proof of Lemma~\ref{apriori} is a sequence of reductions; at the end it will be proved for the corresponding model RHP associated with $\lambda_j$.

For technical reason, from now on we'll assume that
\begin{eqnarray}
\label{uvest}
\begin{cases}\limsup_{\lambda \to \infty} \Big(\big|\lambda(M_0)_{12}(\lambda,t)\big| +\big|\lambda(M_0)_{21}(\lambda,t)\big|\Big) = O_t(1)\\
 \limsup_{\lambda \to \infty} \Big(\big|\lambda(M_1)_{12}(\lambda,t)\big| +\big|\lambda(M_1)_{21}(\lambda,t)\big|\Big) = O_t(1)
 \end{cases}
\end{eqnarray}
Since the reduction of $u(t),v(t)$ implies the reduction of these estimates, this assumption is harmless.

Assuming these a priori estimates, we will show below that the leading asymptotics of $u(t)$ and $v(t)$ are essentially the sum of the leading asymptotics of $u_0(t) + u_1(t)$ and $v_0(t) + v_1(t)$. As usual, we'll suppress $t$ for brevity.

\begin{proposition} Assuming the a priori estimates (\ref{pwmujest}),(\ref{uvest}), we have
\begin{eqnarray}
\label{contribest1} \limsup_{\lambda \to \infty} |\lambda\Big(M_{12}(\lambda)- (M_0)_{12}(\lambda) - (M_1)_{12}(\lambda)\Big)| \lesssim_\epsilon t^{-\frac{3}{2(k_\theta+1)}+\epsilon}\\
\label{contribest2} \limsup_{\lambda \to \infty} |\lambda\Big(M_{21}(\lambda)- (M_0)_{21}(\lambda) - (M_1)_{21}(\lambda)\Big)| \lesssim_\epsilon t^{-\frac{3}{2(k_\theta+1)}+\epsilon}
\end{eqnarray}
\end{proposition}

\proof
We'll prove (\ref{contribest1}, \ref{contribest2}) using a similar argument as in \cite{varzuginyd}. It is not hard to see that $m(\lambda):=M_1(\lambda) M_0(\lambda)$ solves the following normalized RHP:
$$m_+(\lambda) = m_-(\lambda) J_m(\lambda,t)$$
$$J_m(\lambda,t) =  \Big((I-w^-_1)M_{0,-}\Big)^{-1}\Big((I+w^+_1)M_{0,+}\Big)$$

To show the proposition, we'll show that as $\lambda\to\infty$ nontangentially
\begin{eqnarray*}
\begin{cases}
\limsup_{\lambda \to \infty}|\lambda\big(M(\lambda)- m(\lambda)\big)| \;\;\lesssim_\epsilon\;\; t^{-\frac{3}{2(k_\theta+1)}+\epsilon}\\
\limsup_{\lambda \to \infty} |\lambda \ \text{Off}\big(m(\lambda) - M_1(\lambda)- M_0(\lambda)\big)| \;\;=\;\; 0
\end{cases}
\end{eqnarray*}
Here, for any $2\times 2$ matrix $A$, we denote by $\text{D}(A)$ the diagonal part of $A$ and $\text{Off}(A)$ the off-diagonal part of $A$. Both $\text{D}(A)$ and $\text{Off}(A)$ are $2\times 2$ matrices.

The first task is to find a suitable factorization for $J_m$ that will help us prove the first estimate. Recall that $M_{0,\pm} = \mu_0 (I \pm w_0^\pm)$, here
$$\mu_0 = I + C_{w_0}\mu_0$$

If matrix multiplication \emph{was} commutative then we would have $$(I+w^+_1)M_{0,+}=\mu_0 (I + w^+_0) (I+ w^+_1) = \mu_0(I+ w^+_1 + w^+_0)$$ thanks to disjoint support,
and similarly $(I-w^-_1)M_{0,-} = \mu_0(I- w^-_1 - w^-_0)$. However, we can only estimate $\|\mu_0-I\|_2$ by  $t^{-\frac{1}{2(k+1)}}$ which is far from the right-hand side of (\ref{contribest1}, \ref{contribest2}) if $k \sim k_\theta$. Consequently, to effectively compare $m$ with $M$ using weight comparison (i.e. the first scheme), we don't want to use the immediate choice $I \pm w^\pm_m = (I \pm w^\pm_1)M_{0,\pm}$.

Instead, we'll use the factorization $J_m=(I-w^-_m)^{-1}(I+w^+_m)$, where:
\begin{eqnarray*}
I + w_m^+ := \mu_0^{-1} (I+ w^+_1) M_{0,+}, \;\;\;
I - w_m^- := \mu_0^{-1} (I- w^-_1) M_{0,-}
\end{eqnarray*}
The idea is to somehow cancel out the effect of $\mu_0$. Now,
\begin{eqnarray*}
I + w_m^+ &=& \mu_0^{-1} (I+ w^+_1) \mu_0 (I + w^+_0)\\
              &=& (I + w^+_1)(I + w^+_0) +  \mu_0^{-1} [w^+_1,\mu_0](I + w^+_0)\\
              &=& (I+ w^+) +  \Delta w^+\\
   \Delta w^+ &:=& \mu_0^{-1} [w^+_1,\mu_0](I + w^+_0) = \mu_0^{-1}[w^+_1,\mu_0]
\end{eqnarray*}
because $w_0 = 0$ on the support of $w_1$.

Similarly, $I + w_m^+ = (I-w^-) - \Delta w^-$ where $\Delta w^- : =  \mu_0^{-1} [w^-_1,\mu_0]$. We will see that with this choice of weights for $J_m$, we can achieve $\|\Delta w\|_p \lesssim t^{-\frac{1}{k_\theta+1}}$ for any $1\leq p\leq\infty$.

From the proof of (\ref{direct_est}), $\limsup_{\lambda \to \infty}|\lambda(M (\lambda)- m(\lambda))|$ can be controlled by
\begin{eqnarray*}
&&\|\Delta w\|_2 H_2(w) + \|\mu - \mu_m\|_2\|w_m\|_2 + \limsup_{\lambda \to \infty} \big| \lambda C(\Delta w)(\lambda) \big|\\
&\equiv& \text{X} + \text{Y} + \text{Z}
\end{eqnarray*}
Here, $\|\mu-\mu_m\|_2$ can be controlled by
\begin{eqnarray}
\label{mumumest} \|C_{\Delta w}I\|_2 + \|\Delta w\|_\infty H_2(w)
\end{eqnarray}
using $\|(1-C_w)^{-1}\|_{L^2 \to L^2}\lesssim 1$ (by the previous proposition). Now by Corollary~\ref{Hardyrough},
$$H_2(w) \lesssim  t^{-\frac{1}{2(k_\theta+1)}}$$

Thus it remains to estimate $\Delta w$. For simplicity of notation, let
$$\mu^0_0(x,t):=\begin{pmatrix} 0 & \frac{u_0(t)}{x-\lambda_0}\cr \frac{v_0(t)}{x-\lambda_0} & 0\end{pmatrix} = O(t^{-\frac{1}{k+1}})$$
For $\lambda \in \text{supp}(w_1)$, we have $\text{\emph{distance}}(\lambda,\text{supp}(w_0) \cup \{ \lambda_0 \}) \gtrsim 1$, and (\ref{pwmujest}) gives
\begin{eqnarray}
\label{muNpwest1} \mu_0(\lambda)  = I + \mu^0_0(\lambda) + O_\epsilon(t^{-\frac{3}{2(k+1)}+\epsilon})
\end{eqnarray}
Since $\det M_0(\lambda) = 1$ on $\mathbb R$, we have $\det(\mu_0) \equiv 1$. Thus $\mu^{-1}_0$ is essentially a rearrangement of entries of $\mu_0$ with additional minus signs in appropriate places. Thus, (\ref{muNpwest1}) implies:
\begin{eqnarray}
\label{muNpwest2} \mu^{-1}_0(\lambda) = I - \mu^0_0(\lambda) + O_\epsilon(t^{-\frac{3}{2(k+1)}+\epsilon})
\end{eqnarray}

Now, notice that $\Delta w$ is supported on a small set, so using estimates (\ref{muNpwest1}, \ref{muNpwest2}) we have:
\begin{eqnarray*}
 \Delta w^\pm &=& (I - \mu^0_0)[w^\pm_1,I+ \mu^0_0] + O_\epsilon(t^{-\frac{3}{2(k+1)}+\epsilon})\\
              &=& [w^\pm_1,\mu^0_0] + O_\epsilon(t^{-\frac{3}{2(k+1)}+\epsilon})\\
\end{eqnarray*}
Since $w_1$ and $\Delta w$ are supported on small sets, the $O(t^{-\frac{3}{2(k+1)}+\epsilon})$ term remains small in $L^p$ for any $1\leq p\leq \infty$. In particular, $\|\Delta w\|_2, \|\Delta w\|_\infty \lesssim t^{-\frac{1}{k + 1}}$, and thus
\begin{eqnarray*}
   \text{X} &\lesssim& t^{-\frac{3}{2(k_\theta+1)}}\\
   \text{Y} &\lesssim& \|C_{\Delta w} I\|_2 + O(t^{-\frac{3}{2(k_\theta+1)}})
\end{eqnarray*}

Below we'll estimate $\|C_{\Delta w} I\|_2$. Notice that $[w^\pm_1,\mu^0_0]$ are algebraic sums of terms like $\frac{C(t)w^\pm_1(x,t)}{x-\lambda_0}$ or $\frac{w^\pm_1(x,t)C(t)}{x-\lambda_0}$ where $C(t)$ denote matrix valued functions of $t$ (independent of $x$) and have order $t^{-\frac{1}{k+1}}$. In addition, the weights $w^\pm_1$'s have the correct phase-weight relation. Since $|x-\lambda_0| \gtrsim 1$ for $x \in$ supp$(w^\pm_1)$, Corollary~\ref{Hardyrough} now gives:
$$\|C_{\Delta w}I\|_2 \lesssim |C(t)| t^{-\frac{1}{2(k_\theta+1)}} + O_\epsilon(t^{-\frac{3}{2(k_\theta+1)}+\epsilon}) = O_\epsilon (t^{-\frac{3}{2(k_\theta+1)}+\epsilon})$$
Thus, $\text{Y} \lesssim_\epsilon t^{-\frac{3}{2(k_\theta+1)}+\epsilon}$. Now to estimate $Z$ notice that $\Delta w \in L^1$, thus
$$\text{Z} = \limsup_{\lambda \to \infty} \big| \lambda \int_{-\infty}^\infty \frac{\Delta w(x,t)}{x-\lambda} dx \big| = \frac{1}{2\pi} \big| \int_{-\infty}^\infty \Delta w(x,t) dx \big|$$
so using linear theory of oscillatory integrals (see Lemma~\ref{lineartheory}) the contribution of $[w^\pm_1,\mu^0_0]$ can be controlled by $t^{-\frac{1}{k+1}-\frac{1}{k_\theta+1}}$, here $\frac{1}{k+1}$ comes from constants like $C(t)$ and $\frac{1}{k_\theta+1}$ is contributed by the respective oscillatory integral. Thus, $Z$ and hence $\limsup_{\lambda \to \infty}|\lambda(M(\lambda)- m(\lambda))|$ can be controlled by $t^{-\frac{3}{2(k_\theta+1)}+\epsilon}$.

To finish the proof of our proposition, we need to show the second estimate
\begin{eqnarray}
\label{offdiagest}    \limsup_{\lambda \to \infty} |\lambda \ \text{Off}\big(m(\lambda) - M_1(\lambda)- M_0(\lambda)\big)|  = 0
\end{eqnarray}
Observe that for any two $2\times 2$ matrices $A,B$,
\begin{eqnarray*}
    |\text{Off}(AB)|
    &=& |\text{Off}(A)\text{D}(B) + \text{D}(A)\text{Off}(B)| \\
    &\lesssim& |\text{Off}(A)| |\text{D}(B)| + |\text{D}(A)| |\text{Off}(B)|
\end{eqnarray*}
(\ref{offdiagest}) now follows by applying this inequality for $A=M_1(\lambda)-I$ and $B = M_0(\lambda)-I$, noticing
\begin{eqnarray*}
     \lim_{\lambda \to \infty}|\text{D}(M_0(\lambda)-I)| = \lim_{\lambda \to \infty} |\text{D}(M_1(\lambda)-I)| = 0\\
     \limsup_{\lambda \to \infty}|\lambda  \text{Off}(M_0(\lambda)-I)| = \limsup_{\lambda \to \infty}|\lambda  \text{Off}(M_0(\lambda))| = O_t(1) \\
     \limsup_{\lambda \to \infty}|\lambda  \text{Off}(M_1(\lambda)-I)| = \limsup_{\lambda \to \infty}|\lambda  \text{Off}(M_1(\lambda))|=O_t(1)
\end{eqnarray*}
\endproof

\subsubsection{Reduction of the a priori estimate (\ref{pwmujest}) in Lemma~\ref{apriori}} \label{aprioriphasesect}

For future reference, we'll set up a framework to reduce (\ref{pwmujest}) along our reduction of $u(t)$, $v(t)$.

Let $\mu$ and $\mu_0$ be associated with two pairs of differentiable weights $w$ and $w_0$ that are localized to a given stationary point, which we'll assume $0$ of order $k_0$. Assume that for each $x$ the triangularity of $w(x)$ and $w_0(x)$ agree, which in turn is invariant on each $D_+$ and $D_-$. For convenient, let $\Delta \mu = \mu - \mu_0$ and $\Delta w = w-w_0$.

Our setting will be as usual: the restriction on $D_+$, $D_-$ of both $\Delta w^+$ and $\Delta w^-$ are oscillatory functions of the form $fe^{it\Theta}$, where $f(0) = 0$ and $f$ has one $L^2$ derivative. Since we have already done the phase reduction, in our setting $\Theta$ is of the form $a+bx^{k_0+1}$, but the argument works for more general phases.

Recall that $P$ is away from $\{0\} \cup \text{supp}(w)$. In this localized case, we'll also have $$\text{\emph{distance}}(P, \text{supp}(w_0)) \gtrsim 1$$
In our final reduction (to model RHPs), the model weights will not be locally supported. To overcome this, we'll exploit certain analytic continuation of the model weights to deform the Riemann-Hilbert contour $\mathbb R$. Effectively, this deformation ``moves'' $P$ away from the supports of the (deformed) model weights (so that the current argument can be reused).

\begin{proposition}\label{apriorireduct} If $w$ and $w_0$ have the correct phase-weight relation and are supported away from $P$ then
$$\|\Delta \mu\|_{L^\infty(P)} \lesssim \| \Delta \mu\|_2 + t^{-\frac{1}{2(k_0+1)}}H_2(\Delta w) + t^{-\frac{2}{k_0+1}}$$
\end{proposition}

In applications of this proposition, observe that thanks to the correct phase-weight relation in $w$ and $w_0$, $H_2(\Delta w)$ can be estimated by $t^{-\frac{1}{k_0+1}+\epsilon}$ or $t^{-\frac{3}{2(k_0+1)}+\epsilon}$ depending on availability of regularity. Also, as a consequence of the reduction of $u(t)$, $v(t)$, we'll have an {\it a posteriori} estimate for $\|\Delta \mu\|_2$.  To see that, note that $\Delta \mu \equiv 0$ in the second scheme and small for large $t$ in the first scheme, using Lemma~\ref{direct_lem}. Often we'll have $\|\Delta \mu\|_2 = O(t^{-\frac{3}{2(k_0+1)}+\epsilon})$, which is enough for the desired reduction.

\proof
Writing $\Delta \mu = C_{\Delta w} \mu_0  + C_w \Delta \mu$ and using Cauchy-Schwarz's inequality in the second term, we have:
\begin{eqnarray}
\label{mupwest} \|\Delta \mu\|_{L^\infty(P)} \lesssim \|C_{\Delta w} \mu_0 \|_{L^\infty(P)} + \|\Delta \mu\|_2 \|w\|_2
\end{eqnarray}
Thus, it comes down to estimating $\big\|C_{\Delta w} \mu_0 \big\|_{L^\infty(P)}$. Intuitively, we would expect this term to be small if $\mu_0-I$ {\it was} sufficiently small (in some $L^p$), since by Lemma~\ref{lineartheory} it is not hard to see that $C_{\Delta w}I$ is small in $L^\infty(P)$. However, from Lemma~\ref{mainlemma} (more precisely, Corollary~\ref{Hardyrough}) it is only known that
$$\|\mu_0 - I \|_2 \lesssim t^{-\frac{1}{2(k_0+1))}}$$
which is not enough. To get the desired estimate, the main idea is to exploit the fact that $P$ is away from $\{0\}\cup \text{supp}(w) \cup \text{supp}(w_0)$.

Now, on $P$, write $C_{\Delta w} \mu_0$ as
\begin{eqnarray*}
C_{\Delta w} I + C_{\Delta w}\Big(C_{w_0} \mu_0\Big)
&=& O(t^{-\frac{2}{k_0+1}}) + C_+\Big(C_{w_0}\mu_0 \;\Delta w^-\Big) + C_-\Big(C_{w_0}\mu_0 \;\Delta w^+\Big)\\
&\equiv& O(t^{-\frac{2}{k_0+1}}) + \text{X} + \text{Y}
\end{eqnarray*}
Below we'll prove that $X$ is small on $P$; the proof for $Y$ is similar.

Write $X =C_+(E) + C_+(F)$, where $E \equiv C_+(\mu_0 w^-_0) \Delta w^-$ and $F\equiv C_-(\mu_0 w^+_0) \Delta w^-$. Notice that $E$ and $F$ are not symmetric in general because $C_+$ is being applied to them. However, if $C_+(X)$ is evaluated on $P$ then this symmetry is available as both $E$ and $F$ are supported away from $P$.

Keeping in mind that $\|C_+(\Delta w^-)\|_2$ and $\|C_+(\Delta w^-)\|_{L^\infty(P)}$ are small, we decompose $F$ to establish the appearance of these terms in $C_+(F)$:
\begin{eqnarray*}
 F &=&C_-(\mu_0 w^+_0) C_+(\Delta w^-)  -  C_-(\mu_0 w^+_0) C_-(\Delta w^-)\\
   &=& C_+(\mu_0 w^+_0) C_+(\Delta w^-) - \mu_0 w^+_0 C_+(\Delta w^-)  -  C_-(\mu_0 w^+_0) C_-(\Delta w^-)
\end{eqnarray*}
Under $C_+$, the last term vanishes, while the first term is unaffected. Using the distance condition, we can estimate $\|C_+(F)\|_{L^\infty(P)}$ by
\begin{eqnarray*}
&& C_+(\mu_0 w^+_0) C_+(\Delta w^-) - C_+ \Big((\mu_0-I) w^+_0 C_+(\Delta w^-)\Big) - C_+\Big(w^+_0 C_+(\Delta w^-)\Big)\\
 &=& O\Big(\|C_+(\Delta w^-)\|_{L^\infty(P)} + H_2(w_0) H_2(\Delta w)\Big)  - C_+\Big(w^+_0 C_+(\Delta w^-)\Big)
\end{eqnarray*}

For $\lambda_0\in P$, let $g(x) = \frac{w^+_0(x)}{x-\lambda_0}$. Then
\begin{eqnarray*}
C_+(w^+_0 C_+(\Delta w^-))(\lambda_0)
&=&\int g  C_+(\Delta w^-)\\
&=&  \int C_-(g) \;\;  C_+(\Delta w^-)\\
&\lesssim& \|C_-(g)\|_2 \;\; H_2(\Delta w)
\end{eqnarray*}

Consequently, using the correct phase-weight relation in both $w$ and $w_0$ we have
$$\|C_+(F)\|_{L^\infty(P)} \lesssim t^{-\frac{2}{k_0+1}} + t^{-\frac{1}{2(k_0+1)}} H_2(\Delta w)$$

\endproof

\section{Reduction to model cases (II): Deift-Zhou's steepest descent argument} \label{reductmodelsection2}

In this section, we consider a RHP that is localized to a small neighborhood of a stationary point and has a nice analytic phase. The goal of this section is to reduce this RHP to a model RHP which will be explicitly studied in Section~\ref{modelsection}. Without loss of generality we can assume that this stationary point is $0$ and is of order $k$, and the phase of this RHP is $\Theta(x) = a + bx^{k+1}$ for $a,b\in\mathbb R$, $b\neq 0$. Let $w^\pm$ denote the weights of this RHP.

Let $\delta_0(\lambda)$ be the local approximation of $\delta(\lambda)$ near $0$, and let $D_- := \{\Theta'<0\}$ and $D_+ := \{\Theta'>0\}$. Recall that $\delta_{0}$ satisfies the following scalar RHP:
\begin{eqnarray*}
   \delta_{0+}(x) = \delta_{0-}(x) \Big(1_{D_{+}}(x) + [1+p_0q_0]1_{D_{-}}(x)\Big), \ \; x\in \mathbb R
\end{eqnarray*}
here $p_0 \equiv p(0)$ and $q_0 \equiv q(0)$. Also, $D_\pm$ can be made explicit by considering parity of $k$ and sign of $b$; each of them is either $\emptyset$, $\mathbb R_-$, $\mathbb R_+$, or $\mathbb R\setminus \{0\}$.

The current pair of weights $w=(w^-, w^+)$ is of the form
\begin{eqnarray*}
\label{weightwj} (w^-, w^+)
&=& \begin{cases}
             \Bigl(\begin{pmatrix} 0 & \phi \delta^2_-p e^{-it\Theta} \cr 0 & 0 \end{pmatrix}, \begin{pmatrix} 0 & 0 \cr \phi \delta^{-2}_+q e^{it\Theta}& 0   \end{pmatrix}\Bigr), & \text{if $x \in D_+$;}\\
             \Bigl(\begin{pmatrix} 0 &  0\cr \phi \delta^{-2}_-\frac{q}{1+pq} e^{it\Theta}& 0 \end{pmatrix}, \begin{pmatrix} 0 & \phi \delta^2_+\frac{p}{1+pq}e^{-it\Theta} \cr 0 & 0   \end{pmatrix}\Bigr), & \text{if $x \in D_-$}
    \end{cases}
\end{eqnarray*}
here $\phi$ is a normalized cutoff supported near $0$. Note that we implicitly used the jump relation of $\delta$ to write $w$ in the above form. $w$ will be reduced to:
\begin{eqnarray*}
\label{weightwMj}  (w^-_M,w^+_M) &=&
 \begin{cases}
             \Bigl(\begin{pmatrix} 0 & \delta^2_{0-}p_0 e^{-it\Theta} \cr 0 & 0 \end{pmatrix}, \begin{pmatrix} 0 & 0 \cr  \delta^{-2}_{0+}q_0 e^{it\Theta}& 0   \end{pmatrix}\Bigr), & \text{if $x \in D_+$;}\\
             \Bigl(\begin{pmatrix} 0 &  0\cr \delta^{-2}_{0-}\frac{q_0}{1+p_0q_0} e^{it\Theta}& 0 \end{pmatrix}, \begin{pmatrix} 0 & \delta^2_{0+} \frac{p_0}{1+p_0q_0}e^{-it\Theta}\cr 0 & 0   \end{pmatrix}\Bigr), & \text{if $x \in D_-$}
    \end{cases}
\end{eqnarray*}
Intuitively, this means that as $t\to\infty$ the RHP $(w^-, w^+)$ localizes at $0$: $\phi(x) \to 1$, $p(x) \to p_0$, $q(x) \to q_0$, and $\delta(x) \to \delta_{0}(x)$.

The main difficulty in the above reduction is the lack of $L^2$ integrability of the model weights $w_M$, despite the fact that they are still in $L^\infty$. This prevents Beals-Coifman's operator formulation which is essential to our perturbation schemes. To get around this issue, the idea introduced in \cite{deiftzhouMKdVyd} is to exploit analytic continuation of $w_M$ and deform the model RHP to a suitably chosen contour $\Gamma$, on which $w_M$ has strong decay. We note that this approach was not taken by \cite{varzuginyd}.

We'll use the following adaptation of a notation in \cite{deiftzhouNLSyd}: For any continuous $f$, let
\begin{eqnarray*}
   [f](x) = \frac{f(0)}{1+ix^N}
\end{eqnarray*}
here $N$ is a large natural number that depends only on $k$, so that $[f]$ decays fast enough on $\mathbb R$. Large choice of $N$ in particular ensures the applicability  to $[f]$ of Lemma~\ref{localprpl} and related propositions. Intuitively, $[f]$ is an analytic approximation at $0$ of $f$ that stays in $L^p$ spaces.

For simplicity of notation, let $\mathbb R_+$ denote $[0,\infty)$ and $\mathbb R_-$ denote $(-\infty,0]$. Consider $\Gamma_0,\Gamma_1,\Gamma_2,\Gamma_3,\Gamma_4,\Gamma_5$ six rays originating at $0$ such that $\Gamma_0 = \mathbb R_+$, $\Gamma_3 = \mathbb R_-$ and the others form small angles with $\mathbb R$. More specifically, we will take
$$\Gamma_1 = e^{i\alpha}\mathbb R_+, \; \Gamma_5 = e^{-i\alpha}\mathbb R_+,$$
$$\Gamma_2 = e^{-i\alpha}\mathbb R_-, \;\Gamma_4 = e^{i\alpha}\mathbb R_-$$
for a sufficiently small $\alpha>0$, say $\alpha = \frac{\pi}{3N}$. Small value of $\alpha$ ensures that $[f](x)$ is analytic and decays strongly in the two angles formed by $\Gamma_5, \Gamma_1$ and $\Gamma_2,\Gamma_4$.

For any (reasonable) oriented contour $\Sigma$, we can define (nontangential) $\pm$ boundary values on $\Sigma$ of a function analytic on $\mathbb C\setminus \Sigma$ using the following standard convention: Along $\Sigma$ following the local direction, the left side is $+$ and the right side is $-$. Consequently, for any $f\in L^p(\Sigma)$, $1<p<\infty$, we can define the nontangential boundary values $C^\pm_\Sigma f$ of the Cauchy operator $C_\Sigma f$.

Let $\Gamma = \bigcup_{i=0}^5 \Gamma_i$. On $\Gamma$, we'll orient $\Gamma_0,\Gamma_2,\Gamma_4$ outwards (i.e. $\to \infty$), and $\Gamma_1,\Gamma_3,\Gamma_5$ inwards (i.e. $\to 0$). In particular, $\mathbb R\subset \Gamma$ is oriented naturally. This orientation makes $\Gamma$ a \emph{complete} contour in the sense that it divides $\mathbb C$ into two (disconnected) regions: one stays entirely on the $-$ side of $\Gamma$, and the other stays entirely on the $+$ side of $\Gamma$ (it is equivalent to saying that $\Gamma$ is the common boundary of these regions, clockwise for one and counter-clockwise for the other) \cite{Durenyd}. For convenience, we'll use $-\Gamma_j$ to refer to the reverse orientation on $\Gamma_j$, and $|\Gamma_j|$ whenever orientation is not taken in account.

The main advantage of $\Gamma$ being \emph{complete} is: the Cauchy operators $C^\pm_\Gamma$ satisfy standard properties of $C^\pm_{\mathbb R}$ (see for instance, Zhou \cite{zhouyd}):
$$C^+_\Gamma - C^-_\Gamma = I, \; C^-_\Gamma C^+_\Gamma = C^+_\Gamma C^-_\Gamma = 0$$

For any $i\neq j$, let $\Gamma_{ij} \subset \mathbb R^2$ denote the angle formed by rotating $|\Gamma_i|$ {\it counter-clockwise} to $|\Gamma_j|$. With this notation, $\Gamma_{01} \cup \Gamma_{23}\cup \Gamma_{45}$ is the $+$ side of $\Gamma$, and $\Gamma_{12} \cup \Gamma_{34}\cup \Gamma_{50}$ is the $-$ side of $\Gamma$.

The following observation is crucial in the analysis of this section: $w^+_M$ and $w^-_M$ have analytic continuation to small angles (pivoted at $0$) respectively in the upper and lower half planes, where they decay strongly. More precisely, $w^+_M$ can be nicely continued to $\{z\in \mathbb C: \text{Arg}z \in (0,\frac{\pi}{k+1}) \cup (\pi - \frac{\pi}{k+1}, \pi)\}$
and $w^-_M$ can be nicely continued to $\{z\in \mathbb C: \text{Arg}z \in (-\frac{\pi}{k+1},0) \cup (-\pi, -\pi + \frac{\pi}{k+1})\}$. To verify these claims, it may be convenience to consider different cases based on the parity of $k$ and the sign of $\Theta^{(k+1)}(0)$.

Intuitively, the above strong decay of the deformed model weights is not a coincidence, it is rather expected. If instead of real-variable methods, Lemma~\ref{mainlemma} and Lemma~\ref{localprpl} \emph{were proved} using a steepest descent argument (assuming relevant things were analytic) then this type of decay would be exactly what we need. Given a complete contour $\Sigma$, in order to have $C^+_\Sigma(w^-)$ to be small it would be convenient that $w^-$ has a decaying analytic continuation to the $-$ side of $\Sigma$. Similarly, for $C^-_\Sigma(w^+)$ to be small we want $w^+$ to have a decaying analytic continuation to the $+$ side of the contour. Whenever a pair of weights has these properties, we'll say that they are {\it natural}. The naturality of our weights are indeed guaranteed by the correct phase-weight behavior in $(w^-, w^+)$, which is a consequence of the studies in Section~\ref{hardyspacesection}.

The reduction in this section consists of two major steps:

\underline{\emph{Step 1:  Preparation for steepest descent.}} In this step, we'll reduce $w \to \widetilde w_M$, an analytic approximation of $w_M$ that stays in $L^2(\mathbb R)$:
\begin{eqnarray*}
 \begin{cases}
             \Bigl(\begin{pmatrix} 0 & \delta^2_-[p]e^{-it\Theta} \cr 0 & 0 \end{pmatrix}, \begin{pmatrix} 0 & 0 \cr  \delta^{-2}_+[q] e^{it\Theta}& 0   \end{pmatrix}\Bigr), & \text{if $x \in D_{0+}$;}\\
             \Bigl(\begin{pmatrix} 0 &  0\cr \delta^{-2}_-[\frac{q}{1+pq}] e^{it\Theta}& 0 \end{pmatrix}, \begin{pmatrix} 0 & \delta^2_+[\frac{p}{1+pq}]e^{-it\Theta}\cr 0 & 0   \end{pmatrix}\Bigr), & \text{if $x \in D_{0-}$}
    \end{cases}
\end{eqnarray*}

\underline{\emph{Step 2: Steepest descent.}} In this step, we'll reduce $\widetilde w_M$ to the $\Gamma$-deformation of the model RHP associated with $w_M$. To do this we'll deform the former to $\Gamma$. In this reduction, we largely follow the steepest descent method of Deift and Zhou \cite{deiftzhouMKdVyd,deiftzhouNLSyd}, with several adaptations to our settings.

\subsection{Preparation for steepest descent}

Recall that the set $P$ used in (\ref{pwmujest}) is of distance $\gtrsim 1$ from $0$. Without loss of generality, we can assume that
$$distance(P,\text{supp}(\phi))\gtrsim 1$$

The reduction is broken up into two sub-steps
$$w \to \phi \widetilde w_M \to \widetilde w_M$$
In the second sub-step, the localization argument as in Section~\ref{nbhreductsect} can be used to control the effect on $u(t)$, $v(t)$ by terms decaying highly as $t\to\infty$; hence the only task is to prove the reduction of (\ref{pwmujest}); observe that the same argument as in Section~\ref{aprioriphasesect} can not be repeated because $\widetilde w_M$ is not localized.

\noindent \underline {\it Sub-step 1:} ($w \to \phi \widetilde w_M$). Let $\Delta w = w - \phi \widetilde w_M$. Using the two perturbation schemes described in Section~\ref{localizationsection} and implemented in Section~\ref{reductmodelsection}, the following estimates are enough for our reduction:

\begin{proposition} (i) If $p$ and $q$ has two $L^2$ derivatives near $0$ then for any sign combinations
\begin{eqnarray}
\label{wMLp} H_p(1_{D_{\pm}}\Delta w) &\lesssim_p& t^{-(1+\frac{1}{p})\frac{1}{k+1}} \ln t,\;\;\;2\leq p < \infty \\
\label{wMLinfty} H_\infty(1_{D_{\pm}}\Delta w) &\lesssim_\epsilon& t^{-\frac{1}{k+1}+\epsilon}  \;  \forall \epsilon>0\\
\label{wMpw} \limsup_{\lambda \to \infty} |\lambda C \Big(1_{D_{\pm}}\Delta w^\pm \Big)| &\lesssim& t^{-\frac{2}{k+1}}
\end{eqnarray}

(ii) If $p$ and $q$ has only one $L^2$ derivative near $0$ then (\ref{wMpw}) remains true, and
\begin{eqnarray}
\label{wMLp2} H_p(1_{D_{\pm}}\Delta w) &\lesssim_p& \max(t^{-\frac{1}{p}}, t^{-(\frac{1}{2}+\frac{1}{p})\frac{1}{k+1}}),\;\;\;2\leq p < \infty \\
\label{wMLinfty2} H_\infty(1_{D_{\pm}}\Delta w) &\lesssim& O(1) + O(t^{-\frac{1}{2(k+1)}+\epsilon})
\end{eqnarray}
in (\ref{wMLinfty2}) the $O(1)$ term is (modulo some positive power) proportional to the size of $supp(\phi)$.
\end{proposition}

\proof The above estimates are indeed direct consequences of Lemma~\ref{mainlemma} (in particular Corollary~\ref{Hardyvanish}) and Lemma~\ref{lineartheory}. For instance, we'll consider the setting in (i). Each $1_{D_{\pm}}\Delta w^\pm$ has at most one nonzero entry $1_{D_{\pm}}\Omega f e^{it\theta}$, with:

(i) correct monotonicity for $\theta$ on the support of this entry;

(ii) $f$ vanishes at $0$ with multiplicity $1$ up to the second derivative, and

(iii) $\Omega$ has an $A_2+B_2$ decomposition.

\noindent For example, on $D_{+}$ the nontrivial entry of $\Delta w^+_M$ is
$$(\phi \delta^{-2}_+ q -\phi \delta^{-2}_+[q])e^{it\Theta} \equiv \Omega^{-1}_0 \phi \Big(q - [q]\Big)e^{it\Theta}$$
The argument is similar for the setting in (ii), note that the corresponding multiplicity will be $\frac{1}{2}$.
\endproof

Using the above proposition, the effect on $u(t),v(t)$ can be controlled by $O_\epsilon(t^{-\frac{2}{k+1}+\epsilon})$ if $p$ and $q$ have two $L^2$ derivatives near $0$, or $O_\epsilon(t^{-\frac{3}{2(k+1)}+\epsilon})$ if they have only one $L^2$ derivative. In the latter case, choosing $supp(\phi)$ small enough will ensure that the reduction doesn't harm the unique solvability and the corresponding resolvent bound of the RHP associated with $w$ (see Section~\ref{reductphasesect}); alternatively we can perform a suitable approximation as in Corollary~\ref{coruniqsolv}.

Regarding the a priori estimate (\ref{pwmujest}), we'll use the argument of Section~\ref{aprioriphasesect}. Note that this a priori estimate is not needed if $N=1$ (thus one $L^2$ derivative for $p, q$ near $0$ is enough). When $N>1$, for the reduction of (\ref{pwmujest}) we want $\|\Delta \mu\|_2  = t^{-\frac{3}{2(k+1)}}$, and this is satisfied if $p$ and $q$ have two $L^2$ derivatives near $0$.

\noindent \underline{\it Sub-step 2:} ($\phi \widetilde w_M$ to $\widetilde w_M$). The only thing we have to show is the reduction of (\ref{pwmujest}). As usual, $P$ denotes a subset of $\mathbb R$ that is away from $0$ and supp$(\phi)$, and $\lambda_0$ is a generic point in $P$. The main idea is to exploit analytic continuation, which was not available for us in Section~\ref{aprioriphasesect}, to move the weights away from $\mathbb R$. Effectively, this re-establishes the distance condition between $P$ and the support of the weights, which was essential in our reduction of (\ref{pwmujest}) in the localized case.

\begin{proposition}\label{apriorisd} If (\ref{pwmujest}) is true for $\mu$ (of $\widetilde w_M$) then it is true for $\mu$ (of $\phi \widetilde w_M$).
\end{proposition}

\proof We'll use a variant of a steepest descent argument in \cite{varzuginyd} to show this proposition.

For simplicity of notation, in this proof the pair of weights $\widetilde w_M$ will be simply denoted as $w$. Let $\mu$ be associated with $w$, and the solution of this RHP is denoted by $M$. We'll use $\Delta \mu, \Delta w^\pm$ to denote the respective differences of two RHPs.

Since $\text{\emph{distance}}(P,\text{supp}(\phi w)) \gtrsim 1$, Cauchy-Schwarz's inequality gives
\begin{eqnarray*}
\|\Delta \mu\|_{L^\infty(P)}
&\lesssim& \|C_{\Delta w} \mu \|_{L^\infty(P)} + \|\Delta \mu\|_2 \|\phi w\|_2 = \|C_{\Delta w} \mu \|_{L^\infty(P)} + O(t^{-\frac{3}{2(k+1)}})
\end{eqnarray*}
To estimate $C_{\Delta w} \mu$ on $P$, we'll write it as $C_+(\mu \Delta w^-) + C_-(\mu \Delta w^+) \equiv \text{X} + \text{Y}$. Below we'll prove that $\|\text{X}\|_{L^\infty(P)}$ is small; $\|\text{Y}\|_{L^\infty(P)}$ can be estimated similarly.

Using $\mu= M_+(I+ w^+)^{-1} = M_-(I-w^-)^{-1}$, an explicit computation gives
\begin{eqnarray}
\label{mu11D+} 1_{D_+}\mu_{11}  &= 1_{D_+}M_{11-} &= 1_{D_+}\Big( M_{11+} - \delta^{-2}_+[q(1+pq]e^{it\Theta} M_{12+}\Big)\\
\label{mu11D-} 1_{D_-}\mu_{11}  &= 1_{D_-} M_{11+} &= 1_{D_-}\Big( M_{11-} + \delta^{-2}_-[q/(1+pq]e^{it\Theta}M_{12-}\Big)
\end{eqnarray}
Thus, $\mu_{11} 1_{D_+}$ has analytic continuation to a \emph{good} sector (with pivot at $0$) on which it is normalized to $1$ as $z\to\infty$ nontangentially. Similar observations can be made for $\mu_{11}1_{D_-}$.

Without loss of generality, we can assume that $\Theta^{(k+1)}(0)>0$. The picture (and our argument) is reflected across the imaginary axis if $\Theta^{(k+1)}(0)<0$.

Under this assumption, the continuation sector can be made precise as follows:

(i) $k$ is odd: $D_+ = \mathbb R_+, D_- = \mathbb R_-$, and a good continuation sector for $\mu_{11}1_{D_+}$ is $\{ -\pi< \text{Arg}z < \frac{\pi}{2N} \}$, and a good continuation sector for $\mu_{11}1_{D_-}$ is $\{0 <  \text{Arg} z < \pi+\frac{\pi}{2N} \}$.

(ii) $k$ is even: $D_+ = \mathbb R\setminus \{0\}$, $D_-=\emptyset$, and a good continuation sector for $\mu_{11}1_{D_+}$ is $\{ -\pi - \frac{\pi}{2N}< \text{Arg}z < \frac{\pi}{2N} \}$.

Repeating this computation, we'll deduce similar conclusions for every entry of $\mu$ (one small difference: the diagonal entries converge to $1$ in their continuation sectors, while the off-diagonal entries vanish at $\infty$).

The computation also reveals that {\it (the restrictions to $D_\pm$ of) different entries of $\mu$ do not always have the same continuation sectors}. Since our estimation of $X$ involves deformation of $\mathbb R$ by nontrivial angles, the knowledge of the continuation sector is important. This is the reason why we'll need to estimate $X$ {\it entry-wise}, unlike our previous argument for the localized case.

An explicit computation, with triangularity taken into account, gives
$$\mu \Delta w^- =
\begin{pmatrix}
     \mu_{12} \Delta w^-_{21}1_{D_-}  & \mu_{11} \Delta w^-_{12}1_{D_+} \cr
     \mu_{22} \Delta w^-_{21}1_{D_-}  & \mu_{21} \Delta w^-_{12}1_{D_+}
\end{pmatrix}$$

Repeating the above computation (for $\mu_{11}$) on other entries of $\mu$, it is not hard to see that: in {\it every entry} of $\mu \Delta w^-$, the respective restriction of the relevant entry of $\mu$ has a {\it good continuation sector that contains the lower half plane}.

Intuitively, when estimating $C_+(\mu \Delta w^-)$, this observation allows us to reduce the respective $\Delta w^-_{ij} 1_{D_\pm}$ to their Hardy projections onto $\mathbb C_+$ (i.e. $C_+(\Delta w^-_{ij}1_{D_\pm})$) which are known to be small. Combining this observation with a contour deformation argument, we will be able to show the desired estimate for $X$.

Below we'll estimate $X_{12}$; detailed computation for other entries can be done similarly. Let $f$ be the analytic continuation of $\mu_{11}1_{D_+}$ mentioned above. Then $f$ agrees with $M_{11-} \in 1 + H_2(\mathbb R_-)$ on $\{\text{Im}(z)<0\}$, therefore
\begin{eqnarray*}
\text{X}_{12} &=& C_+(f \Delta w^-_{12}1_{D_+})\\
              &=& C_+(\Delta w^-_{12}1_{D_+}) + C_+\Big((f-1) \Delta w^-_{12}1_{D_+}\Big)\\
              &=& O(t^{-\frac{2}{k+1}}) +  C_+\Big((f-1) C_+(\Delta w^-_{12}1_{D_+})\Big)
\end{eqnarray*}
Let $g = (f-1) C_+(\Delta w^-_{12}1_{D_+})$. Consider two cases:

\noindent {\bf Case 1:} $k$ is odd. Then $D_-=\mathbb R_-$ and $D_+=\mathbb R_+$.

Define $\Gamma^0_+ = \Gamma_1$, $\Gamma^0_- = \Gamma_4$, and $\Gamma^0 = \Gamma^0_+ \cup \Gamma^0_-$, orientated by increasing order of the real part.  Notice now $g$ has analytic continuation from $D_+$ to $\{0 <\text{Arg} z <\frac{\pi}{2N}\}$ and from $D_-$ to $\{-\pi <\text{Arg} z <-\pi + \frac{\pi}{2N}\}$ (observe that $C_+(\Delta w^-_{12}1_{D_+})$ is analytic outside $D_+$). We'll show that:

\noindent \underline{\emph{Claim:}} For $\lambda$ in the upper half plane,
\begin{eqnarray*}
\label{contourinteg}
\frac{1}{2\pi i}\int_{\mathbb R}\frac{g(z)}{z-\lambda} dz &=& \frac{1}{2\pi i}\int_{\Gamma^0} \frac{g(z)}{z-\lambda} dz +
\begin{cases}
g(\lambda) & \text{\emph{if $0<\text{Arg}\lambda <\alpha$;}}\\
0 & \text{\emph{if $\alpha<\text{Arg}\lambda < \pi$.}}
\end{cases}
\end{eqnarray*}

\proof[Proof of claim] We'll use contour integration (twice). For $R>0$ large, consider a triangular contour $\gamma$ consists of three edges: $\gamma_1:=[0, R]$, $\gamma_2 := e^{i\alpha}[0,R]$, and $\gamma_3$ straight line segment connecting $R$ and $Re^{i\alpha}$. This contour is oriented counter clockwise, and
$$\frac{1}{2\pi i}\int_\gamma \frac{g(z)}{z-\lambda} dz = \begin{cases}
g(\lambda) & \text{if $0<\text{Arg}(\lambda)<\alpha$;}\\
0 & \text{if $\alpha<\text{Arg}(\lambda)<\pi$.}
\end{cases}$$

We'll show that the contribution of $\gamma_3$ vanishes as $R\to \infty$. Indeed, a simple application of Carleson's measure theorem (see for instance \cite{garnettbafyd}) shows that: for any $0\leq \alpha\leq \pi$ and $h \in L_2(\mathbb R)$, and for $\Gamma_R:=\{R+e^{i\alpha}x: x\in [0,\infty)\}$ (which contains $\gamma_3$):
$$\|C_+(h)\|_{L^2(\Gamma_R)} \lesssim \|C_+(h)\|_{H_2(\mathbb C_+)} \sim \|C_+(h)\|_{L^2(\mathbb R)} $$
Consequently, by (\ref{mu11D+}), by the boundedness of $\delta$, and by the fact that the respective oscillating phase of $\widetilde w^+$ is bounded (indeed it decays) on $\gamma_3$, we have,
$$\|f-1\|_{L^2(\gamma_3)} \lesssim \|\widetilde M_{11+}\|_{L^2(\mathbb R)} + \|\widetilde M_{12+}\|_{L^2(\mathbb R)} \lesssim 1$$
Cauchy-Schwarz's inequality now gives
$$\|g\|_{L^1(\gamma_3)} \lesssim \|f-1\|_{L^2(\gamma_3)} \|C_+(\Delta^-_{12}1_{D_+})\|_{L^2(\gamma_3)} \lesssim 1$$
so as $R\to\infty$, the contribution of $\gamma_3$ becomes zero, and we get
\begin{eqnarray*}
\frac{1}{2\pi i}\int_0^\infty \frac{g(z)}{z-\lambda} dz = \frac{1}{2\pi i}\int_{\Gamma^0_+} \frac{g(z)}{z-\lambda} dz +
\frac{1}{2\pi i}\int_\gamma \frac{g(z)}{z-\lambda} dz
\end{eqnarray*}
Similar argument shows that, for $\lambda$ in the upper half plane,
\begin{eqnarray*}
\frac{1}{2\pi i}\int_{-\infty}^0 \frac{g(z)}{z-\lambda} dz &=& \frac{1}{2\pi i}\int_{\Gamma^0_-} \frac{g(z)}{z-\lambda} dz
\end{eqnarray*}
Adding these two equalities, we get the desired claim. \endproof

\noindent Taking limit of (\ref{contourinteg}) from the upper half plane to any $\lambda_0 \in P$, we have
$$\text{X}_{12}(\lambda_0) = O(t^{-\frac{2}{k+1}}) +  g(\lambda_0)1_{\mathbb R_+} + \frac{1}{2\pi i}\int_{\Gamma^0} \frac{g(z)}{z-\lambda_0} dz$$
Now, $|\mu(\lambda_0) -I| \lesssim t^{-\frac{1}{k+1}}$ since $\mu$ satisfies the a priori estimate (\ref{pwmujest}). Thus, using the definition of $g$, we have
\begin{eqnarray*}
g(\lambda_0)1_{\mathbb R_+}
  &\lesssim& t^{-\frac{1}{k+1}} \|C_+(\Delta w^-_{12}1_{D_+})\|_{L^\infty(P)}\\
  &\lesssim& t^{-\frac{2}{k+1}}
\end{eqnarray*}
On the other hand, by Cauchy-Schwarz's inequality,
\begin{eqnarray*}
\frac{1}{2\pi i}\int_{\Gamma^0} \frac{g(z)}{z-\lambda_0} dz &\lesssim&
\frac{1}{|\lambda_0|}\|f-1\|_{L^2(\Gamma^0)} \|C_+(\Delta w^-_{12}1_{D_+})\|_{L^2(\Gamma^0)}\\
&\lesssim& \|C_+(\Delta w^-_{12}1_{D_+})\|_{L^2(\Gamma^0)}
\end{eqnarray*}
Applying Lemma~\ref{variantMainLemma}, we have
$$\|C_+(\Delta w^-_{12}1_{D_+})\|_{L^2(\Gamma^0)}\lesssim t^{-\frac{3}{2(k+1)}}$$
In fact, write $\Delta w^-_{12}1_{D_+} = h(x) \delta_-^2 e^{-it\Theta(x)} 1_{D_+}$ with
$$h = (1-\phi) [p]$$
Note that $\delta^2_-$ has the $A_1+B_1$ decomposition, $1-\phi$ is smooth and supported away from $0$, thus $h$ vanishes at $0$ with multiplicity $1$ up to the first derivative. Now Lemma~\ref{variantMainLemma} gives us the above estimate.

\noindent {\bf Case 2:} $k$ is even. Then $D_-=\emptyset$ and $D_+=\mathbb R\setminus \{0\}$.

In this case, define $\Gamma^0_+ :=\Gamma_1$, $\Gamma^0_- :=\Gamma_2$. The rest of the argument is essentially the same. Similar contour integration argument gives
\begin{eqnarray*}
\frac{1}{2\pi i}\int_{\mathbb R}\frac{g(z)}{z-\lambda} dz = \frac{1}{2\pi i}\int_{\Gamma^0} \frac{g(z)}{z-\lambda} dz +
\begin{cases}
g(\lambda), & \text{\emph{$0<\text{\emph{Arg}} \lambda <\alpha$;}}\\
g(\lambda), & \text{\emph{$\pi-\alpha<\text{\emph{Arg}}\lambda <\pi$;}}\\
0, & \text{\emph{$\alpha <\text{\emph{Arg}} \lambda <\pi-\alpha$.}}
\end{cases}
\end{eqnarray*}

\noindent From there, for $\lambda_0 \in P$ we have
\begin{eqnarray*}
\text{X}_{12}(\lambda_0)
     &=& O(t^{-\frac{2}{k+1}}) +  g(\lambda_0) + \frac{1}{2\pi i}\int_{\Gamma^0} \frac{g(z)}{z-\lambda_0} dz\\
\implies   \text{X}_{12}(\lambda_0)
     &=& O(t^{-\frac{2}{k+1}}) +  O(\|C_+(\Delta w^-_{12}1_{D_+})\|_{L^2(\Gamma^0)})
\end{eqnarray*}
By Lemma~\ref{variantMainLemma}, $\|C_+(\Delta w^-_{12}1_{D_+})\|_{L^2(\Gamma^0)}$ is small. Eventually $X_{12}(\lambda_0)$ is small.
\endproof

\subsection{Steepest descent reduction to the $\Gamma$-deformed model case} The reduction consists of two steps:
\begin{eqnarray*}
&\text{Deformation of the pre-model RHP to $\Gamma$}&\\
&\downarrow&\\
&\text{Reduction to the $\Gamma$-deformed model RHP}&
\end{eqnarray*}
However, regarding unique solvability and resolvent bounds, we can reduce a step further to the model RHP $w_M$ on $\mathbb R$ (the undeformed model RHP). The lack of $L^2(\mathbb R)$ integrability of $w_M$ is the main reason that prevents the reduction to this actual model RHP of the leading asymptotics of $u,v$ and the a priori estimate (\ref{pwmujest}). On the other hand, we can still talk about invertibility of $1-C_{w_M}$ because $w^\pm_M$ remains in $L^\infty(\mathbb R)$.

The following analogue is very useful to keep in mind:
$$\text{contour deformation} \;\; \sim \;\;\text{second perturbation scheme},$$
$$\text{reduction on $\Gamma$}\;\; \sim \;\;\text{first perturbation scheme}.$$
We'll see many similarities which support this intuition, one of which was already mentioned above: in the first perturbation scheme and in the reduction on $\Gamma$, we estimate the effect on $u(t), v(t)$ by directly estimating the weight differences.

\noindent \underline{\bf I. Contour deformation}

In this section, for simplicity the weights associated with the pre-model RHP will be denoted as $(w^-,w^+)$:
\begin{eqnarray*}
(w^-,w^+) = \begin{cases}
             \Bigl(\begin{pmatrix} 0 & \delta^2_-[p]e^{-it\Theta} \cr 0 & 0 \end{pmatrix}, \begin{pmatrix} 0 & 0 \cr  \delta^{-2}_+[q] e^{it\Theta}& 0   \end{pmatrix}\Bigr), & \text{if $x \in D_{+}$;}\\
             \Bigl(\begin{pmatrix} 0 &  0\cr \delta^{-2}_-[\frac{q}{1+pq}] e^{it\Theta}& 0 \end{pmatrix}, \begin{pmatrix} 0 & \delta^2_+[\frac{p}{1+pq}]e^{-it\Theta}\cr 0 & 0   \end{pmatrix}\Bigr), & \text{if $x \in D_{-}$}
    \end{cases}
\end{eqnarray*}
Here, $\delta$ is the solution to the scalar RHP defined in Section~\ref{factorizesection}, $\Theta(x) = a+bx^{k+1}$ is a real valued phase, and $D_+ = \{\Theta >0\}$, $D_-=\{\Theta <0\}$.

Our goal in this section is to deform this RHP to the following RHP:
\begin{eqnarray*}
m_{+}(z) &=& m_{-}(z) j_\Gamma(z), \;\;\; z \in \Gamma\\
m_+ - I \in H_2(\Gamma_+), && m_- - I \in H_2(\Gamma_-)
\end{eqnarray*}
The jump matrix $j_\Gamma$ is defined on $\Gamma$ and equals to

(i) $I$ on $\Gamma_0 \cup \Gamma_3$; and

(ii) $(I+w^+)^{-1}$ on $\Gamma_1 \cup \Gamma_2$; and

(iii) $I-w^-$ on $\Gamma_4 \cup \Gamma_5$.

\noindent Here $w^\pm$ on $\Gamma$ should be understood as the analytic continuation of the respective restriction of these weights. For instance, on $\Gamma_1$, $w^+$ is the continuation of $w^+|_{\mathbb R_+}$.

Indeed, if $n$ solves the pre-model RHP ($w^-,w^+$) then $m(z):=n(z)(I+\Phi(z))^{-1}$ solves the above RHP, where:
\begin{eqnarray*}
   \Phi(z) &=& \begin{cases}
                   0, & \text{ if $z\in \Gamma_{12} \cup \Gamma_{45}$;}\\
                   w^+, & \text{ if $z\in \Gamma_{01}\cup \Gamma_{23}$;}\\
                   -w^-,  & \text{ if $z\in \Gamma_{50}\cup \Gamma_{34}$.}
              \end{cases}
\end{eqnarray*}
(as usual, $t$ is suppressed for simplicity). To see that $m$ satisfies the $L^2$ normalization condition, it might be convenient to exploit strong decay of the above $\Phi$ (which goes back to strong decay of the respective oscillating terms $e^{\pm it\Theta}$). However, as we'll see, uniform boundedness of $\Phi$ is good enough for this purpose.

Since $\Gamma$ is complete, the same theory on $\mathbb R$ is applicable: For a $L^2\cap L^\infty$ factorization $j_\Gamma = (I-w^-_\Gamma)^{-1} (I+ w^+_\Gamma)$, as long as $(1-C_{w_\Gamma})^{-1}$ exists on $L^2(\Gamma)$, we can find a function $\mu_\Gamma$ on $\Gamma$ such that $\mu_\Gamma = I + C_{w_\Gamma} \mu_\Gamma$, and the above RHP has the following unique $L^2$ normalized solution:
$$m_{\pm}(\lambda):=I + C^\pm_{\Gamma}(\mu_\Gamma(w^+_\Gamma + w^-_\Gamma))(\lambda) , \; \lambda \not\in \Gamma$$
The following relations also hold $\forall z\in \Gamma$:
$$m_{+}(z) = \mu_\Gamma (I + w^+_\Gamma), \;\; m_{-}(z) = \mu_\Gamma (I - w^-_\Gamma)$$

Since $w^-, w^+ \in L^2 \cap L^1$, we can recover $u$, $v$ from the limit as $z\to\infty$ of $z n(z)$ along {\it any} non-tangential direction in $\mathbb C\setminus \mathbb R$. Our choice of factorization will ensure that $w^\pm_\Gamma \in L^\infty\cap L^2\cap L^1$, so we can recover $u_\Gamma(t), v_\Gamma(t)$ by taking limit of $z m(z)$ as $z\to\infty$ along {\it any} non-tangential direction in $\mathbb C\setminus \Gamma$. Since $m(z)$ and $n(z)$ agree inside $\Gamma_{12} \cup \Gamma_{45}$, we get:

\begin{corollary} Assuming unique solvability (of both RHPs), the deformation doesn't change the recovered potentials $u,v$ of the RHP associated with $(w^-,w^+)$.
\end{corollary}

\noindent {\it Remarks}: 1. The above deformation can be carried out for the model RHP (\ref{weightwMj}).

2. Our second scheme is based on a conjugation of the jump matrix, which corresponds to a right multiplication on each Riemann-Hilbert factor. This is algebraically similar to the above contour deformation; the main difference is that the deformed and original RHPs are on different contours. Thus, contour deformation can be seen as \emph{a contour extension followed by the second scheme}.

To enable contour extension, we'll show that the invertibility of the respective $1-C_w$ are not affected when extending the weights to trivially to a larger contour; furthermore the old and new resolvent norms are comparable.

\subsubsection{Contour extension}

Our setting is general: $\Sigma = \Sigma_1 \cup \Sigma_2$ is a union of two oriented contours (intersecting at finitely many points);  $w = (w^+,w^-)$ is a pair of $L^\infty$ weights defined on $\Sigma$ such that $w|_{\Sigma_2} \equiv 0$; $C_{\Sigma,w}$ and $C_{\Sigma_1,w}$ are Beals-Coifman operators on $\Sigma$ and $\Sigma_1$ respectively. Note that we do not require the weights to be in $L^2$; this allows us to apply the following proposition to both model and pre-model weights.

\begin{proposition}[Contour extension] \label{extendcontour} Let $1<p<\infty$. The invertibility of $1-C_{\Sigma_1,w}$ and $1-C_{\Sigma,w}$ respectively on $L^p(\Sigma_1)$ and $L^p(\Sigma)$ are equivalent, furthermore the resolvent norms are ``comparable" in the sense that the followings are equivalent:
$$\|(1-C_{\Sigma,w})^{-1}\|_{L^p(\Sigma) \to L^p(\Sigma)} \lesssim 1$$
$$\|(1-C_{\Sigma_1,w})^{-1}\|_{L^p(\Sigma_1) \to L^p(\Sigma_1)} \lesssim 1$$
\end{proposition}
\proof ($\Rightarrow$): Suppose that $1-C_{\Sigma,w}$ is invertible on $L^p(\Sigma)$. Let $f\in L^p(\Sigma_1)$ and $F \in L^p(\Sigma)$ is its trivial extension to $\Sigma$. Then $\exists G\in L^p(\Sigma)$ such that $(1-C_{\Sigma,w})G = F$. Since $w$ is supported on $\Sigma_1$, we can write
$$G = F + C_{\Sigma,w} G = F + C_{\Sigma_1,w} (G|_{\Sigma_1})$$
Letting $g = G|_{\Sigma_1} \in L^p(\Sigma_1)$, we have
\begin{eqnarray*}
(1-C_{\Sigma_1,w}) g &=& f
\end{eqnarray*}
so $1-C_{\Sigma_1,w}$ is invertible on $L^p(\Sigma_1)$. Furthermore, we can easily show $\|(1-C_{\Sigma_1,w})^{-1}\| \lesssim \|(1-C_{\Sigma,w})^{-1}\|$, by noticing $\|g\|_{L^p(\Sigma_1)} \leq \|G\|_{L^p(\Sigma)}$, $\|F\|_{L^p(\Sigma)} =\|f\|_{L^p(\Sigma_1)}$, and
$$\|G\|_{L^p(\Sigma)} \leq \|(1-C_{\Sigma,w})^{-1}\|\cdot \|F\|_{L^p(\Sigma)}$$

($\Leftarrow$) Now, suppose that $1-C_{\Sigma_1,w}$ is invertible on $L^p(\Sigma_1)$. Let $F \in L^p(\Sigma)$ and $f=F|_{\Sigma_1}$. Let $g:= (1-C_{\Sigma_1,w})^{-1}f$. If $(1-C_{\Sigma,w})G = F$ then similar computations as before imply that
\begin{eqnarray*}
  G(z) = \begin{cases}
              g(z), & \text{ if $z \in \Sigma_1$;}\\
              F(z) + \big(C_{\Sigma_1,w}g\big)(z), &  \text{ if $z \in \Sigma_2$.}
         \end{cases}
\end{eqnarray*}
Thus $G$ exists and is unique, furthermore using $\|w\|_\infty \lesssim 1$ we have:
$$\|G\|_{L^p(\Sigma)} \lesssim \|g\|_{L^p(\Sigma_1)} + \|F\|_{L^p(\Sigma_2)} \lesssim \Big(1+\|(1-C_{\Sigma_1,w})^{-1}\|\Big) \|F\|_{L^p(\Sigma)}$$
as desired.
\endproof

\noindent {\it Remarks:} 1. The above proof reveals that if both $1-C_{\Sigma_1,w}$ and $1-C_{\Sigma_1,w}$ are invertible and if two functions $f\in L^p(\Sigma_1)$ and $F\in L^p(\Sigma)$ agree on $\Sigma_1$, then
$$\Big((1-C_{\Sigma_1,w})^{-1}f\Big)(z) = \Big((1-C_{\Sigma,w})^{-1}F\Big)(z) \; \text{ for any $z\in\Sigma_1$}$$
(note that the left-hand side is in $L^p(\Sigma_1)$ while the right-hand side is in $L^p(\Sigma)$, however the above equality is about values on $\Sigma_1$). Consequently,
$$\Big((1-C_{\Sigma_1,w})^{-1} C_{\Sigma_1,w}I \Big)(z) = \Big((1-C_{\Sigma,w})^{-1} C_{\Sigma,w}I \Big)(z) \; \text{ for any $z\in\Sigma_1$}$$
Thus, if the weights are in $L^2$, the new $\mu$ will agree with the old $\mu$ on the old contour, provided that the new contour inherits the orientation of the old contour. As a corollary, we have

\begin{corollary} The trivial extension of $w$ from $\mathbb R$ to $\Gamma$ doesn't affect the a priori estimate (\ref{pwmujest}) (which is an estimate on $\mathbb R$).
\end{corollary}

We'll take our common contour to be $\Gamma$.

\subsubsection{Factorization of $j_\Gamma$} Below, we'll determine a natural factorization for $j_\Gamma$. Let $j=(I-w^-)^{-1}(I+w^+)$ be the jump matrix associated with $w$. Inspired by the second scheme, our factorization of $j_\Gamma \equiv (I+\Phi_-)j (I+\Phi_+)^{-1}$ will be chosen such that
$$I+\Phi_+ = (I+ w^+_\Gamma)^{-1}(I+ w^+)$$
$$I+\Phi_- = (I- w^-_\Gamma)^{-1}(I- w^-).$$
This will allow us to show the $\Gamma$ analogue of (\ref{factorize1} $\to$ \ref{factorize3}) (using completeness of $\Gamma$), and hence prove the norm equivalence claim. Detailed computation give:
\begin{eqnarray*}
 (w^-_\Gamma, w^+_\Gamma) &=&
  \begin{cases}
     (0,0), & \text{ on $\Gamma_0\cup\Gamma_3$;}\\
    (0,(I+ w^+)^{-1}-I) \equiv (0,-w^+), & \text{ on $\Gamma_1 \cup \Gamma_2$;}\\
     (I-(I- w^-)^{-1},0) \equiv (-w^-, 0), & \text{ on $\Gamma_4 \cup \Gamma_5$.}
 \end{cases}
\end{eqnarray*}
This factorization of $j_\Gamma$ is natural, it is consistent with our intuition that $w^+_\Gamma$ should have a decaying analytic continuation to the $+$ side of $\Gamma$, while $w^-_\Gamma$ should have a decaying analytic continuation to the $-$ side of $\Gamma$ (see the discussion near the beginning of Section~\ref{reductmodelsection2}).

Using the standard relation $\mu \equiv M_+(I+w^+)^{-1} \equiv M_-(I-w^-)^{-1}$, the above choice of factorization will ensure that conjugation from (the trivial extension to $\Gamma$) $j$ to $j_\Gamma$ leaves $\mu$ unchanged.

\begin{corollary} The contour deformation of the pre-model RHP doesn't affect the a priori estimate (\ref{pwmujest}).
\end{corollary}

\subsubsection{Contour deformation}

Now, we are ready to show that our deformation of the pre-model/model RHPs doesn't affect the invertibility of the respective resolvent operator:

\begin{proposition} [Contour deformation] \label{contourdeformation} Let $1<p<\infty$.

(i) For large $t$, if one of $1-C_{w_\Gamma}$ and $1-C_w$ is invertible (on the respective $L^p$ space) then so is the other. In that case, their inverse norms are ``comparable" in the sense that the following bounds are equivalent:
$$\|(1-C_{w_\Gamma})^{-1}\|_{L^p(\Gamma)\to L^p(\Gamma)} \lesssim 1$$
$$\|(1-C_{w})^{-1}\|_{L^p(\mathbb R)\to L^p(\mathbb R)}\lesssim 1$$

(ii) The same conclusions are true for the deformation of the model RHP.
\end{proposition}

\proof Note that by Proposition~\ref{extendcontour}, it suffices to show Proposition~\ref{contourdeformation} for the extension of the weights to $\Gamma$. The following argument will be essentially the same as in the proof of Lemma~\ref{boundedfactorize}.

(i) We want to show that on $L^p(\Gamma)$, the invertibility of $1-C_{\Gamma,w_\Gamma}$ and $1-C_{\Gamma,w}$ are equivalent and the resolvent norms are comparable. For simplicity, we'll suppress $\Gamma$ when writing these operators. The desired claim is now an immediate consequence of the following analogue of (\ref{factorize1} $\to$ \ref{factorize3})
$$(1- C_{\Phi})\circ(1-C_{w_\Gamma}) = 1-C_{w}$$
$$(1- C_{\Phi^{-1}})\circ(1-C_{w}) = 1-C_{w_\Gamma}$$
$$(1- C_{\Phi})\circ(1- C_{\Phi^{-1}}) = (1- C_{\Phi^{-1}})\circ(1- C_{\Phi}) =  1$$
where $C_{\Phi}$ is the Beals-Coifman operators on $\Gamma$ with weights $(-\Phi_-, \Phi_+)$, and $C_{\Phi^{-1}}$ is the Beals-Coifman operators on $\Gamma$ with weights $(I-(I +  \Phi_-)^{-1},(I + \Phi_+)^{-1} - I) \equiv (\Phi_-, - \Phi_+)$ (thanks to triangularity). The proof of these identities can be done similarly, using completeness of $\Gamma$ and the following analogue of (\ref{multiplyCauchyPhi1}): For any $h\in L^p(\Gamma)$,
\begin{eqnarray}
\label{multiplyCauchyPhi2} (C_\Gamma h)(z)  \Phi(z) = C_\Gamma\Big((C^+_\Gamma h) \Phi_+ -(C^-_\Gamma h)  \Phi_-\Big)(z),\;\;\; z\not\in \Gamma
\end{eqnarray}
Using (\ref{multiplyCauchyPhi2}) and orthogonality of $C^\pm_\Gamma$, the rest of the argument is purely algebraic, similar to the proof of (\ref{factorize1} $\to$ \ref{factorize3}).

The key idea in the proof of (\ref{multiplyCauchyPhi2}) is to exploit the uniform boundedness of $\Phi$ on $\mathbb C\setminus \Gamma$ (as in the proof of Lemma~\ref{boundedfactorize}) to make up for the lack of ``Hardy space continuation'', plus a contour integration argument.

It suffices to show that for every $\lambda\in \Gamma$,
$$(C^+_\Gamma h)(\lambda) \Phi_+(\lambda) = C^+_\Gamma\Big((C^+_\Gamma h)  \Phi_+\Big)(\lambda)$$
$$(C^-_\Gamma h)(\lambda) \Phi_-(\lambda) = - C^-_\Gamma\Big((C^-_\Gamma h)  \Phi_-\Big)(\lambda)$$
We'll only show this for the $+$ sign and for $\lambda\in \Gamma_2 \cup \Gamma_3$, the other is similar.

Let $z\in \{\pi-\alpha<\text{Arg}(z)<\pi\}$. For any $\gamma_R$ positively-oriented triangular contour with two edges of length $R$ lying on $\Gamma_2$ and $\Gamma_3$ (so one vertex is $0$), we have
$$(C_\Gamma h)(z) \Phi(z) = \frac{1}{2\pi i}\int_{\gamma_R} \frac{(C_\Gamma h)(\lambda) \Phi(\lambda)}{\lambda-z} d\lambda$$
for $R$ sufficiently large. Using Carleson's measure theorem and boundedness of $\Phi$, it is not hard to see that the contribution of the other edge of $\gamma_R$ (which we'll denote by $\widetilde\gamma_R$) vanishes as $R\to\infty$
\begin{eqnarray*}
\int_{\widetilde \gamma_R} \frac{(C_\Gamma h)(\lambda) \Phi(\lambda)}{\lambda-z} d\lambda \lesssim \|\frac{1}{z-\lambda}\|_{L^{p'}_\lambda(\widetilde \gamma_R)}\|h\|_{L^p(\Gamma)} \to 0 \;\;\;\text{ as $R\to\infty$}
\end{eqnarray*}
Consequently
$$(C_\Gamma h)(z) \Phi(z) = \frac{1}{2\pi i}\int_{\Gamma_2\cup\Gamma_3} \frac{(C_\Gamma h)(\lambda) \Phi(\lambda)}{\lambda-z} d\lambda \equiv C_\Gamma\Big(1_{\Gamma_2\cup\Gamma_3}(C^+_\Gamma h) \Phi_+\Big)(z)$$
Similar argument for similar contours for other connected components of the positive part of $\mathbb C\setminus \Gamma$) (which will give $0$ (on the left-hand side) instead of $(C_\Gamma h)(z)  \Phi(z)$ because $z$ is outside these components), we obtain
$$(C_\Gamma h)(z) \Phi(z) = C_\Gamma\Big((C^+_\Gamma h) \Phi_+\Big)(z)$$
so taking nontangential limit to $\Gamma$ we get
$$(C^+_\Gamma h)(\lambda) \Phi_+(\lambda) = C^+_\Gamma\Big((C^+_\Gamma h) \Phi_+\Big)(\lambda)$$
for $\lambda\in \Gamma_2 \cup \Gamma_3$, as desired.

(ii) Similar to (i).
\endproof

\noindent {\it Remarks:} The given proof of (\ref{multiplyCauchyPhi1}) uses the fact that $\Phi_\pm$ are bounded on $\mathbb R$ and at the same time belong to some Hardy spaces of the upper/lower half planes (so that we can use (\ref{multiplycauchyop})). This argument can certainly be reapplied to show (\ref{multiplyCauchyPhi2}) in part (i), because in the {\it pre-model} case $\Phi$ decays uniformly (allowing us to write $\Phi=C_\Gamma(\Phi_+ - \Phi_-)$). However, we avoided that approach because it is not reusable for the deformation of the model weights of part (ii), where strong decay is not available for the corresponding $\Phi$.

\proof[Alternative proof using equivalence of families of RHPs]
We largely follows \cite{deiftzhouNLSyd}. As before, we only prove (i), the argument for (ii) is similar. For any pair of weights $W^-,W^+ \in L^\infty(\Gamma)$, the invertibility of $1-C_W$ on $L^p(\Gamma)$ is connected to the unique solvability of the following family of RHP:
\begin{eqnarray*}
m_+ &=& m_- (I-W^-)^{-1}(I+W^+)\\
m_\pm -f &\in& H_p(\Gamma_\pm)
\end{eqnarray*}
which is parameterized by $f\in L^p(\Gamma)$. Here $H_p(\Gamma_\pm)$ denotes the complex Hardy space $H_p$ of the positive/negative component of $\mathbb C\setminus \Gamma$. Essentially, for a fixed $f$, the above RHP can be reduced to the functional equation
$$\mu = f + C_W \mu, \; \mu \in L^p(\Gamma)$$
($m_\pm$ can be computed from $\mu$ by $m_\pm = f + C^\pm_\Gamma(\mu(W^- + W^+)) \equiv \mu(I\pm W^\pm)$.) Thus, the unique solvability of the above family is equivalent to the invertibility of $1-C_W$.

To show that the invertibility of $1-C_{w_\Gamma}$ and $1-C_{w}$ are equivalent, we'll exploit analyticity of both pairs of weights to conjugate back and forth between the respective RHP families.  Recall the relations
$$I + \Phi_+ = (I+  w^+_\Gamma)^{-1}(I+ w^+)$$
$$I + \Phi_- = (I-  w^-_\Gamma)^{-1}(I- w^-)$$
The extensions of our jump matrices are therefore connected via the conjugation:
$$j = (I+\Phi_-)^{-1} j_\Gamma (I+\Phi_+)$$
For any $f\in L^p(\Gamma)$, we'll transform the RHP normalized to $f$ in the $w_\Gamma$ family,
$$m_{+} = m_{-} j_\Gamma, \;\; m_{\pm} \in f + H_p(\Gamma_\pm),$$
to a RHP in the $w$ family by defining
$$n_{+} = m_{+} (I + \Phi_+), \;\; n_{-} = m_{-} (I + \Phi_- )$$
(the inverse transformation can be done similarly). Clearly $n_{\pm}$ satisfies the jump relation for $w$, what's left to show is the existence of some $F \in L^p(\Gamma)$ such that $n_{\pm} \in F + H_p(\Gamma_\pm)$, and {\it $F$ depends only on $f$ and the jump matrices}. The last requirement is crucial, it basically ensures that $F$ remains a {\it coefficient} (and $n$ is the {\it variable}) of the new RHP; otherwise for any two functions $f_\pm \in L^p(\Gamma)$ we can always write
\begin{eqnarray}
\label{fpm} f_\pm = \Big(C^+_\Gamma(f_-) - C^-_\Gamma(f_+) \Big) + C^\pm_\Gamma(f_+-f_-)
\end{eqnarray}

To find $F$, observe that we can write $m_{\pm} = f + C^\pm_\Gamma h$ for $h := m_{\Gamma+} -m_{\Gamma-}\in L^p(\Gamma)$ (here the completeness of $\Gamma$ is important, it allows for cancelations $C^+_\Gamma C^-_\Gamma = C^-_\Gamma C^+_\Gamma = 0$).  Now,
$$n_{+} = f (I+ \Phi_+) + C^+_\Gamma h + (C^+_\Gamma h) \Phi_+ $$
$$n_{-} = f (I+\Phi_-) + C^-_\Gamma h + (C^-_\Gamma h) \Phi_- $$
Now (\ref{multiplyCauchyPhi2}) implies that $(C^\pm_\Gamma h) \Phi_\pm  \in H^p(\Gamma_\pm)$. Applying (\ref{fpm}) to two $h$-independent functions $f \Phi_\pm$, we can take
$$F = f + C^+_\Gamma(f \Phi_-) -C^-_\Gamma(f \Phi_+) \equiv (1- C_{\Phi}) f$$
completing our transformation.

Now, as discussed above, the conjugation on $\Gamma$ from $w \to w_\Gamma$ doesn't alter $\mu$, just like the second perturbation scheme. We therefore obtain
$$(1-C_{w_\Gamma})^{-1}f = (1-C_{w})^{-1} F = (1-C_{w})^{-1} (1- C_{\Phi}) f$$
Similarly, the reverse direction gives
$$(1-C_{w})^{-1} f = (1-C_{w_\Gamma})^{-1}(1- C_{\Phi^{-1}}) f$$
and again, these identities complete the proof.
\endproof

\noindent \underline{\bf II. Reduction on $\Gamma$.}
In this section, the deformation of the pre-model weights will be denoted by $\widetilde{w}_\Gamma$ and the deformation of the model weights will be denoted by $w_\Gamma$.
\begin{eqnarray*}
 (\widetilde w^-_\Gamma, \widetilde w^+_\Gamma) =
  \begin{cases}
    (0,0), & \text{ on $\Gamma_0\cup \Gamma_3$;}\\
    (0, -\widetilde w^+_M), & \text{ on $\Gamma_1 \cup \Gamma_2$;}\\
    (-\widetilde w^-_M, 0), & \text{ on $\Gamma_4 \cup \Gamma_5$.}
 \end{cases} &&  (w^-_\Gamma, w^+_\Gamma) =
  \begin{cases}
    (0,0), & \text{ on $\Gamma_0\cup \Gamma_3$;}\\
    (0, -w^+_M), & \text{ on $\Gamma_1 \cup \Gamma_2$;}\\
    (-w^-_M, 0), & \text{ on $\Gamma_4 \cup \Gamma_5$.}
 \end{cases}
\end{eqnarray*}
For future reference, we'll explicitly compute the potentials $\widetilde u_\Gamma, \widetilde v_\Gamma$ associated with the RHP $\widetilde w_\Gamma$:
\begin{eqnarray*}
\begin{pmatrix}0& \widetilde u_\Gamma(t)\cr \widetilde v_\Gamma(t) & 0\end{pmatrix}
&=&  -\frac{1}{2\pi i}\text{Off}\int_\Gamma \widetilde \mu_\Gamma(z) (\widetilde w^+_\Gamma(z) + \widetilde w^-_\Gamma(z)) dz
\end{eqnarray*}
here $\text{Off}(A)$ stands for the off-diagonal part of any $2\times 2$ matrix $A$. Similarly,
\begin{eqnarray*}
\begin{pmatrix}0& u_\Gamma(t)\cr v_\Gamma(t) & 0\end{pmatrix}
&=& -\frac{1}{2\pi i}\text{Off}\int_\Gamma \mu_\Gamma(z) (w^+_\Gamma(z) + w^-_\Gamma(z)) dz
\end{eqnarray*}
(notice that the $\Gamma$-deformed model weights have strong decay, unlike the original model weights on $\mathbb R$.) Assuming
$$\|(1-C_{w_\Gamma})^{-1}\|_{L^2(\Gamma)\to L^2(\Gamma)} \lesssim 1,$$
we can easily show following analogue of the direct perturbation estimate (\ref{direct_est}):
\begin{eqnarray*}
&&|\int_\Gamma \widetilde \mu_\Gamma (\widetilde w^+_\Gamma + \widetilde w^-_\Gamma) - \int_\Gamma \mu_\Gamma (w^+_\Gamma + w^-_\Gamma)| \\
&\lesssim& \|\widetilde w_\Gamma - w_\Gamma\|_2 \|C_{w_\Gamma}I\|_2 + \|\widetilde \mu_\Gamma-\mu_\Gamma\|_2 \|\widetilde w_\Gamma\|_2+ \|\widetilde w_\Gamma - w_\Gamma\|_1
\end{eqnarray*}
here $\|\widetilde \mu_\Gamma-\mu_\Gamma\|_2$ can be similarly controlled by
$$\|\widetilde w_\Gamma - w_\Gamma\|_2 \ + \  \|\widetilde w_\Gamma - w_\Gamma\|_q \|\mu_\Gamma-I\|_\frac{2q}{q-2},
\text{ for any } 2< q \leq \infty$$
Note that the above norms are taken on $\Gamma$. Consequently, for the reduction involving leading asymptotics of $u(t),v(t)$ from the deformed pre-model RHP to the deformed model RHP, it suffices to show:

\begin{proposition} [Weight estimates on $\Gamma$] \label{weightestGamma} (i) If $p$ and $q$ have one $L^2$ derivative near $0$ then for $1\leq p\leq \infty$ and large $t$:
\begin{eqnarray*}
\|\widetilde w_\Gamma - w_\Gamma\|_{L^p(\Gamma)} &\lesssim& t^{-(\frac{1}{2}+\frac{1}{p})\frac{1}{k+1}}\\
\|\widetilde w_\Gamma\|_{L^p(\Gamma)} + \|w_\Gamma\|_{L^p(\Gamma)} &\lesssim& t^{-\frac{1}{p(k+1)}}
\end{eqnarray*}
(ii) If $p$ and $q$ have two $L^2$ derivatives near $0$ then the first estimate can be improved to $t^{-(1+\frac{1}{p})\frac{1}{k+1}+\epsilon}$.
\end{proposition}
\noindent {\it Remarks:} Let $p=\infty$, this proposition implies that for large $t$ the resolvent norms $\|(1-C_{w_\Gamma})^{-1}\|$ and $\|(1-C_{\widetilde w_\Gamma})^{-1}\|$ comparable on any $L^q$.

\proof (i) We'll estimate $\|\widetilde w_\Gamma - w_\Gamma\|_{L^p(\Gamma_1)}$ below. Contribution of other $\Gamma_i$ can be estimated similarly.

Let $\theta$ be the common phase of the weights on $\Gamma_1$ ($\theta = \pm \Theta$). Then $\text{Im}(\theta) \gtrsim |z|^{k+1}$ on $\Gamma_1$ thanks to decaying property of our weights. It suffices to show that on $\Gamma_1$:
$$\Big(\delta^{-2}_+ [q] - \delta^{-2}_{0+}q(0)\Big) e^{it\theta} \lesssim t^{-\frac{1}{2(k+1)}}e^{-ct\text{ Im}\theta}$$
for some absolute constant $c>0$ (taking $L^p(\Gamma_1)$ and using a change of variable, we'll get the desired estimate).

Since $|\delta^{-2}_+(z) - \delta^{-2}_{0-}(z)| \lesssim |z|^\frac{1}{2}$ on $\Gamma_1$ (by our previous study of $\delta$), and $\delta^{-2},\delta^{-2}_0$ are bounded, we have
\begin{eqnarray*}
\Big(\delta^{-2}_+ [q] - \delta^{-2}_{0+}q(0)\Big) e^{it\theta}
&=& \delta^{-2}_+ ([q] -q(0)) e^{it\theta}  + (\delta^{-2}_+ - \delta^{-2}_{0+})q(0) e^{it\theta}\\
&\lesssim& \Big(|z|^{N} + |z|^\frac{1}{2}\Big)e^{-t\text{ Im}\theta(z)}
\end{eqnarray*}
Let $u:=t^\frac{1}{k+1}|z|$, then on $\Gamma_1$ we can write $t\text{Im}\theta(z)$ as $ C(k)u^{k+1}$ for some $C(k)>0$. Consequently, we can control the above right-hand side by
\begin{eqnarray*}
&&t^{-\frac{1}{2(k+1)}}  \Big((u^{N} + u^\frac{1}{2}) e^{-\frac{1}{2}C(k)u^{k+1}}\Big) e^{-\frac{1}{2}t\text{ Im}\theta(z)}\\
&\lesssim& t^{-\frac{1}{2(k+1)}} e^{-\frac{1}{2}t\text{ Im}\theta(z)},\;\; \text{ as desired.}
\end{eqnarray*}

The estimation of $\|\widetilde w_\Gamma\|_p$ and $\|w_\Gamma\|_p$ can be done similarly. Using boundedness of $\delta_0$ and $\delta$ on $\mathbb C$, similarly it comes down to showing $$\|e^{-tC|x|^{k+1}}\|_{L^p(\mathbb R)} \lesssim t^{-\frac{1}{p(k+1)}}, \;\;\text{$C>0$,}$$
which follows by a change of variable.

\noindent (ii) Similar argument, note that we can estimate $|\delta(z) - \delta_0(z)|$ by $|z|^{1-\epsilon}$ on $\Gamma_i$.
\endproof
To sum up, regarding resolvent bounds (i.e bounds on $(1-C_w)^{-1}$) we completely reduce our localized RHP to the model RHP with weights $w^\pm_M$, while for asymptotics of $u,v$ we can reduce it to the $\Gamma$-deformed model RHP $w^\pm_\Gamma$.

\noindent \underline{\it Reduction of the a priori estimate (\ref{pwmujest}) on $\Gamma$.} Notice that our deformation intuitively pushes every $x\in \mathbb R$ with $|x|\gtrsim 1$ away from the supports of the deformed weights, which lie in $\Gamma\setminus \mathbb R$. This separation allows us to reapply some ideas from the argument in Section~\ref{aprioriphasesect}.

The situation is indeed simpler in our case. Observe that the deformed (model/pre-model) weights vanishes on $\mathbb R$, so for any $x\in\mathbb R$ s.t. $|x|\gtrsim 1$ we can write
$$\mu_\Gamma(x) = I + C_{\Gamma\setminus \mathbb R}(\mu_\Gamma(w^+_\Gamma + w^-_\Gamma))(x)$$
$$\widetilde \mu_\Gamma(x) = I + C_{\Gamma\setminus \mathbb R}(\widetilde \mu_\Gamma(\widetilde w^+_\Gamma + \widetilde w^-_\Gamma))(x)$$

Assuming boundedness of $\|(1-C_{w_\Gamma})^{-1}\|_{L^2(\Gamma)\to L^2(\Gamma)}$ (which will be shown in the next section), Proposition~\ref{weightestGamma} easily implies
$$\|\mu_\Gamma-I\|_{L^2(\Gamma)} \lesssim t^{-\frac{1}{2(k+1)}}$$
Now, for $x\in \mathbb R$ such that $|x|\gtrsim 1$, we can write $\widetilde \mu_\Gamma(x) -\mu_\Gamma(x) $ as
\begin{eqnarray*}
&& C_{\Gamma\setminus \mathbb R}(\widetilde \mu_\Gamma \widetilde w_\Gamma)(x) - C_{\Gamma\setminus \mathbb R}(\mu_\Gamma w_\Gamma)(x)\\
&=& C_{\Gamma\setminus \mathbb R}((\widetilde \mu_\Gamma - \mu_\Gamma) \widetilde w_\Gamma)(x) + C_{\Gamma\setminus \mathbb R}((\mu_\Gamma-I)(\widetilde w_\Gamma - w_\Gamma))(x) + C_{\Gamma\setminus \mathbb R}(\widetilde w_\Gamma - w_\Gamma)(x)\\
&\lesssim& \frac{1}{|x|}\Big(\|\widetilde \mu_\Gamma - \mu_\Gamma\|_{L^2(\Gamma)} \|\widetilde w_\Gamma\|_{L^2(\Gamma)} + \|\mu_\Gamma-I\|_{L^2(\Gamma)} \|\widetilde w_\Gamma - w_\Gamma\|_{L^2(\Gamma)} + \|\widetilde w_\Gamma - w_\Gamma\|_{L^1(\Gamma)}\Big)\\
&\lesssim& t^{-\frac{3}{2(k+1)}}, \;\; \text{ as desired.}
\end{eqnarray*}

\section{Model Riemann-Hilbert problems} \label{modelsection}

In previous sections, we reduced our oscillatory RHP to the $\Gamma$-deformation of $N$ model RHPs, one for each stationary point. Consider one such model RHP. Modulo symmetry, we can assume the stationary point is $0$ of order $k$. In this section and in the following weights, $p$ and $q$ are constants, they are the values of the original $p, q$ at the given stationary point:
\begin{eqnarray*}
  (w^-_M,w^+_M) &=&
 \begin{cases}
             \Bigl(\begin{pmatrix} 0 & \delta^2_{-}pe^{-it\Theta} \cr 0 & 0 \end{pmatrix}, \begin{pmatrix} 0 & 0 \cr  \delta^{-2}_{+}q e^{it\Theta}& 0   \end{pmatrix}\Bigr), & \text{if $x \in D_{+}$;}\\
             \Bigl(\begin{pmatrix} 0 &  0\cr \delta^{-2}_{-}\frac{q}{1+pq} e^{it\Theta}& 0 \end{pmatrix}, \begin{pmatrix} 0 & \delta^2_{+} \frac{p}{1+pq}e^{-it\Theta}\cr 0 & 0   \end{pmatrix}\Bigr), & \text{if $x \in D_{-}$}
    \end{cases}
\end{eqnarray*}
Here $\Theta$ is of the form $a + b x^{k+1}$ with $a,b\in\mathbb R$ and
$$D_{+} = \{x\in\mathbb R: \Theta'(x) >0\}, \;\;\;\; D_{-} = \{x\in\mathbb R: \Theta'(x)<0\},$$
$\delta(z) := \exp(i\omega + \beta(z))$ is defined on $\mathbb {C}\setminus \mathbb {R}$ with boundary values $\delta_\pm(x)$ on $\mathbb R$. Here,
\begin{eqnarray*}
\beta(\lambda) &=& \begin{cases}
0, &\text{if $0$ is an exterior point of $\overline{D_-}$;}\\
i\epsilon \nu\ln\big[\epsilon \lambda\big], &\text{if $0$ is an endpoint of $\overline{D_-}$;}\\
-\pi\nu\text{sgn}[\text{Im}(\lambda)], &\text{if $0$ is an interior point of $\overline{D_-}$.}
\end{cases}\\
\epsilon &=& \begin{cases}0, & \text{ if $k$ is even;}\\
                                     \text{sgn }(b), & \text{ if $k$ is odd.}
                         \end{cases}\\
\nu &=& -\frac{1}{2\pi} \ln (1+pq)
\end{eqnarray*}
and $\omega$ is a real number, which we can assumed $0$ by incorporating it into $p,q$. The deformed weights were naturally chosen to be
\begin{eqnarray*}
 (w^-_\Gamma(z), w^+_\Gamma(z)) &=&
  \begin{cases}
    (0,\; 0), & \text{ $z\in \Gamma_0 \cup \Gamma_3$;}\\
    (0,\; \text{analytic cont. of } -w^+_{M}1_{\mathbb R_+}), & \text{ $z\in \Gamma_1$;}\\
    (0,\; \text{analytic cont. of } -w^+_{M}1_{\mathbb R_-}), & \text{ $z\in \Gamma_2$;}\\
    (\text{analytic cont. of } w^-_{M}1_{\mathbb R_-},\; 0), & \text{ $z\in \Gamma_4$.}\\
    (\text{analytic cont. of } w^-_{M}1_{\mathbb R_+},\; 0), & \text{ $z\in \Gamma_5$.}
 \end{cases}
\end{eqnarray*}
These can be made explicit by considering the parity of $k$ and the sign of $b$.

Our task is to show the boundedness of $(1-C_{w_\Gamma})^{-1}$ in $L^2(\Gamma)$, the a priori estimate for $\mu_\Gamma:= I + (1-C_{w_\Gamma})^{-1}C_{w_\Gamma}I$ on a set $P\subset \mathbb R$ bounded away from $\lambda_j$, and verify the desired asymptotics for those potentials recovered from $M_\Gamma$. Recall that the boundedness of $(1-C_{w_\Gamma})^{-1}$ on $L^2(\Gamma)$ is equivalent to the boundedness of $(1-C_{w_M})^{-1}$ on $L^2(\mathbb R)$ as proved in the last section (Proposition~\ref{contourdeformation}).

Recall that $1+pq >0$. We'll divide our analysis into the following cases:

1. Degenerate (or linear) case, when $pq=0$.

2. Defocusing case, when $-1 < pq<0$ (which corresponds to a defocusing system \cite{zhouL2Sobolevyd}: $q=-\overline{p}$).

3. Focusing case, when $pq>0$.

\subsection{Degenerate and defocusing cases}

 It is probably not too hard to check the boundedness of $(1-C_{w_M})^{-1}$ in the degenerate case $pq=0$, however we'll prove the unique solvability of our model RHP in these cases using the following proposition:

\begin{proposition}[Resolvent bound] \label{resolventbound} If $-1 < p q < 1$ then
$$\|(1-C_{w_M})^{-1}\|_{L^2(\mathbb R) \to L^2(\mathbb R)} \lesssim_{p,q} 1$$
\end{proposition}
\proof
We call a pair $(p,q)$ ``good" if the above conclusion is true. Then the proposition is a consequence of the following two lemmas:
\begin{lemma} If $(p,q)$ is good then so is $(cp,c^{-1}q)$ for any $c > 0$.
\end{lemma}

\begin{lemma} If $|p|,|q|<1$ then $(p,q)$ is good.
\end{lemma}

First, it is not hard to see that $(w^-_M,w^+_M)$ can be written as:
\begin{eqnarray}
\label{wMjsimpler} \begin{cases}
             \Bigl(\begin{pmatrix} 0 & p e^{-A(t,x)} \cr 0 & 0 \end{pmatrix}, \begin{pmatrix} 0 & 0 \cr  q e^{A(t,x)} & 0   \end{pmatrix}\Bigr), & \text{if $x \in D_{+}$;}\\
             \Bigl(\begin{pmatrix} 0 &  0\cr q e^{A(t,x)} & 0 \end{pmatrix}, \begin{pmatrix} 0 & p e^{-A(t,x)} \cr 0 & 0   \end{pmatrix}\Bigr), & \text{if $x \in D_{-}$}
    \end{cases}
\end{eqnarray}
where $A(t,x) = i\Big(t \Theta(x)  - 2\epsilon\nu\ln|x|\Big)$. Note that $A$ is purely imaginary on the real line, and depends on the product $pq$ rather than individual $p,q$, so it remains invariant under the change $(p,q) \to (cp,c^{-1}q)$.

The first lemma now follows from three simple observations:

(i) $\forall c> 0$, the equality $f=(1-C_{w_M})^{-1}F$ is equivalent to
$$c^{\sigma_3/2}fc^{-\sigma_3/2} = \Big(1-C_{c^{\sigma_3/2}w_Mc^{-\sigma_3/2}}\Big)^{-1}\Big(c^{\sigma_3/2}Fc^{-\sigma_3/2}\Big),\;\;\; \sigma_3 = \begin{pmatrix}1 & 0 \cr 0 & -1\end{pmatrix}.$$

(ii) The pair of weights $c^{\sigma_3/2}w_Mc^{-\sigma_3/2}$ can be obtained from the pair of weight $w_M$ if we replace $(p,q)$ by $(cp,c^{-1}q)$.

(iii) For any $g\in L^p$, the $L^p$ norms of $g$ and $c^{\sigma_3/2}gc^{-\sigma_3/2}$ are comparable up to a harmless constant depending on $c$.

Below we'll show the second lemma. Observe that the model jump matrix can be rewritten as
$$J_M = \delta^{\sigma_3}_{-} \begin{pmatrix}1 +pq &  pe^{-it\Theta}\cr q e^{it\Theta} & 1\end{pmatrix} \delta^{-\sigma_3}_{+}$$
Using analyticity and uniform boundedness of $\delta^{\pm 1}$, we can use a similar argument as before (Proposition~\ref{contourdeformation}) to deform the RHP associated with $J_M$ to the RHP associated with the following jump matrix:
$$J = \begin{pmatrix}1 +pq &  p e^{-it\Theta}\cr qe^{it\Theta} & 1\end{pmatrix}$$
More precisely, we can show that $\|(1-C_{w_M})^{-1}\|_{L^2(\mathbb R) \to L^2(\mathbb R)}$ is comparable to $\|(1-C_{W_M})^{-1}\|_{L^2(\mathbb R) \to L^2(\mathbb R)}$, where the $L^\infty$ pair of weights $W_M$ is defined by
$$I + W^+_M = (I + w^+_M)\delta^{-\sigma_3}_{+},\;\;\; I - W^+_M = (I - w^-_M)\delta^{-\sigma_3}_{-}$$
We can also think of $W^\pm_M$ as being obtained from a factorization of $J$, and recall that the invertibility and the inverse norm of $(1-C_{W_M})$ doesn't depend on the particular factorization (as long as every factors and their inverses remain bounded). So it suffices to show the desired bound (on the resolvent norm) for one factorization only, which we'll take to be
$$J = \begin{pmatrix}1 &  pe^{-it\Theta}\cr 0 & 1\end{pmatrix}\begin{pmatrix}1 &  0\cr qe^{it\Theta} & 1\end{pmatrix} \equiv (I-W^-)^{-1}(I+W^+)$$
$$\text{i.e. }\;\;\; W^- = \begin{pmatrix}0 &  pe^{-it\Theta}\cr 0 & 0\end{pmatrix},\;\; W^+ = \begin{pmatrix}0 &  0\cr qe^{it\Theta} & 0\end{pmatrix}$$
Now assume that $|p|,|q|<1$. For any $L^2$ matrix-valued function $f = (f_{ij})_{i,j= 1}^2$, we can explicitly compute $C_W f$ as
\begin{eqnarray*}
C_W f = \begin{pmatrix}C_-(f_{12}q e^{it\Theta}) & C_+(f_{11}p e^{-it\Theta}) \cr C_-(f_{22}q e^{it\Theta}) & C_+(f_{21}p e^{-it\Theta})\end{pmatrix}
\end{eqnarray*}
Consequently, using the fact that $\|C_\pm\|_{L^2(\mathbb R)\to L^2(\mathbb R)} = 1$, we have
$$\|C_W f\|_2 \leq \max(|p|,|q|)\|f\|_2$$
so $1-C_W$ is invertible on $L^2(\mathbb R)$, with
$$\|(1-C_W)^{-1}\|_{L^2(\mathbb R) \to L^2(\mathbb R)} \leq \frac{1}{1-\max(|p|,|q|)}$$
as desired. We note that if $|pq|$ are sufficiently small then $L^p$ boundedness is also true.
\endproof

\noindent {\it Remarks}: When $k=1$ it is possible to show the above resolvent bounds on $L^p$, $2\leq p<\infty$ using a suitable hypergeometric function that solves a {\it timeless} model RHP equivalent to the above RHP. This is the approach taken in \cite{deiftzhouMKdVyd,deiftzhouNLSyd} and \cite{varzuginyd}.

\subsubsection{Computation of model potentials}\label{computmodeluv}
Assuming unique solvability of the model cases, we'll show the following result:
\begin{proposition} \label{modelpotentials} The model potentials recovered from $M_\Gamma$ are of the form
\begin{eqnarray*}
\begin{pmatrix}0& u_\Gamma(t)\cr v_\Gamma(t) & 0\end{pmatrix}
= \begin{pmatrix}0& pU t^{-\frac{1}{k+1}}\exp(-2\Delta)\cr qV t^{-\frac{1}{k+1}}\exp(2\Delta) & 0\end{pmatrix}
\end{eqnarray*}
\begin{eqnarray*}
\text{with } \Delta := \frac{1}{2}ita  + \frac{i\epsilon\nu\ln t}{k+1} - i\omega
\end{eqnarray*}
for absolute constants $U,V$ depending on $pq$, $b$, and $k$ such that $U = -\overline{V}$.
\end{proposition}
\proof
To compute the potentials recovered from $M_\Gamma$, we'll make the following change of variable:
$$z \mapsto \xi = t^\frac{1}{k+1}z$$
$$M_\Gamma \mapsto N_\Gamma(\xi) = \exp\big(\Delta\sigma_3\big)M_\Gamma(z)\exp\big(-\Delta\sigma_3\big)$$
here $\Delta := \frac{1}{2}ita + \frac{i\epsilon \nu \ln t}{k + 1}$ is a $t$-dependent constant. We'll abuse notation and use the same symbol $\Gamma$ for the respective contour of $\xi$. Also, $D_{\pm}$ still denote the respective sets for $\xi$.

Then $N_\Gamma$ solves the $L^2$-normalized RHP on $\Gamma$ whose weights are analytic continuation to $\Gamma$ of the pair $(b^-_M(\xi),b^+_M(\xi))$ defined as follows:
\begin{eqnarray}
\label{bMj} \begin{cases}
             \Bigl(\begin{pmatrix} 0 & p e^{-A(\xi)} \cr 0 & 0 \end{pmatrix}, \begin{pmatrix} 0 & 0 \cr  q e^{A(\xi)} & 0   \end{pmatrix}\Bigr), & \text{if $\xi \in D_{+}$;}\\
             \Bigl(\begin{pmatrix} 0 &  0\cr q e^{A(\xi)} & 0 \end{pmatrix}, \begin{pmatrix} 0 & p e^{-A(\xi)} \cr 0 & 0   \end{pmatrix}\Bigr), & \text{if $\xi \in D_{-}$}
    \end{cases}
\end{eqnarray}
where $A(\xi) = ib \xi^{k+1}- 2i\epsilon \nu \ln |\xi|$.
The analytic continuation rule is exactly as before, and the deformed jump matrix for $\xi$ is
\begin{eqnarray}
\label{bGamma} (b^-_\Gamma(\xi), b^+_\Gamma(\xi)) &=&
  \begin{cases}
    (0,\;0), & \text{ $\xi\in \Gamma_0 \cup \Gamma_3$;}\\
    (0,\text{ analytic cont. of } -b^+_M 1_{\mathbb R_+}), & \text{ $\xi\in \Gamma_1$;}\\
    (0,\text{ analytic cont. of } -b^+_M 1_{\mathbb R_-}), & \text{ $\xi\in \Gamma_2$;}\\
    (\text{analytic cont. of } -b^-_M 1_{\mathbb R_-},\;0), & \text{ $\xi\in \Gamma_4$.}\\
    (\text{analytic cont. of } -b^-_M 1_{\mathbb R_+},\;0), & \text{ $\xi\in \Gamma_5$.}
 \end{cases}
\end{eqnarray}
The analytic continuation of $A(\xi)$ from (the respective side of) $\mathbb R$ to $\Gamma$ can be computed explicitly, in particular, we'll have $\overline{A(\xi)} = -A(\overline{\xi})$ for any $\xi\in\Gamma$.

The model RHP associated with $N_\Gamma$ is still normalized in $L^2(\Gamma)$ sense with $L^2\cap L^\infty$ weights, furthermore it is {\it timeless} i.e. the weights are independent of $t$. This RHP is clearly uniquely solvable, thanks to the unique solvability of $M_\Gamma$. Indeed, the $O(1)$ bound on $\|(1-C_{b_\Gamma})^{-1}\|$ can be seen as a particular case of the same bound on $\|\|(1-C_{w_M})^{-1}\|$ when $t=1$. Consequently, we can recover the following ``potentials'' from $N_\Gamma$ by sending $\xi\to\infty$ nontangentially in $\mathbb C\setminus \Gamma$:
\begin{eqnarray}
\label{uvN} \begin{pmatrix}0& U_{N\Gamma}\cr V_{N\Gamma} & 0\end{pmatrix} \;\;=\;\; \lim_{\xi\to\infty}\xi N_\Gamma(\xi)
\;\;=\;\; -\frac{1}{2\pi i}\text{Off}\int_\Gamma \mu_{N\Gamma}(\xi) (b^+_\Gamma(\xi) + b^-_\Gamma(\xi)) d\xi \;\;\;\;
\end{eqnarray}
Note that $U_{N\Gamma},V_{N\Gamma}$ are constants depending only on $p,q,k,b$  (it is not hard to see that these constants are independent of the contour $\Gamma$, as long as the deformed angle is small enough - anyway our chosen angle depends on $k$).

Consequently
\begin{eqnarray*}
\begin{pmatrix}0& u_\Gamma(t)\cr v_\Gamma(t) & 0\end{pmatrix}
&=& \lim_{z\to\infty}z M_\Gamma(z) \\
&=& t^{-\frac{1}{k+1}}\exp\big(-\Delta \sigma_3\big) \lim_{\xi\to\infty}\xi N_\Gamma(\xi)\exp\big(\Delta \sigma_3\big)\\
&=& t^{-\frac{1}{k+1}}\exp\big(-\Delta \sigma_3\big)\begin{pmatrix}0& U_{N\Gamma}\cr V_{N\Gamma} & 0\end{pmatrix}\exp\big(\Delta \sigma_3\big)\\
&=& \begin{pmatrix}0& U_{N\Gamma}t^{-\frac{1}{k+1}}\exp(-2\Delta)\cr V_{N\Gamma}t^{-\frac{1}{k+1}}\exp(2\Delta) & 0\end{pmatrix}
\end{eqnarray*}
To finish up, we'll show that there exist constants $U,V$ such that
$$U_{N\Gamma} = p U,\;\;\; V_{N\Gamma} =  q V,$$
and $U, V$ depend only on $p q$, $k$, $b$ and are anti-complex conjugates:
$$U = -\overline{V}$$
The existence of $U,V$ is indeed natural from (\ref{uvN}) and the strict triangular structures of the weights (where $p$ is always part of the top right and $q$ is always part of the bottom left entry).

To see that $U,V$ depend on $pq$ (instead of individuals $p$ and $q$), observe that $\forall c\in\mathbb C\setminus \{0\}$, the transformation $(p,q) \mapsto (c^2 p, c^{-2}q)$ corresponds to the following changes:
$$b_\Gamma \mapsto c^{\sigma_3} b_\Gamma c^{-\sigma_3}, \;\;\;\; \mu_{N\Gamma} \mapsto c^{\sigma_3} \mu_{N\Gamma} c^{-\sigma_3}$$
which imply that the diagonal entries of $\mu_{N\Gamma}$ are unaffected, so it leaves $U,V$ invariant.

Below, we'll show that $U,V$ are anti-complex conjugates. By (\ref{uvN}),
\begin{eqnarray*}
U &=& \frac{1}{2\pi i}\int_\Gamma  \mu_{N\Gamma}^{11}e^{-A(z)} g_1(z) dz\\
V &=& \frac{1}{2\pi i}\int_\Gamma  \mu_{N\Gamma}^{22}e^{A(z)}  g_2(z) dz
\end{eqnarray*}
where $g_1$, $g_2$ take values in $\{-1,0,1\}$ and they depend on $D_{-}, D_{+}$:
\begin{itemize}
\item If $D_{-} = \mathbb R_-$ and $D_{+} = \mathbb R_+$ then $g_1 = 1_{\Gamma_2}+1_{\Gamma_5}$, and $g_2 = 1_{\Gamma_1} + 1_{\Gamma_4}$.
\item If $D_{-} = \emptyset$ and $D_{+} = \mathbb R_- \cup \mathbb R_+$ then $g_2 = 1_{\Gamma_1} + 1_{\Gamma_5}$ and $g_1 = 1_{\Gamma_2} +  1_{\Gamma_4}$
\end{itemize}
The remaining cases are similar. By considering all cases it is not hard to see that $\overline{g_1(\xi)} = g_2(\overline{\xi})$. Thus, it suffices to show that
$$\overline{\mu^{11}_{N\Gamma}(z)} = \mu^{22}_{N\Gamma}(\overline{z}),\;\;\; \forall z\in\Gamma$$

For simplicity, we'll suppress the subscript $N\Gamma$ in $\mu_{N\Gamma}$ in our computation. Since $b_\Gamma=0$ on $\Gamma_0\cup\Gamma_3$, the operator $C_{b_\Gamma}$ is technically operating on functions supported on $\Gamma^0:=\Gamma\setminus \mathbb R$.

For convenience of notation (which will be explained shortly), we'll orient $\Gamma^0$ by orienting $|\Gamma_1|,|\Gamma_4|$ outwards and $|\Gamma_2|,|\Gamma_5|$ inwards. Then,
$$C_{b_\Gamma}f = C_{b_{\Gamma^0}}f \;\;\;\; \forall f \in L^p(\Gamma)$$
here $C_{b_{\Gamma^0}}$ is the Beals-Coifman operators on $\Gamma^0$ with the following pair of weights:
\begin{eqnarray*}
 (b^-_{\Gamma^0}(\xi), b^+_{\Gamma^0}(\xi)) &=&
  \begin{cases}
    (\text{analytic cont. of } b^+_M 1_{\mathbb R_+},\;0), & \text{ $\xi\in \Gamma_1$;}\\
    (\text{analytic cont. of } b^+_M 1_{\mathbb R_-},\;0), & \text{ $\xi\in \Gamma_2$;}\\
    (\text{analytic cont. of } -b^-_M 1_{\mathbb R_-},\;0), & \text{ $\xi\in \Gamma_4$.}\\
    (\text{analytic cont. of } -b^-_M 1_{\mathbb R_+},\;0), & \text{ $\xi\in \Gamma_5$.}
 \end{cases}
\end{eqnarray*}
As can be seen, $b^+_{\Gamma^0} \equiv 0$, so it is more notationally convenient to work with $C_{b_{\Gamma^0}}$. This is the main advantage of $b_{\Gamma^0}$ over $b_\Gamma$.

Let $g:\Gamma^0 \to \{\pm 1\}$ defined by $g = 1_{\Gamma_1\cup\Gamma_2} - 1_{\Gamma_4\cup\Gamma_5}$. Unraveling the equality $\mu = I + C_{b_\Gamma}\mu = I + C_{b_{\Gamma^0}}\mu$ we get
\begin{eqnarray*}
\mu^{11} = 1 + qT_2 \mu^{12},& \mu^{12} = pT_1 \mu^{11}\\
\mu^{22} = 1 + pT_1\mu^{21}, &\mu^{21} = qT_2 \mu^{22}
\end{eqnarray*}
where $T_1,T_2$ denote the following operators
$$T_1 f = C^+_{\Gamma^0}\Big(f e^{-A} g_1 g\Big),\;\;\; T_2 f = C^+_{\Gamma^0}\Big(f e^{A} g_2 g\Big)$$
From here it is not hard to see that $\mu^{11} = 1 + pqT_2T_1 \mu^{11}$ and $\mu^{22} = 1 + pqT_1T_2 \mu^{22}$. For any $f\in L^p(\Gamma)$ ($1<p<\infty$) we denote by $\widetilde{f}$ the function $\widetilde{f}(z) = \overline{f(\overline{z})}$. We'll show that
\begin{eqnarray}
\label{communiteT1T2} \widetilde{(T_1T_2 f)}(z) = (T_2 T_1 \widetilde{f})(z)
\end{eqnarray}
Since $pq\in \mathbb R$, this equality and uniqueness of $\mu$ automatically implies $\widetilde{\mu^{11}} = \mu^{22}$, as desired. To see (\ref{communiteT1T2}), notice that our previous discussion gives $\widetilde A = -A$, while clearly $\widetilde g = - g$. Consequently, together with the previous observation $g_1 = \widetilde{g_2}$ we have
$$\widetilde{T_1 f}(z) = -T_2 \widetilde{f}(z), \text{ and }\widetilde{T_2 f}(z) = -T_1 \widetilde{f}(z)$$
for any $z\in\Gamma$, which easily imply (\ref{communiteT1T2}).
\endproof

\subsubsection{A priori estimate for $\mu_\Gamma$}\label{muGammaverify}
We'll show the following estimate:
\begin{proposition}For any $x\in\mathbb R$ such that $|x|\gtrsim 1$, we have
$$\mu_\Gamma(x,t) = I + \begin{pmatrix} 0 &\frac{u(t)}{x}\cr \frac{v(t)}{x} &0\end{pmatrix}+ O(t^{-\frac{2}{k+1}})$$
\end{proposition}
\proof We'll make the same change of variable as in the proof of Proposition~\ref{modelpotentials}.
$$z \mapsto \xi = t^\frac{1}{k+1}z$$
$$M_\Gamma \mapsto N_\Gamma(\xi) = \exp\big(\sigma_3\Delta \big)M_\Gamma(z)\exp\big(-\sigma_3\Delta \big)$$
Under this change of variable the weights will become $b^\pm_\Gamma$ which are defined by (\ref{bGamma}). As before, we use the same symbol $\Gamma$ for the respective contour of $\xi$ (which is really a dilation of $\Gamma_z$). Since $b^\pm_\Gamma = 0$ on $\mathbb R = \Gamma_0 \cup \Gamma_3$, $\forall \xi_0\in\mathbb R\setminus \{0\}$ we have
\begin{eqnarray*}
\mu_{N\Gamma}(\xi_0)
&=& I + \frac{1}{2\pi i} \int_{\Gamma\setminus \mathbb R} \frac{\mu_{N\Gamma}(\xi)b_\Gamma(\xi)}{\xi-\xi_0}d\xi \\
&=& I + \frac{1}{2\pi i} \frac{\int_{\Gamma\setminus \mathbb R} \mu_{N\Gamma}(\xi)b_\Gamma(\xi)d\xi}{-\xi_0} + \frac{1}{2\pi i} \int_{\Gamma\setminus \mathbb R} \frac{\mu_{N\Gamma}(\xi)b_\Gamma(\xi)\xi}{\xi_0(\xi-\xi_0)}d\xi\\
&=& I + \begin{pmatrix} 0 &\frac{U_{N\Gamma}}{\xi_0}\cr \frac{V_{N\Gamma}}{\xi_0} & 0\end{pmatrix} + O(\frac{1}{|\xi_0|^2})
\end{eqnarray*}
notice that, thanks to boundedness of $(1-C_{b_\Gamma})^{-1}$ and strong decay of $b_\Gamma$,
$$\|\mu_{N\Gamma}(\xi)b_\Gamma(\xi)\xi\|_{L_\xi^1(\Gamma)} \lesssim \|\mu_{N\Gamma}-I\|_{L^2(\Gamma)} \|\xi b_\Gamma\|_{L^2(\Gamma)} + \|\xi b_\Gamma\|_{L^1(\Gamma)} \lesssim 1$$
and the trivial inequality
$$\frac{1}{\xi_0(\xi-\xi_0)} \lesssim \frac{1}{|\xi_0|^2} \;\;\;\forall \xi\in\Gamma\setminus \mathbb R$$
Finally, recover $x$ from $\xi_0$ we obtain the desired estimate (notice that $\Delta$ is purely imaginary so multiplication by $e^{\pm 2\Delta}$ won't destroy our error bound).
\endproof

\subsection{Focusing cases}
As mentioned above, the quadratic (focusing) case when $k=1$ can be studied using parabolic cylinder functions, see for instance \cite{varzuginyd}. In this section, we'll discuss the cubic (focusing) case, i.e. when $k=2$ and $pq>0$.

The main idea is to look at our deformed model RHP as an inverse monodromy problem (in the sense of Birkhoff \cite{birkhoffyd}) with one irregular singularity at $\infty$ of Poincar\'e rank $k+1$ \cite{JMUyd}). Essentially, the solvability of this inverse problem can be understood by studying the location of the poles of suitable deformation equations, which exist if $k\geq 2$ \cite{FIKNyd}.

In the cubic case (i.e. $k=2$), it is classical that these deformation equations form a system of two Painlev\' e II equations, and under the {\it imaginary} assumption $pq>0$ they can be shown to be free of real poles. This approach and the corresponding result go back to Its-Novokshenov \cite{ItsNovokshenovNotesyd} (see also \cite{BoIKyd} for a more direct perspective).

To see how to formulate our deformed model RHP as an inverse monodromy problem, we first repeat the same variable change $z\mapsto \xi$ to obtain a timeless normalized RHP $(N_\Gamma,j_\Gamma)$ with jump matrix $j_\Gamma(\xi)  = (I-b^-_\Gamma(\xi))^{-1}(I+b^+_\Gamma(\xi))$, where $b_\Gamma$ are defined by (\ref{bGamma}). In other words,
\begin{eqnarray*}
   j_\Gamma
  &=&
  \begin{cases}
     I, & \text{$\xi \in \Gamma_0\cup\Gamma_3$;}\\
    \text{ analytic cont. of } I-b^+_M 1_{\mathbb R_+}, & \text{$\xi\in \Gamma_1$;}\\
    \text{ analytic cont. of } I-b^+_M 1_{\mathbb R_-}, & \text{$\xi\in \Gamma_2$;}\\
    \text{ analytic cont. of } I-b^-_M 1_{\mathbb R_-}, & \text{$\xi\in \Gamma_4$;}\\
    \text{ analytic cont. of } I-b^-_M 1_{\mathbb R_+}, & \text{$\xi\in \Gamma_5$.}
  \end{cases}
\end{eqnarray*}
with $b_M$ defined by (\ref{bMj}). For simplicity of notation, let $c=\frac{b}{2}$.

If $c>0$, using (\ref{bMj}) the jump matrix $j_\Gamma$ can be written explicitly as
\begin{eqnarray*}
  j_\Gamma
  &=&
  \begin{cases}
     I, & \text{$\xi \in \Gamma_0\cup\Gamma_3$;}\\
    \begin{pmatrix}1 & 0\cr -q e^{2ic\xi^{k+1}} & 1\end{pmatrix}, & \text{$\xi\in \Gamma_1\cup\Gamma_2$;}\\
    \begin{pmatrix}1 &  -p e^{-2ic\xi^{k+1}}\cr 0 & 1\end{pmatrix}, & \text{$\xi\in \Gamma_4\cup\Gamma_5$.}
  \end{cases}
\end{eqnarray*}
while if $c<0$ the picture if essentially reflected across $\mathbb R$ and can be handled similarly. So without loss of generality we'll assume that $c>0$. We can think of our model RHP as a search for a $2\times 2$ matrix valued function
$$\Psi(\xi):=N_\Gamma(\xi)e^{-ic\xi^{k+1}\sigma_3}$$
analytic in the six sectors formed by $\Gamma$ with asymptotics $e^{-ic\xi^{k+1}\sigma_3}$, so that we can travel from one sector to the other using connection matrix $e^{ic\xi^{k+1}\sigma_3}j_\Gamma(\xi) e^{-ic\xi^{k+1}\sigma_3}$. If we fix the travel direction to be {\it counter clockwise} (instead of from the $+$ side to the $-$ of $\Gamma$ as we usually do), then we have six connection matrices defined on $\Gamma_0,\dots, \Gamma_5$:
$$I,   \begin{pmatrix}1 & 0\cr q & 1\end{pmatrix},    \begin{pmatrix}1 & 0\cr -q  & 1\end{pmatrix}, I,
    \begin{pmatrix}1 & -p\cr 0 & 1\end{pmatrix},    \begin{pmatrix}1 & p\cr 0 & 1\end{pmatrix}$$
(observe that $\Gamma_1$ and $\Gamma_5$ are oriented inwards while $\Gamma_2$ and $\Gamma_4$ are oriented outwards). This suggests an inverse monodromy problem in the sense of Birkhoff \cite{birkhoffyd}, although some adaptations will be required. Let $\Omega_k$ denote the sector:
$$\Omega_k = \{\xi\in\mathbb C: \frac{(k-2)\pi}{3}<\text{arg}z <\frac{k\pi}{3}\}$$
for each $1\leq k\leq 7$. Consider the following set of Stokes matrices
$$S_{2l} = \begin{pmatrix}1 & s_{2l}\cr 0 & 1\end{pmatrix}, S_{2l-1} = \begin{pmatrix}1 & 0\cr s_{2l-1} & 1\end{pmatrix},\;\;\; l=1,2,3$$
satisfying the cyclic relation $S_1S_2\dots S_6 = I$ (in general the right-hand side is $e^{-2\pi i \gamma \sigma_3}$ for a monodromy exponent $\gamma$, but in our case $\gamma$ will be $0$).
\begin{theorem}[Its-Novokshenov] For all but a discrete subset of $x\in\mathbb R$, there exists uniquely a family of ($2\times 2$ matrix valued) functions $\Psi_k(.,x)$ holomorphic in $\mathbb C$, such that $\Psi_k$ has the following asymptotics inside $\Omega_k$
$$\Psi_k(\xi,x) = (I+ \frac{m_1(x)}{\xi} + \frac{m_2(x)}{\xi^2} + \dots)e^{-i(c\xi^3+x\xi)\sigma_3}$$
and $\Psi_k$ are related via $\Psi_{k+1}(\xi)  = \Psi_k(\xi)S_k$ in $\Omega_k \cap \Omega_{k+1}$ for any $1\leq k\leq 7$. Furthermore, if $u(x) := 2m^{12}_1(x)$ and $v(x) := 2m^{21}_1(x)$ then $u,v$ satisfy
$$\begin{cases}
u_{xx} = 2u^2v + \frac{4x}{3c}u\\
v_{xx} = 2v^2u + \frac{4x}{3c}v
\end{cases}$$
and the above discrete subset of $\mathbb R$ is exactly the set of real poles of $u(x)$.
\end{theorem}
\noindent {\it Remarks}: 1. The above asymptotics should be understood as follows: Given any strict subsector $\Omega'_k\subset \Omega_k$ and any sufficiently small neighborhood in $\mathbb C$ of $x$, we have
$$\Psi_k(\xi,x)e^{i(c\xi^3 +x\xi)\sigma_3} - (I+ \frac{m_1(x)}{\xi} + \dots + \frac{m_n(x)}{\xi^n}) = O(\xi^{-n-1}),\;\; n\geq 0$$
More over, this asymptotics is differentiable, i.e. similar bounds are true for the derivatives of $\Psi_k$.

2. The above system of two Painlev\' e equations is called the deformation equation for the respective isomonodromy problem.

\proof[Sketch of proof \cite{BoIKyd,ItsNovokshenovyd}] For convenience, we sketch the proof of the above theorem, more details can be found in the references. First, using a classical theorem of Sibuya \cite{sibuyayd}, we can always find a weak solution to the above inverse problem, i.e. a family $\Phi_k$ satisfying the required properties in a neighborhood of infinity. Fix $k$. Using Birkhoff-Grothendieck's theorem, we can factorize
$$\Phi_k(\xi) = T_k^{-1}(\xi)\xi^{n_k\sigma_3}\Psi_k(\xi)$$
where $n_k$ is a nonnegative integer depending on $x$, and $T_k,\Psi_k$ are holomorphic invertible functions on $\overline{\mathbb C}\setminus \{0\}$ and $\mathbb C$ respectively. By a suitable left multiplication and an appropriate normalization, we can assume that $\Psi_k$ has an asymptotics of the form
\begin{eqnarray}
\label{Psiasymp} \Psi_k(\xi,x) = (I+ \frac{m_1(x)}{\xi} + \frac{m_2(x)}{\xi^2} + \dots)e^{-i(c\xi^3+x\xi - i n_k\ln \xi)\sigma_3}
\end{eqnarray}
inside $\Omega_k$. Here, as a result of the left multiplication, the $12$-entry of $m_j$ will be zero for every $1\leq j < 2n_k$. If $2n_k > k_0+1=3$ it is possible to show that the formal monodromy exponent $n_k$ associated with the above $\Psi_k$ is exactly $0$ by using a recursive formula of Jimbo, Miwa, and Ueno \cite{JMUyd}, giving a contradiction. Consequently, $n_k=0$ or $1$.

On the other hand, using a generalized version of the Birkhoff-Grothendieck theorem (initially proved by Malgrange \cite{malgrangeyd}, elementary proof given later by Bolibruch \cite{BoIKyd}) we can show that $n_k(x)$ are the same for all but a discrete set of $x\in\mathbb R$. This means for these $x$'s the monodromy of the associated ODE are the same, hence using a compatibility condition we can obtain the following deformation system
$$\begin{cases}
u_{xx} = 2u^2v + \frac{4x}{3c}u\\
v_{xx} = 2v^2u + \frac{4x}{3c}v
\end{cases}$$
with $u(x) := 2m^{12}_1(x)$ and $v(x) := 2m^{21}_1(x)$. Using a suitable Schlesinger transformation and the Painlev\' e property of the above system, we can show that the set of $x$ where $n_k(x)=1$ are exactly the poles of the above system. Now, for those $x$ with $n_k(x)=0$, it is not hard to see that $\Psi_k$ are exactly what we need.
\endproof

To apply this theorem to our case, we first notice that by a symmetry observation as in the proof of Proposition~\ref{resolventbound}, we can assume that $p=\overline{q}$. Let $\Psi_k$'s satisfy the conclusion of the above theorem for the following values of $S_1,\dots,S_6$:
$$\begin{pmatrix} 1 & 0\cr q & 1\end{pmatrix}, \begin{pmatrix} 1 & 0\cr 0 & 1\end{pmatrix}, \begin{pmatrix} 1 & 0\cr -q & 1\end{pmatrix}, \begin{pmatrix} 1 & -p\cr 0 & 1\end{pmatrix}, \begin{pmatrix} 1 & 0\cr 0 & 1\end{pmatrix}, \begin{pmatrix} 1 & p\cr 0 & 1\end{pmatrix}$$
Recall the second Pauli matrix $\sigma_2 = \begin{pmatrix} 0 & -i\cr i & 0\end{pmatrix}$. Using $p=\overline{q}$, we observe that
$$\sigma_2 \overline{S_k}\sigma_2 =S^{-1}_{7-k}\;\;\;\; \forall\; 1\leq k\leq 6$$
Thus, if we define $\widetilde \Psi_k(\xi)$ by $\widetilde \Psi_k(\xi):= \sigma_2 \overline{\Psi_{7-k}(\overline{\xi})}\sigma_2 \;\;\forall 1\leq k\leq 6$ and $\widetilde \Psi_7 \equiv \widetilde \Psi_1$,
then $\widetilde \Psi_k$ also satisfies the conclusion of the theorem. Consequently, $\widetilde \Psi_k = \Psi_k$, hence $m_1(x) = \sigma_2 \overline{m_1(x)}\sigma_2$, which implies
$$u(x) = -\overline{v(x)}\;\;\;\text{ for any $x\in\mathbb R$}.$$
Thus, $u$ solves the following {\it imaginary} Painlev\' e II equation:
$$u_{xx} = -2|u|^2 u + \frac{4x}{3c}u$$
By looking at the Laurent series of $u$, it is clear from the above equation that $u$ doesn't have any pole on $\mathbb R$. Consequently, the above inverse monodromy problem is solvable for any real $x$, in particular for $x=0$ (which is what we are interested in). Since $S_2 = S_5 = I$, the asymptotics (\ref{Psiasymp}) is true on $\mathbb C_+$ for $\Psi_2\equiv \Psi_3$, and on $\mathbb C_-$ for $\Psi_5\equiv \Psi_6$. Now, define $N_\Gamma$ on $\mathbb C\setminus \Gamma$ by
$$N_\Gamma(\xi)e^{-ic\xi^{k+1}\sigma_3} =
\begin{cases}
\Psi_1(\xi), & \text{for $\xi \in \Gamma_{50} \cup \Gamma_{01} $;}\\
\Psi_2(\xi), & \text{for $\xi \in \Gamma_{12}$;}\\
\Psi_4(\xi), & \text{for $\xi \in \Gamma_{23} \cup \Gamma_{34}$;}\\
\Psi_5(\xi), & \text{for $\xi \in \Gamma_{45}$.}
\end{cases}
$$
It is clear that $N_\Gamma$ satisfies the required jump relation on $\Gamma$ with jump matrix $j_\Gamma$,
furthermore using strong decay of $N_\Gamma$ and contour integration, we can easily show
\begin{eqnarray}
\label{NGamma} N_\Gamma(\xi) = I + (C_\Gamma g)(\xi),\;\;\text{ for }\xi\in\mathbb C\setminus \Gamma
\end{eqnarray}
here $g(x) :=  N_{\Gamma+}(x) - N_{\Gamma-}(x)$. Clearly $g \in L^p(\Gamma)$ for any $1<p\leq\infty$. Let $\mu_{N\Gamma}$ be defined on $\Gamma$ by
$$\mu_{N\Gamma} = N_{\Gamma+}\big(I+b^+_\Gamma\big)^{-1} \equiv N_{\Gamma-}\big(I-b^-_\Gamma\big)^{-1}$$
From (\ref{NGamma}), it is not hard to see that $\mu_{N\Gamma}$ satisfies:
$$\mu_{N\Gamma} = I + C_{b_\Gamma} \mu_{N\Gamma}.$$
Furthermore, $\mu_{N\Gamma}\in I + L^p(\Gamma)$ for any $1<p\leq \infty$. This is because $\mu_{N\Gamma}$ is asymptotically $I + O(\frac{1}{|\xi|})$ and at the same time continuous on $\Gamma$. We can recover $\mu_\Gamma(x,t)$ from $\mu_{N\Gamma}(\xi)$ by
$$\mu_\Gamma(x,t) = \mu_{N\Gamma}(t^{\frac{1}{k+1}}x)$$
Notice that $\mu_\Gamma$ is invertible as a matrix (indeed its determinant is 1) and $\mu_\Gamma, \mu^{-1}_\Gamma \in L^\infty(\Gamma)$. The invertibility (and boundedness of the inverse) of $1-C_{w_\Gamma}$ on $L_x^2(\Gamma)$ now follow from the $\Gamma$ analogue of Proposition~\ref{dzconverse}, which can be summarized briefly as:

\begin{corollary} For any $1<p<\infty$, $1-C_{w_\Gamma}$ is invertible on $L_x^p(\Gamma)$ with bounded norm:
$$\|1-C_{w_\Gamma}\|_{L_x^p(\Gamma)\to L_x^p(\Gamma)} \lesssim_p 1$$
\end{corollary}
In particular, when $t=1$ we have $\|1-C_{b_\Gamma}\|_{L_\xi^p(\Gamma)\to L_\xi^p(\Gamma)} \lesssim_p 1$. The rest of the computations (i.e. computing asymptotics of $u(t)$, $v(t)$ and proving the a priori estimates for $\mu$) can be done exactly as in Section~\ref{computmodeluv} and Section~\ref{muGammaverify}.

{\bf Acknowledgement.} The author would like to thank his advisor, Christoph Thiele, for suggesting this project and for his invaluable guidance and support throughout the course of this project. The author is indebted to John Garnett for useful discussions on complex analysis. The author would like to thank the referee for valuable corrections and suggestions, and for suggesting wording for the abstract.

\end{document}